\documentclass[12pt,a4paper,twoside]{article}
\usepackage{bbm}
\usepackage{color}
\usepackage{mathtools}
\usepackage[usenames,dvipsnames,svgnames,table]{xcolor}
\usepackage[numbers]{natbib}
\usepackage{authblk}
\usepackage{amsmath,amssymb,mathrsfs,amsthm}
\usepackage{algorithm}
\usepackage{algpseudocode}
\usepackage{lipsum}
\usepackage[margin=2.5cm]{geometry}
\usepackage[colorlinks,linkcolor=black, citecolor=blue,urlcolor=blue]{hyperref}
\usepackage{caption}
\usepackage{amsthm,amsmath,amsfonts,amssymb,bbm,multirow,makecell, tabularx}
\usepackage{cleveref}
\usepackage{graphicx}
\newtheorem{theorem}{Theorem}
\newtheorem{lemma}{Lemma}

\newtheorem{corollary}{Corollary}
\newtheorem{proposition}{Proposition}
\newtheorem{example}{Example}
\newtheorem{remark}{Remark}
\newtheorem{condition}{Condition}
\newtheorem{definition}{Definition}

\newcommand{\bbo}{\mathbbm{1}}
\newcommand{\bbR}{\mathbb{R}}
\newcommand{\bbE}{\mathbb{E}}
\newcommand{\bbP}{\mathbb{P}}

\newcommand{\caP}{\mathcal{P}}
\newcommand{\caN}{\mathcal{N}}
\newcommand{\al}{\alpha}

\newcommand{\ds}{\displaystyle}
\newcommand{\sgn}{{\rm sign}}
\newcommand{\diff}{{\rm d\,}}

\newcommand{\factorone}{{~}}

\title{Minimax and adaptive estimation of general linear functionals under sparsity}
\author{Jie Xie}
\author{Dongming Huang}
\date{}
\affil{Department of Statistics and Data Science, National University of Singapore}

\begin{document}
\maketitle
\begin{abstract}
We study nonasymptotic minimax estimation of the linear functional $L(\theta)=\eta^\top \theta$ for a high-dimensional $s$-sparse mean vector with an arbitrary loading vector $\eta$. For symmetric noise with exponentially decaying tails, we derive the sharp minimax rate, explicit in $s$, $\eta$, the tail parameter, and the noise level. The proposed estimator combines plug-in estimation for coordinates with large loadings and thresholding for coordinates with small loadings, and the matching lower bound is obtained via a loading-dependent sparse prior. For unknown sparsity, we construct an $\eta$-dependent Lepski-type procedure and show that, for a broad verifiable class of loading vectors, its risk matches the oracle rate up to the optimal logarithmic factor. Explicit examples illustrate how heterogeneity in $\eta$ changes both the minimax and adaptive rates. We also extend the analysis to non-symmetric noise, hypothesis testing, and estimation with unknown noise variance, where we show that asymmetry can increase the minimax rate in certain examples of $\eta$. Among these results, the two main technical novelties are the following. First, we extend the sharp lower-bound theory beyond the Gaussian setting via a new $\chi^2$ bound for generalized Gaussian distributions. Second, for possibly non-symmetric noise, we derive new lower bounds through a worst-case asymmetric construction.
\end{abstract}

\section{Introduction}\label{sec:introduction}

We consider the model
\begin{equation}\label{eq:sequence model}
y_j = \theta_j + \sigma \xi_j, \qquad j = 1, \ldots, d,
\end{equation}
where \( \theta = (\theta_1, \ldots, \theta_d) \in \mathbb{R}^d \) is an unknown parameter vector, \( \xi_j \) are independent, centered noise random variables with unit variance, and \( \sigma>0 \) is the noise level.

We study the problem of estimating the following linear functional:
\begin{equation}\label{eq:linear functional}
L(\theta) = \eta^\top \theta = \sum_{j=1}^d \eta_j \theta_j,
\end{equation}
based on the observations $y_1,\ldots,y_d$, where \( \eta = (\eta_1, \ldots, \eta_d) \in \mathbb{R}^d \) is a fixed loading vector.
Because coordinates with zero loadings can be dropped, we assume that all $\eta_j$ ( $j=1,\ldots,d$) are nonzero and sorted such that \( |\eta_1| \ge |\eta_2| \ge \cdots \ge |\eta_d| > 0 \).

Throughout this paper, we assume that $\theta$ is $s$-sparse, i.e.,
\[
\|\theta\|_0 = \sum_{j=1}^d \mathbbm{1}\{\theta_j \neq 0\} \le s,
\]
for some integer $s \in \{1, \ldots, d\}$.
This sparsity assumption arises in many applications, including spectroscopy, astronomy, and interferometry \cite{donoho1992maximum}. In these applications, the observed signal is typically close to zero with a few rare spikes, often described as a nearly black object. Most theoretical work has focused on Gaussian or sub-Gaussian noise, although in practice other noise distributions may also be relevant.

In this paper, we denote by \( \mathcal{P}_{j} \) the distribution of the noise variable \( \xi_j \)  ($j=1,\ldots,d$)
and by \( \mathcal{P}_\xi = (\mathcal{P}_1, \ldots, \mathcal{P}_d) \) the joint noise distribution.
We denote by \( \mathbb{P}_{\theta,\mathcal{P}_\xi} \) the distribution of \( (y_1,\ldots,y_d) \) when the signal is \( \theta \) and the joint noise distribution is \( \mathcal{P}_\xi \). We also write \( \mathbb{E}_{\theta,\mathcal{P}_\xi} \) for the corresponding expectation.
The classical Gaussian sequence model corresponds to the case where the noise variables $\xi_i$ are standard Gaussian, i.e., $\caP_1 = \cdots = \caP_d = \mathcal{N}(0,1)$, in which case we write $\caP_\xi = \mathcal{N}^{\otimes}$.

Beyond the Gaussian setting, we study a broader family of $\mathcal{P}_\xi$ that retains symmetry and fast-decaying tails. Specifically, we introduce the following class.

\begin{definition}[Symmetric distributions with exponentially decaying tails]\label{def: symmetric exponential}
For some $\al, \tau > 0$, let $\mathcal{G}_{\al,\tau}$ denote the class of distributions on $\mathbb{R}$ such that for any $P \in \mathcal{G}_{\al,\tau}$ and any random variable $W \sim P$,

$W$ is symmetric around $0$, $\mathbb{E}(W^2) = 1$,
and
\[
\forall t \ge 0, \quad \mathbb{P}(|W| \ge t) \le 2 \exp\!\left\{ - 2 \left(\frac{t}{\tau}\right)^\al \right\}.
\]
\end{definition}

Distributions in $\mathcal{G}_{\al,\tau}$ exhibit sub-Weibull tails and have recently attracted considerable attention \cite{comminges2021adaptive,li2023robust,kuchibhotla2022moving}.
The additional symmetry condition is important for our analysis.
The class parameter $\al$ specifies the tail behavior, and the parameter $\tau$ specifies the spread of the tail.
In particular, when $\al = 2$, the class reduces to a sub-Gaussian class; when $\al=2$ and $\tau=2$, it includes the standard Gaussian distribution.

We denote by $\mathcal{G}_{\al,\tau}^{\otimes}$ the class of (independent) product distributions on $\mathbb{R}^d$ whose marginals all belong to $\mathcal{G}_{\al,\tau}$. In other words, $\mathcal{G}_{\al,\tau}^{\otimes}=\{\otimes_{j=1}^{d}P_j: P_j \in \mathcal{G}_{\al,\tau}, ~~ j=1, \ldots, d\}$.
We assume the noise distribution $\mathcal{P}_\xi$ lies in $\mathcal{G}_{\al,\tau}^{\otimes}$ with fixed $\al$ and $\tau$, but the marginals of $\mathcal{P}_\xi$ may differ and remain unknown.

As a measure of the quality of an estimator $\hat{T}$ of the functional $L(\theta)$, we consider the maximum mean squared error over $s$-sparse vectors:
$$\sup_{\theta\in \Theta_{s}}\sup_{{\caP}_\xi\in  \mathcal{G}_{\al,\tau}^\otimes}\bbE_{\theta, {\caP}_\xi}\left(\hat{T}-L(\theta)\right)^2, $$
where $\Theta_{s}=\{\theta: \|\theta\|_0 \le s\}$.
In this paper, we propose rate-optimal estimators in a nonasymptotic minimax sense, that is, estimators $\widetilde{T}$ such that
$$\sup_{\theta\in \Theta_{s}}\sup_{{\caP}_\xi\in  \mathcal{G}_{\al,\tau}^\otimes}\bbE_{\theta, {\caP}_\xi}\left(\widetilde{T}-L(\theta)\right)^2\asymp \inf_{\hat{T}}\sup_{\theta\in \Theta_{s}}\sup_{{\caP}_\xi\in  \mathcal{G}_{\al,\tau}^\otimes}\bbE_{\theta, {\caP}_\xi}\left(\hat{T}-L(\theta)\right)^2,$$
where $\inf_{\hat{T}}$ denotes the infimum over all estimators.

The minimax rate for estimation under the Gaussian sequence model has been extensively studied in the literature; see, for example, \cite{ibragimov1985nonparametric,cai2004minimax,cai2005adaptive,golubev2004method,golubev2004oracle,laurent2008adaptive}.
The most closely related work to ours is \cite{collier2017minimax}, which established the nonasymptotic minimax rate for estimating the linear functional with a homogeneous loading vector (i.e., $\eta_j = 1$ for all $j$) in the Gaussian sequence model.
Specifically, for estimating the sum $\sum_{j=1}^d\theta_j$, they show that
\begin{equation}\label{eq:collier 2017}
  \inf_{\hat{T}} \sup_{\theta\in\Theta_{s}}
  \mathbb{E}_{\theta, {\cal N}^\otimes}\!\biggl(\hat{T}-\sum_{j=1}^d\theta_j\biggr)^2
  \;\asymp\; \sigma^2 s^2 \log\!\left(1 + \frac{d}{s^2}\right),
\end{equation}
and the optimal rate is attained by an estimator of the form
\begin{equation}\label{eq:collier estimator}
    \hat{T}_{s}=\begin{cases}
      \sum_{j=1}^d y_j \, \mathbbm{1}\!\left\{|y_j|>\sigma\sqrt{2\log\!\left(1+\frac{d}{s^2}\right)}\right\}, & s<\sqrt{d}, \\[6pt]
      \sum_{j=1}^d y_j, & s\ge \sqrt{d}.
    \end{cases}
  \end{equation}
This estimator includes all $y_j$ when $s$ is large and applies thresholding when $s$ is small.
However, their analysis assumes knowledge of the sparsity level $s$, which is typically unknown in practice.
To address this issue,  \cite{collier2018optimal} revisited the problem in the same setting but without assuming knowledge of $s$, and established that the optimal rate is
\begin{equation}\label{eq:collier 2018}
\sigma^2 s^2 \log\!\left(1 + \frac{d \log d}{s^2}\right).
\end{equation}

While the aforementioned theoretical works on the homogeneous loading case yield valuable insights into the roles of sparsity and dimensionality, they do not yield sharp minimax rates for estimating linear functionals with heterogeneous loadings (i.e., unequal $\eta_j$ values).
Heterogeneous loadings commonly arise in practice, as illustrated in the following two examples.

\begin{example}[Estimation of Integrals]\label{eg: integral estimation}
Let $\mathcal{H}$ be the class of square-integrable functions on $[0,1]$.
Suppose $f\in \mathcal{H}$ is unknown and $Y$ is an observation from the white noise model
$$
\diff Y(t)=f(t) \, \diff t+\sigma \, \diff W(t)
$$
where $W(t)$ is a standard Brownian motion.
It is often of interest to estimate the integral $\int fg $ for some test function $g\in \mathcal{H}$.
To see the connection to our problem, let $\{\varphi_i(t)\}$ be any orthonormal basis for $\mathcal{H}$ and define
$$
y_j = \int_0^1 \varphi_j\, \diff Y,
~~ \theta_j=\int_0^1 f \varphi_j,
~~\xi_j=\int_0^1 \varphi_j\, \diff W,
~~\text{ and } \eta_j=\int_0^1 \varphi_j g.
$$
It then follows that $\{\xi_i\}$ are independent standard normal r.v.s and $\int fg =\sum_j \eta_j \theta_j$.
If it is believed that $\theta_i=0$ for $i>d$ and $\{\theta_i\}_{i=1}^{d}$ has at most $s$ nonzero components, this model coincides with \eqref{eq:sequence model} and the estimation reduces to estimating $L(\theta)$ in \eqref{eq:linear functional}; see \cite{Johnstone2017GaussianEstimation} for further discussion. In many problems, the test function $g$ is arbitrary, so the induced sequence $\{\eta_j\}$ may lack any assumed structural properties.
\end{example}

\begin{example}[Prediction in linear regression]\label{eg: linear model prediction}
Given $n$ covariate vectors $x_i\in \mathbb{R}^d$, suppose we observe  $z_i=x_i^\top \beta+ \sigma \varepsilon_i$ where $\beta\in \mathbb{R}^{d}$ is the unknown parameter and $\varepsilon_i$ are i.i.d. samples from $\mathcal{N}(0,1)$.
A common task is to predict the response for a new fixed covariate vector $x_0$ by estimating $x_{0}^\top\beta$.
To see the connection to our problem, consider an orthogonal design, where the design matrix $\mathbf{X}$ can be written as $\mathbf{X}=\sum_{j \le d} b_j u_j e_j^\top$, with singular values $b_j>0$, singular vectors $u_j\in \mathbb{R}^n$, and canonical basis vectors $e_j\in \mathbb{R}^{d}$.
For any $j\le d$, define the transformed variables and parameters as
\begin{equation*}
 y_j=\sum_{i\le n} u_{j,i}z_i,
 ~~
 \theta_j=b_j \beta_j,
 ~~
 \xi_j=\sum_{i\le n} u_{j,i}\varepsilon_i,
 \text{ and }
 ~~
 \eta_{j}=b_j^{-1} x_{0,j}.
\end{equation*}
It follows that $y_j=\theta_j+\sigma \xi_j$, $\{\xi_j\}$ are i.i.d. from $\mathcal{N}(0,1)$, and $x_{0}^\top\beta=\sum_{j=1}^d \eta_j \theta_j=L(\theta)$.
The standard assumption that $\beta$ is $s$-sparse is equivalent to $\theta$ being $s$-sparse.
In this case, the model coincides with the sequence model in \eqref{eq:sequence model}, and the prediction problem reduces to estimating $L(\theta)$ with a loading vector $\eta = (b_j^{-1} x_{0,j})$ that depends on the target $x_0$ and the singular values $b_j$.
\end{example}

In both examples and many other applications, the loading vector is arbitrary and heterogeneous.
This presents a critical gap in the existing minimax theory:
the sequence model literature \cite{cai2004minimax,collier2017minimax,collier2018optimal} assumes equal loadings, while the linear model literature \cite{cai2017confidence,cai2021optimal} considers special loading structures (specifically, loadings with bounded coordinate ratios or polynomial decay rates).
Consequently, the optimal rates for general $\eta$ remain unknown. To our knowledge, the current paper is the first to address this problem and to provide a sharp nonasymptotic characterization of the minimax rate for estimating $L(\theta)$ with arbitrary loadings under symmetric sub-Weibull noise.

\textit{Contributions.}
We study estimation of $L(\theta)=\eta^{\top} \theta$ for $s$-sparse $\theta$ with an arbitrary loading vector $\eta$ under any noise distribution $\mathcal{P}_\xi \in \mathcal{G}_{\al,\tau}^{\otimes}$.
All results are nonasymptotic and explicit in $(s, d, \eta, \al, \sigma)$.

In the oracle setting where $s$ is known, we establish the minimax rate $\Phi_{o}(s; \eta)$, which is governed by an important quantity $\lambda_o$ that depends on $(s, \eta, \al)$ in \Cref{eq:lambda nonadaptive}.
Our theory not only recovers the homogeneous case, but also reveals how heterogeneity in $\eta$ changes the minimax rate.
To obtain sharp upper and lower bounds, the derivation must take into account the structure of the loading vector. However, this aspect is completely missing in existing works on homogeneous loadings \cite{cai2004minimax,collier2017minimax,collier2018optimal}. To fill this gap, we make the following technical contributions:

\begin{itemize}
  \item \textit{Upper bounds.}
  The estimator in \Cref{eq:collier estimator} for homogeneous loadings is no longer optimal for heterogeneous loadings.
  To achieve minimax optimality with general $\eta$, it is necessary to treat large and small loadings differently.
  We use $\lambda_o$ as the cutoff value for $|\eta_j|$ and propose a hybrid estimator that handles coordinates with large $|\eta_j|$ by plug-in and those with small $|\eta_j|$ by thresholding.
  The cutoff and the threshold are precisely calibrated to the loading vector $\eta$ and the noise tail parameter $\al$; these choices are critical for achieving minimax optimality.
  \item \textit{Lower bounds.}
The lower bound arguments in \cite{collier2017minimax} are not sharp for heterogeneous loadings, as they rely on a uniform prior distribution that ignores the differences in $|\eta_j|$.
Furthermore, their analysis is confined to Gaussian noise.
To address these limitations, we develop new techniques that apply to general noise distributions in $\mathcal{G}_{\al, \tau}$ and construct a non-uniform sparse prior where the nonzero probabilities and the magnitudes of coordinates are specifically designed based on $\eta$ and $\al$.
This refined construction yields sharp minimax lower bounds for general loading vectors and a broader family of noise distributions.
\end{itemize}

For unknown $s$, we investigate the price of adaptation.
Unlike the Lepski-type procedure for homogeneous loadings in  \cite{collier2018optimal}, where the threshold for each candidate $s$ is independent of $\eta$ and is not optimal, we identify a new $\eta$-dependent threshold such that the resulting adaptive estimator attains a rate that matches the minimax rate up to a logarithmic factor in $s$.
We further show that this rate is the unimprovable adaptive rate for a broad class of loading vectors satisfying a mild and verifiable condition (\Cref{cond: loading vector}).

In addition, we study non-symmetric noise, showing that asymmetry can increase the minimax rate in certain examples.
We also consider extensions to the setting with unknown noise variance $\sigma^2$ and to hypothesis testing for the linear functional.

Among the contributions above, the main new technical ingredients are concentrated in two places.
First, in the symmetric setting, we go beyond the Gaussian model and establish a sharp minimax lower bound over the product noise class $\mathcal{G}_{\al, \tau}$.
This extension relies on a new $\chi^2$ control for generalized Gaussian distributions.
Second, for possibly non-symmetric noise, we establish new minimax upper and lower bounds.
On the lower-bound side, we derive an additional asymmetric lower bound by constructing a least favorable asymmetric noise distribution together with a prior tailored to the asymmetric setting. In certain examples, this asymmetric lower bound is strictly larger than the symmetric lower bound, showing that non-symmetry can increase the minimax rate.

\textit{Organization.}
The rest of the paper is organized as follows. In \Cref{sec:nonadaptive minimax}, we study minimax estimation of the linear functional \( L(\theta) \) under symmetric sub-Weibull noise when the sparsity level \( s \) is known. \Cref{sec:adaptive minimax} considers the more realistic case where \( s \) is unknown and establishes the optimal adaptive rate for estimating \( L(\theta) \). In \Cref{sec:exponential minimax}, we extend the analysis to general (not necessarily symmetric) sub-Weibull noise with known sparsity. \Cref{sec:example} provides analytically tractable examples of the derived rates.
Extensions and future directions are discussed in \Cref{sec:discussion} and \Cref{sec:future}, respectively.

The last section contains the proofs of the lower bounds. Other proofs are provided in the Supplementary Material \cite{xiehuang2026supplement}.

\textit{Notation:}
We denote by \( \mathbb{N}^* \) the set of positive integers. For any \( k \in \mathbb{N}^* \), we write $[k]=\{1,\ldots,k\}$.
For $\al \in \bbR$, we write $\lfloor \al\rfloor$ and $\lceil\al\rceil$ for the floor and ceiling functions of $\al$, that is, the greatest integer not exceeding $\al$ and the smallest integer not less than $\al$, respectively.
For two real-valued functions \(f\) and \(g\), we write \(f * g\) for their convolution.
We use \(f \lesssim g\) (resp. \(f \gtrsim g\)) to mean that there exists a constant \(C > 0\) (resp. \(c > 0\)) such that
\(f \le C g\) (resp. \(f \ge c g\)).
We write \(f \asymp g\) if both \(f \lesssim g\) and \(g \lesssim f\) hold. We respectively denote by \( x \vee y \) and \( x \wedge y \) the maximum and minimum of the two real values \( x \) and \( y \), and we set \( x_+ = x \vee 0 \).
We use $\mathbbm{1}\{\cdot\}$ to denote an indicator function.
For any \( d \in \mathbb{N}^* \), and for any property \( P(j) \) over index \( j \in [d] \), we set \( \max\left\{ j \in [d] \mid P(j) \right\} = 0 \) if \( P(j) \) does not hold for any \( j \in [d] \).

\section{Minimax estimation with known sparsity}\label{sec:nonadaptive minimax}

In this section, we assume that the sparsity level $s$ and the variance $\sigma^2$ are known.
To present the minimax rate for estimating the linear functional, we first introduce a few definitions.

Given any real number $s\in [1,d]$, let $\beta =\beta(s;\eta, \al) \in \mathbb{R}$ be the unique solution to the equation
\begin{equation}\label{eq:solution nonadaptive}
\frac{\sum_{j=1}^d |\eta_j| \exp(-\beta / |\eta_j|^\al)}{\sqrt{\sum_{j=1}^d \eta_j^2 \exp(-\beta / |\eta_j|^\al)}} = \frac{s}{2},
\end{equation}
and define
\begin{equation}\label{eq:lambda nonadaptive}
    \lambda_{o} = \lambda_{o}(s;\eta, \al) = \beta_+^{1/\al}.
\end{equation}
We first verify that $\beta$ is well defined.
Indeed, the left-hand side of \Cref{eq:solution nonadaptive} is continuous and strictly decreasing
(see Lemma~B.1 in the supplementary material),

diverges to $+\infty$ as $\beta \to -\infty$, and converges to $0$ as $\beta \to +\infty$.
It follows that \Cref{eq:solution nonadaptive} admits a unique solution $\beta$ and we have $\beta_+=\beta\vee 0 \ge 0$.
Since the left-hand side is strictly monotone, $\beta$ can be computed numerically by a standard root-finding method such as bisection.
The parameter $\lambda_{o}$ captures the dependence of the linear functional estimation problem on the three key components: the sparsity level $s$, the loading vector $\eta$, and the noise tail parameter $\al$.
$\lambda_{o}$ plays a crucial role in characterizing the minimax rate, as it appears both in the expression of our lower bound and in the construction of our optimal estimator.

Our estimator is constructed as follows. Define
\(
j_1(s) = \max\left\{ j \in [d] : |\eta_j| \ge \lambda_{o}\right\}
\) and
\begin{equation}\label{eq:minimax nonadaptive estimator}
  \hat{L}_s
  = \sum_{j \le j_1(s)} \eta_j y_j
  + \sum_{j > j_1(s)} \eta_j y_j \,\mathbbm{1}\!\left\{\, |\eta_j y_j| > \factorone \sigma\,\tau\,\lambda_{o} \right\}.
\end{equation}

For every real $s\in[1,d]$, define
\begin{equation}\label{eq:nonadaptive rate}
\Phi_{\rm o}(s;\eta)
=
\left(\lambda_o s+\nu\right)^2,
\quad \text{ and } \quad
\nu
=
\sqrt{
\sum_{j=1}^d
\eta_j^2
\exp\left(-(\lambda_o/|\eta_j|)^\al\right)
}.
\end{equation}
For integer $s\in[d]$, the following theorem shows that
$\sigma^2\Phi_{\rm o}(s;\eta)$ is the minimax rate over the class $\Theta_s$ and that
$\hat L_s$ is minimax rate-optimal.

\begin{theorem}\label{thm: nonadaptive linear}
For any integer $s \in [d]$, let $\lambda_{o}$ and $\nu$ be defined as in \eqref{eq:lambda nonadaptive} and \eqref{eq:nonadaptive rate}, respectively.
\\
(1). [Lower bound] For $\al>0$, there exist constants $c, \tau_\al > 0$ such that for any $\tau\ge \tau_\al$,
it holds that
  \begin{equation}\label{eq:nonadaptive lower bound}
    \inf_{\hat{T}} \sup_{\theta\in\Theta_{s}}\sup_{{\caP}_\xi\in  \mathcal{G}_{\al,\tau}^\otimes}
    \mathbb{E}_{\theta, {\caP}_\xi}\!\left( \hat{T} - L(\theta) \right)^2
    \;\ge\; c \sigma^2\Phi_{\rm o}(s;\eta).
  \end{equation}
  (2). [Upper bound] For the estimator $\hat{L}_s$ defined in \eqref{eq:minimax nonadaptive estimator} and $\al,\tau>0$, there exists a constant $C>0$ such that
  \begin{equation}\label{eq:nonadaptive upper bound}
    \sup_{\theta\in\Theta_{s}}
    \sup_{{\caP}_\xi\in  \mathcal{G}_{\al,\tau}^\otimes}
    \mathbb{E}_{\theta, {\caP}_\xi}\!\left( \hat{L}_s - L(\theta) \right)^2
    \;\le\; C \sigma^2\Phi_{\rm o}(s;\eta).
  \end{equation}
\end{theorem}

\Cref{thm: nonadaptive linear} extends the homogeneous Gaussian result of \cite{collier2017minimax} to arbitrary loading vectors under symmetric exponentially decaying noise with general \(\alpha\). In the benchmark case \(\eta_j\equiv1\) and \(\alpha=2\), it reduces to the classical minimax rate \(\sigma^2 s^2 \log(1+d/s^2)\) established therein; see also \eqref{eq:collier 2017}.
Compared with this homogeneous benchmark, the main change is that the coordinates are no longer exchangeable.
The loading vector determines which coordinates are estimated by plug-in and which are thresholded.
The general tail parameter \(\alpha\) also affects the threshold and requires a different treatment in the analysis.
Consequently, the minimax risk depends on $(s,d,\eta, \al)$ together rather than only on \(s\) and \(d\).
To illustrate these insights concretely, we provide in \Cref{sec:example} explicit expressions for the minimax rate under several representative loading structures, including the classical homogeneous case.

\smallskip
\textbf{Interpretation of $\Phi_o(s ; \eta)$ and $\hat{L}_s$.}
The following proposition gives an alternative expression for \(\Phi_{\rm o}(s;\eta)\), analogous to \cite[Corollary 1]{chhor2024sparse}.

\begin{proposition}\label{prop:nonadaptive linear}
 Recall the definitions of $\lambda_{o}$, $\nu$, and $j_1(s)$.
 For any $s\in [d]$, we have
  \begin{equation*}
         \Phi_{\rm o}(s;\eta)= \left(\lambda_{o} s + \nu\right)^2
    \;\asymp\; \lambda_{o}^2 s^2 + \nu^2
    \;\asymp\; \lambda_{o}^2 s^2 + \sum_{j \le j_1(s)} \eta_j^2.
\end{equation*}
\end{proposition}

Proposition~\ref{prop:nonadaptive linear} shows that the minimax rate has two contributions. The term $\sum_{j\le j_1(s)}\eta_j^2$ corresponds to the variance of estimating the coordinates with the largest loadings, whereas the term $\lambda_o^2 s^2$ captures the squared bias induced by thresholding the remaining coordinates, where sparsity must be exploited.

This decomposition is directly reflected in the estimator $\hat L_s$ in \eqref{eq:minimax nonadaptive estimator}, which separates the coordinates according to the loading magnitude.
For coordinates with $|\eta_j|\ge \lambda_o$ (i.e., $j\le j_1(s)$), the contribution of $\theta_j$ to the target $L(\theta)=\sum_{j=1}^d \eta_j\theta_j$ is large, and imposing a positive threshold would mainly introduce additional bias.
Therefore, these coordinates are estimated by the plug-in rule $\eta_j y_j$.
For coordinates with $|\eta_j|<\lambda_o$ (i.e., $j>j_1(s)$), sparsity becomes useful and we threshold the transformed observations $\eta_j y_j$ at level $\sigma\tau\lambda_o$.
On the original scale this is equivalent to $|y_j|>\sigma\tau\lambda_o/|\eta_j|$,
so smaller loadings receive larger thresholds.
Thus, $\hat L_s$ should be viewed as a coordinate-specific thresholding rule adapted to the loading magnitude.

This dependence of the minimax estimator on $\eta$ is essential.
When $\eta_j\equiv 1$, all coordinates play the same role and we have $j_1(s)\in\{0,d\}$. In this case, $\hat{L}_s$ reduces to
the classical minimax estimator in \cite{collier2017minimax}.
A simple contrasting example is $\eta=
\eta^{(\varepsilon)}=(1,\varepsilon,\ldots,\varepsilon)$.
As $\varepsilon\downarrow 0$, the problem approaches estimation of the single coordinate $\theta_1$. In the limit, $\hat{L}_s$ reduces to the plug-in estimator $y_1$ and is minimax, since thresholding the first coordinate only adds bias and including any other coordinates is useless.

\smallskip

\textbf{Lower bound argument.}
For the lower bound, we apply Le Cam's method (also known as the ``method of two fuzzy hypotheses'' \cite{tsybakov2009introduction}), which is standard in the functional estimation literature. The idea is to construct two priors for $\theta$, say $\mu_1$ and $\mu_2$, that stochastically separate $L(\theta^{(1)})$ from $L(\theta^{(2)} )$ as much as possible for $\theta^{(1)} \sim \mu_1$ and $\theta^{(2)} \sim \mu_2$, while ensuring that the total variation distance between their induced sampling distributions of $(y_j)_{j=1}^{d}$ remains small.

The difficulty of deriving a sharp lower bound lies in maximizing the separation.

In previous works \cite{collier2017minimax,collier2018optimal}, the priors were chosen such that $\mu_1$ is the point mass at $0$, while $\mu_2$ is the uniform distribution over the set of exact $s$-sparse vectors with nonzero entries equal to a fixed constant $\rho$.

This construction satisfies the sparsity constraint and gives \(L(\theta^{(2)})=s\rho\) for \(\theta^{(2)}\sim\mu_2\), while \(L(\theta^{(1)})=0\) for \(\theta^{(1)}\sim\mu_1\). It yields a sharp lower bound in the homogeneous case, but is inadequate for heterogeneous loadings, where coordinates with larger absolute loadings should be more likely to be nonzero to maximize the separation of \(L(\theta^{(2)})\) from \(0\).

To address this, we adopt the random sparsity strategy of \cite{chhor2024sparse}, where each coordinate is independently set to be nonzero with probability $\pi_j$ and, if nonzero, takes value $\gamma_j$.
This approach offers two advantages: (i) each coordinate is treated independently, which enables computation of the $\chi^2$-divergence and thus allows us to control the total variation distance; and (ii) the probabilities $\pi_j$ and magnitudes $\gamma_j$ can be tailored according to the loading vector $\eta$.
A challenge arises in this construction because the values of $\left\|\theta^{(2)}\right\|_0$ and $L\left(\theta^{(2)}\right)$ become random under $\theta^{(2)} \sim \mu_2$, but we can establish probability inequalities to control both quantities.
Full details are provided in \Cref{sec:proof lower}.
Since the random-sparsity idea originates from \cite{chhor2024sparse}, we next clarify how our construction differs from theirs.

\smallskip
\textbf{Comparison with \cite{chhor2024sparse}.}
\Cref{thm: nonadaptive linear} may seem related to Theorem~1 of \cite{chhor2024sparse}, which investigated the minimax separation distance in the heteroscedastic Gaussian sequence model.
After transformation, their testing problem can be written using our notation as
\begin{equation}\label{eq: chhor testing}
H_0:\; \theta = 0
\quad \text{against} \quad
H_1:\; \|\theta \circ \eta \|_t \ge \varepsilon, ~ \|\theta\|_0 \le s,
\end{equation}
where $(\theta \circ \eta)_j=\theta_j\eta_j$ and $\|x\|_t = (\sum_{j=1}^d |x_j|^t)^{1/t}$ for $t \ge 1$.
When $\alpha=2$, \Cref{eq:solution nonadaptive} formally resembles Equation~(5) of \cite{chhor2024sparse} with $t=1$, although their equation is stated only for $t\ge 2$.

Despite the above formal resemblance, the two works address different statistical problems and the lower-bound arguments use different prior constructions.
Our result concerns minimax \emph{estimation} of the signed linear functional
\(
L(\theta)=\sum_{j=1}^d \eta_j\theta_j
\),
whereas \cite{chhor2024sparse} studies \emph{testing} against a sign-invariant $\ell^t$-norm alternative.
In the random-sparsity construction relevant to their Equation~(5), active coordinates are assigned symmetric random signs and the prior is calibrated to separate the alternative from the null in \(\ell^t\)-norm.
In contrast, our prior must generate a directional shift in the signed functional $L(\theta)$, and therefore cannot be chosen in the same sign-symmetric manner.
Instead, the signs must be aligned with \(\eta\) so that active coordinates contribute coherently to \(L(\theta)\).
This different objective changes the induced mixture distribution and the calibration of the inclusion probabilities and signal amplitudes.
Therefore, although our prior uses the same random-sparsity template, its construction is different from that in \cite{chhor2024sparse}.

Another difference is the noise model.
The analysis in \cite{chhor2024sparse} is Gaussian, corresponding to $\alpha=2$, whereas our lower bound is established under symmetric generalized Gaussian noise with general $\alpha$.
In this setting, the $\chi^2$-divergence of the induced mixtures cannot be controlled by the Gaussian calculation used in \cite{chhor2024sparse}.
We therefore develop Lemma~\ref{lem: chi square}, which gives a $\chi^2$-divergence bound for generalized Gaussian mixture distributions and enables control of the total variation distance in the fuzzy-hypotheses argument.
A detailed technical comparison of the lower-bound arguments is given in Appendix~E.

\begin{remark}
There is a useful closed-form expression for the minimax rate $\Phi_{\rm o}(s;\eta)$ in the bounded-sparsity regime.
If \(s\le \overline{s}\) for a fixed constant \(\overline{s}\), then
\[
\Phi_{\rm o}(s;\eta) \asymp \max_{j\in [d]} \eta_j^2 \log^{2/\al} (1+j);
\]
see Appendix D.4.
Interestingly, when $\al=2$, this expression is comparable to $\left[\mathbb E\max_{j\in[d]} | \eta_j Z_j| \right]^2$ for i.i.d. $Z_j\sim N(0,1)$, which coincides with the square of the minimax separation distance in \cite{chhor2024sparse} for testing \eqref{eq: chhor testing} with $t=\infty$.
This agreement is specific to the bounded-sparsity regime,
where the hardest alternatives are essentially one-sparse, so
the estimation error is governed by the magnitude of the largest noise, which relates to the $\ell^\infty$-norm testing problem.

\end{remark}

\section{Adaptation to unknown sparsity}\label{sec:adaptive minimax}
In the previous section, we proposed a rate-optimal estimator that requires knowledge of the sparsity level $s$.
When $s$ is unknown, the estimator is not implementable; if the supplied $s$ is misspecified, it may underperform.

In this section, we turn to the more realistic setting where the sparsity level $s$ is unknown, and we construct an \emph{adaptive estimator} that achieves the optimal rate of convergence for estimating the linear functional $L(\theta)$.

To construct the adaptive estimator, we first define a class of nonadaptive estimators $\{\hat{L}_s^*: s\in [d]\}$ as follows.
For any $s\in [1,d]$, let $\beta_*(s)$ be the solution to the following equation:
\begin{equation}\label{eq:solution linear adaptive}
  \frac{\sum_{j=1}^d|\eta_j| \exp(-\beta_*/|\eta_j|^\al)}{\sqrt{\sum_{j=1}^d \eta_j^2\exp(-\beta_*/|\eta_j|^\al)}}=\frac{s}{2\sqrt{\log (es)}},
\end{equation}
and set $\lambda_*(s)=(\beta_*(s))_+^{1/\al}$ and $j_2(s)=\max\left\{j\in [d]:|\eta_j| \ge \lambda_*(s)\right\}$.
The nonadaptive estimators are defined as
\begin{equation}\label{eq:nonadaptive estimator}
  \hat{L}_s^* = \sum_{j\le j_2(s)}\eta_jy_j+\sum_{j> j_2(s)}\eta_jy_j\bbo\left\{|\eta_jy_j|> \factorone \sigma\tau\lambda_*(s) \right\}, \text{ for } s=1,\ldots,d.
\end{equation}

\Cref{eq:solution linear adaptive} differs from \Cref{eq:solution nonadaptive} (for known sparsity) because its right-hand side involves the additional term $\sqrt{\log(es)}$ in the denominator.
This extra factor is crucial for achieving adaptivity. The existence and uniqueness of $\lambda_*(s)$ can be ensured by
Lemma~B.1 in the supplementary material.

The following result gives an upper bound on the risk of the estimator \(\hat{L}_s^*\) in \eqref{eq:nonadaptive estimator} over \(\Theta_s\).

\begin{theorem}\label{thm: adaptive estimator}
  For any $\al, \tau > 0$, there exists a constant $c_{1,0} > 0$, depending only on $\al$ and $\tau$, such that for any integer $s \in [d]$,  any $\theta\in\Theta_{s}$, and any ${\caP}_\xi \in {\cal G}_{\al,\tau}^\otimes$, it holds that
  \begin{align*}
    \bbE_{\theta, {\caP}_\xi}\!\left(\hat{L}_s^*-L(\theta)\right)^2
    \;\le\; c_{1,0}\,\sigma^2 \left[
      \sum_{j=1}^d \eta_j^2 \exp\!\left(-\left(\frac{\lambda_*(s)}{|\eta_j|}\right)^\al\right)
      + s^2 \lambda_*^2(s)
    \right].
  \end{align*}
\end{theorem}

We now introduce the adaptive estimator that does not rely on knowledge about $\|\theta\|_0$.
Motivated by \cite{collier2018optimal}, the adaptive estimator is selected from the collection of nonadaptive estimators $\{\hat{L}_s^*\}$ via a Lepski-type scheme.
Define $s_*=\max\left\{s\in [d]:\lambda_*(s)>0\right\}$ and $s_0=s_*+1$.
For any $s\in [1, d]$, define
$$\nu_*(s)=\sqrt{\log (es)\sum_{j=1}^d\eta_j^2\exp\left(-\left(\frac{\lambda_*(s)}{|\eta_j|}\right)^\al\right)} \quad\text{and}\quad \Phi_{*}(s;\eta)=\left\{\begin{array}{cc}
s^2\lambda_*^2(s)+\nu_*^2(s),&s\le s_0;\\
\Phi_*(s_0;\eta),&{\rm otherwise.}
\end{array}\right.$$
Our Lepski-type selection uses the thresholds $\omega_s=\left[\zeta\sigma^2{\Phi}_{\rm adp}(s;\eta)\right]^{1/2}$ where
\begin{equation}\label{eq:threshold adaptive}
  {\Phi}_{\rm adp}(s;\eta)=\left\{\begin{array}{cc}
    \Phi_{*}(s;\eta)\vee\Phi_*(1;\eta)\log^2(e(s\wedge s_0)), & 0<\al<2,\\
    \Phi_{*}(s;\eta),   & {\rm otherwise,}
  \end{array}\right.
\end{equation}
and \(\zeta>0\) is a constant that will be chosen large enough.
The selected index $\hat{s}$ is defined as
\begin{equation}\label{eq:selected index}
  \hat{s} = \min\left\{s\in \{1,\cdots,s_*\}:|\hat{L}_s^*-\hat{L}_{s^\prime}^* |\le \omega_{s^\prime}\text{ for all integers}\, s^\prime\in (s, d] \right\}
\end{equation}
with the convention that \(\hat{s} = s_0\) if the set in \eqref{eq:selected index} is empty.
Intuitively, $\hat{s}$ is chosen as the smallest sparsity level $s$ such that the estimate $\hat{L}_s^*$ is already ``stable,'' in the sense that increasing the assumed sparsity level does not lead to substantially different estimates.
The adaptive estimator is then defined as
\begin{equation}\label{eq:adaptive estimator}
  \tilde{L}_* = \hat{L}_{\hat{s}}^*.
\end{equation}
The following theorem establishes an upper bound on the risk of $\tilde{L}_*$.
\begin{theorem}\label{thm: adaptive upper}
For any $\al, \tau > 0$, we can choose $\zeta$ sufficiently large such that for any integer $s\in [d]$,
\begin{equation}\label{eq:adaptive upper bound}
  \sup_{\theta\in\Theta_{s}}\;\sup_{{\caP}_\xi \in {\cal G}_{\al,\tau}^\otimes}
  \bbE_{\theta, {\caP}_\xi}\!\left(\tilde{L}_* - L(\theta)\right)^2
  \;\le\; C \sigma^2 \Phi_{\rm adp}(s;\eta),
\end{equation}
for some constant $C > 0$ depending only on $(\al, \tau)$.
\end{theorem}

The definition of $\Phi_{\rm adp}(s;\eta)$ involves two terms, among which $\Phi_{*}(s;\eta)$ typically dominates, as detailed below.
$\Phi_{*}(s;\eta)$ naturally connects with the quantities $\Phi_{\rm o}(\cdot; \eta)$ introduced in \Cref{sec:nonadaptive minimax}.
Following an argument analogous to that of \Cref{prop:nonadaptive linear}, we can show
\[
\Phi_*(s;\eta)\;\asymp\; \log(es)\sum_{j \le j_2(s)} \eta_j^2 \;+\; s^2 \lambda_*^2(s),\forall s\le s_0.
\]
Since \Cref{eq:solution linear adaptive} is the same as \Cref{eq:solution nonadaptive} if $s$ is replaced by $s/\sqrt{\log(es)}$, we have $j_2(s)=j_1(s/\sqrt{\log(es)})$ and $\lambda_*(s)=\lambda_o(s/\sqrt{\log(es)})$.

It then follows that
\begin{equation}\label{eq:phi star equation}
  \Phi_*(s;\eta)\asymp\log(es)\bigg[ \sum_{j \le j_2(s)} \eta_j^2
   + \bigg(\frac{s}{\sqrt{\log(es)}}\bigg)^2 \lambda_*^2(s) \bigg]
   \asymp\log(es)\,\Phi_{\rm o}\biggl(\frac{s}{\sqrt{\log(es)}};\eta\biggr).
\end{equation}
Note that $\Phi_*(1;\eta)=\Phi_{\rm o}(1;\eta)$.
To ensure that $\Phi_{\rm adp}(s;\eta)\asymp \Phi_*(s;\eta)$, it suffices to require
\begin{equation}\label{eq:growth rate}
    \Phi_{\rm o}(1;\eta)\,\log(es)\;\lesssim\;\Phi_{\rm o}\!\left(\frac{s}{\sqrt{\log(es)}};\eta\right).
\end{equation}

This suggests that if the minimax rate $\Phi_{\rm o}(s;\eta)$ grows at least on the order of $\log(s)$, then we can conveniently equate ${\Phi}_{\rm adp}(s;\eta)$ with $\Phi_*(s;\eta)$.
In \Cref{sec:exponential minimax}, we will have
\[
\Phi_{\rm o}(1;\eta)\;\lesssim\; \sum_{j\le \lceil\log^{2/\al} d\rceil}\eta_j^2
\quad \text{and} \quad
\Phi_{\rm o}(s;\eta)\;\gtrsim\; \sum_{j \le (s^2 \wedge d)} \eta_j^2.
\]
These bounds demonstrate that the growth condition in \Cref{eq:growth rate} is mild and can be verified in some concrete examples, such as those in \Cref{sec:example}.

The following result shows that for any $\eta$,  $\Phi_{\rm adp}(s;\eta)$ is ``almost increasing'' in $s$ while the ratio $\Phi_{\rm adp}(s;\eta)/(s^{2}\log(es))$ is ``almost decreasing'' in $s$.

\begin{proposition}\label{prop:linear decrease}
The following relationships hold for any real numbers
$1\le s\le s^\prime\le s_0\wedge d$:
\[
\frac{\Phi_{*}(s;\eta)}{\log(es)}
\;\lesssim\;
\frac{\Phi_{*}(s^\prime;\eta)}{\log(es^\prime)},
\qquad
\frac{\Phi_{*}(s^\prime;\eta)}
{(s^\prime)^2\log(es^\prime)}
\;\lesssim\;
\frac{\Phi_{*}(s;\eta)}
{s^2\log(es)}.
\]
Consequently, there exists a constant $c_{2,0}>0$ such that
\[
\Phi_{\rm adp}(s;\eta)
\le
c_{2,0}\Phi_{\rm adp}(s^\prime;\eta),
\qquad
1\le s\le s^\prime\le d,
\]
and
\[
\Phi_{\rm adp}(s;\eta)
\le
c_{2,0}s^2\log(es)\Phi_{\rm o}(1;\eta),
\qquad
1\le s\le s_0\wedge d.
\]
\end{proposition}

From Proposition~\ref{prop:linear decrease}, we obtain
\begin{equation}\label{eq:adaptive-oracle-comparison}
    \Phi_*(s;\eta)\;\lesssim\; \Phi_{\rm o}(s;\eta)\,\log(e(s\wedge s_0)).
\end{equation}
Consequently, the risk of the adaptive estimator $\tilde{L}_*$ is within a $\log^2(e(s\wedge s_0))$ factor of the minimax risk when $s$ is known.
Moreover, if $\Phi_*(s;\eta)\asymp \Phi_{\rm adp}(s;\eta)$, then the adaptation cost is at most a factor of $\log(e(s\wedge s_0))$.
This is consistent with the existing literature. Recall that in the homogeneous loading case, where $s_0\asymp \sqrt{d\log d}$, \cite{collier2017minimax} and \cite{collier2018optimal} established the optimal rates for nonadaptive and adaptive estimation, given in \Cref{eq:collier 2017} and \Cref{eq:collier 2018}, respectively. The corresponding adaptation cost is
\[
\frac{\log \bigl(1+\frac{d\log d}{s^2}\bigr)}{\log \bigl(1+\frac{d}{s^2}\bigr)}.
\]
Moreover, this cost is of constant order when $s\lesssim d^b$ for some $b\le 1/2$, and of order $\log d$ when $s\gtrsim \sqrt{d\log d}$; see \cite{collier2018optimal} for further details.

According to the second part of \Cref{prop:linear decrease}, the upper bounds established for nonadaptive estimation in Theorem~\ref{thm: nonadaptive linear} and for the adaptive estimation in Theorem~\ref{thm: adaptive upper} grow with $s$ at most on the order of $s^2\log^{}(es)$. Controlling the upper bounds at this rate is important for the proof of Theorem~\ref{thm: adaptive upper}.

As discussed earlier, \cite{collier2018optimal} studied adaptive estimation under the homogeneous loading vector, i.e., $\eta_j = 1$ for all $j \in [d]$. Theorem~\ref{thm: adaptive upper} extends their results to the case of general loading vectors.
Although both \cite{collier2018optimal} and our work employ Lepski's method to construct the adaptive estimator, extending from the homogeneous case to the heterogeneous case is highly nontrivial and considerably more challenging than in the nonadaptive setting.
The main difficulty lies in the fact that with general loadings, there is no closed-form expression for the threshold $\lambda_*(s)$, and hence the analysis in \cite{collier2018optimal} cannot be directly applied.
To overcome this, we develop a new analysis approach that exploits key properties of the rate $\Phi_{\rm adp}(s;\eta)$ established in \Cref{prop:linear decrease}.

We now turn to the optimality of the upper bound $\sigma^2\Phi_{\rm adp}(s; \eta)$ for $\tilde{L}_*$ established in Theorem~\ref{thm: adaptive upper}.

The following theorem establishes a complementary lower bound: if an estimator achieves a risk sufficiently smaller than
$\Phi_{\rm adp}(s;\eta)$
over $s$-sparse vectors, then its maximal risk over $1$-sparse vectors can be significantly larger than the minimax nonadaptive rate $\Phi_{\rm o}(1;\eta)$.

\begin{theorem}\label{thm: adaptive lower}
  For any $\al>0$ and $\gamma\in (0,2)$, there exist constants \(\tau_\al,C_0,C_2>0\) such that the following holds for all $\tau\ge \tau_\al$.
  If an estimator \(\hat{T}\) satisfies that
  \begin{equation}\label{eq:adaptive lower small}
    \sup_{\theta\in \Theta_{s}}\sup_{{\caP}_\xi\in {\cal G}_{\al,\tau}^\otimes}\bbE_{\theta, {\caP}_\xi}\left(\hat{T}-L(\theta)\right)^2\le \frac{\sigma^2}{C_0}\cdot \Phi_{\rm adp}(s;\eta)\text{ for some } s\in [d],
  \end{equation}
 then its maximal risk over $\Theta_1=\{\theta: \|\theta\|_0\le 1\}$ is lower bounded as
  \begin{equation}\label{eq:adaptive-lower}
    \sup_{\|\theta\|_0\le 1}\sup_{{\caP}_\xi\in {\cal G}_{\al,\tau}^\otimes}\bbE_{\theta, {\caP}_\xi}\left(\hat{T}-L(\theta)\right)^2\ge \sigma^2 C_2\max\left\{\frac{\Phi_{\rm adp}(s;\eta)}{(s \wedge s_0)^{\gamma}},\Phi_{\rm o}(1;\eta)\right\}.
  \end{equation}
\end{theorem}
We illustrate the results of Theorem~\ref{thm: adaptive lower} under the homogeneous loading case.
In this setting, for $s = d^{\gamma^\prime}$ with some $\gamma^\prime \in (0,1/2)$, we can compute
(for more details, see \Cref{sec:example})
\[
\Phi_{\rm adp}(s;\eta) \;\asymp\; s^2 \log^{2/\al}\!\left(1 + \frac{d \log (es)}{s^2}\right)\]
and $s_0\asymp \sqrt{d\log d}$.
Consequently, Theorem~\ref{thm: adaptive lower} (with $\gamma=1$) implies that any estimator
whose risk is substantially smaller than $\Phi_{\rm adp}(s;\eta)$ for $s$-sparse vectors
must incur a maximal risk of at least $d^{\gamma^\prime}$ over $1$-sparse vectors.
Since $d^{\gamma^\prime}\gg \Phi_{\rm o}(1;\eta) \asymp \log d$ as $d \to \infty$,
the lower bound in \Cref{eq:adaptive-lower} is much larger than that of $\Phi_{\rm adp}(1;\eta)$, and is therefore larger than the risk of $\tilde{L}_*$ in the case where $\|\theta\|_0\le 1$.
Hence, such an estimator cannot be considered satisfactory.

Following \cite{tsybakov1998pointwise,collier2018optimal}, we now provide the formal definition of the adaptive rate in the asymptotic context where $d\to \infty$ and $\eta$ is a sequence of dimension-dependent loading vectors.
\begin{definition}\label{def: adaptive rate}
Suppose $d_0$ is a fixed integer and $\{\eta^{(d)}\in \mathbb{R}^d\}_{d=d_0}^{\infty}$ is a deterministic sequence of loading vectors.
For any $d\ge d_0$, we  consider $L(\theta)=\langle \theta, \eta^{(d)}\rangle$.
We call a function \( (s,d) \mapsto \psi_*(s;\eta^{(d)}) \) the \emph{adaptive rate of convergence on the scale of classes} $\left\{ \Theta_{s}\times {\cal G}_{\al,\tau}^\otimes:  s\in [d]\right\}$ if the following holds:
\begin{enumerate}
    \item There exists an estimator \(\hat{L}\) such that, for all \(d\ge d_0\),
    \begin{equation}\label{eq:adaptive_upper_bound}
        \max_{s\in [d]} \sup_{\theta\in \Theta_{s}}\sup_{{\caP}_\xi\in {\cal G}_{\al,\tau}^\otimes}\mathbb{E}_{\theta,{\caP}_\xi}\left( \hat{L} - L(\theta) \right)^2 / \psi_*(s;\eta^{(d)}) \le C,
    \end{equation}
    where \(C > 0\) is a constant.

    \item If there exists another function \(s \mapsto \psi_*^\prime(s;\eta^{(d)})\) and a constant \(C^\prime > 0\) such that, for all \(d\ge d_0\),
    \begin{equation}\label{eq:alternative_rate}
        \inf_{\hat{T}} \max_{s\in [d]}  \sup_{\theta\in \Theta_{s}} \sup_{{\caP}_\xi\in {\cal G}_{\al,\tau}^\otimes}\mathbb{E}_{\theta,{\caP}_\xi}\left( \hat{T} - L(\theta) \right)^2 / \psi_*^\prime(s;\eta^{(d)}) \le C^\prime,
    \end{equation}
    and
    \begin{equation}\label{eq:improvement_condition}
        \min_{s\in [d]} \frac{\psi_*^\prime(s;\eta^{(d)})}{\psi_*(s;\eta^{(d)})} \to 0 \quad \text{as } d \to \infty,
    \end{equation}
    then there exists (a sequence of) \(\bar{s} \in \{1, \ldots, d\}\) such that
    \begin{equation}\label{eq:loss_condition}
        \frac{\psi_*^\prime(\bar{s};\eta^{(d)})}{\psi_*(\bar{s};\eta^{(d)})} \cdot \min_{s\in [d]} \frac{\psi_*^\prime(s;\eta^{(d)})}{\psi_*(s;\eta^{(d)})} \to \infty \quad \text{as } d \to \infty.
    \end{equation}
\end{enumerate}
\end{definition}
Since $\eta^{(d)}$ is deterministic and its dependence on $d$ is clear, we drop the superscript in $\eta$ in the following discussion of adaptive rates.
In \Cref{def: adaptive rate}, the function \(\psi_*(s;\eta)\) is an adaptive rate of convergence if any local improvement over this rate for some \(s\) (cf.~\Cref{eq:improvement_condition}) necessarily incurs a substantially larger loss for at least one other sparsity level \(\bar{s}\) (cf.~\Cref{eq:loss_condition}).

Motivated by the illustration following \Cref{thm: adaptive lower}, we consider $\sigma^2 \Phi_{\rm adp}(s;\eta)$ as the candidate adaptive rate for general loading vectors \(\eta\).
In the homogeneous loading case in \cite{collier2018optimal}, both the minimax rate and the lower bound precluding any local improvement are available in closed form, which readily identifies the adaptive rate.
However, such explicit formulas are typically unavailable for general $\eta$, so it is not evident that an estimator with local improvement must incur a larger loss on $\Theta_1$.
To establish $\sigma^2 \Phi_{\rm adp}(s;\eta)$ as the adaptive rate, we consider the following condition on $\eta$.

\begin{condition}\label{cond: loading vector}
There exist constants $c>0$ and $\gamma_0\in (0,2)$, and
a sequence $s_{\rm cut}=s_{\rm cut}(d)\in[d]$ that diverges as $d\to\infty$ such that

\begin{equation}\label{eq:adaptive lower condition 1}
  \frac{\Phi_{\rm adp}(s;\eta)}{\Phi_{\rm o}(s;\eta)}\le c,\quad \forall 1\le s\le s_{\rm cut},
\end{equation}
and
\begin{equation}\label{eq:adaptive lower condition 2}
  \frac{\Phi_{\rm adp}(s;\eta)}{\Phi_{\rm o}(1;\eta)}\gtrsim (s\wedge s_0)^{\gamma_0},\quad \forall s\ge s_{\rm cut}.
\end{equation}
\end{condition}

In the above condition, \Cref{eq:adaptive lower condition 1} guarantees that the candidate adaptive rate $\sigma^2 \Phi_{\rm adp}(s;\eta)$ is of the same order as the minimax (nonadaptive) rate $\sigma^2 \Phi_{\rm o}(s;\eta)$ when the sparsity level is not too large ($s \le s_{\rm cut}$), while condition \eqref{eq:adaptive lower condition 2} essentially ensures that the ratio $\Phi_{\rm adp}(s;\eta)/\Phi_{\rm o}(1;\eta)$ grows faster than some polynomial in $s$ (with exponent $\gamma_0$) for large sparsity levels ($s \ge s_{\rm cut}$).

\Cref{cond: loading vector} only depends on the loading vector $\eta$ and is thus verifiable.
In \Cref{sec:example}, we will show that \Cref{cond: loading vector} is mild and is satisfied by many loading vectors, including the homogeneous loadings.

\begin{proposition}\label{prop:adaptive rate}
  Under \Cref{cond: loading vector},
  there exists
$\tau_{\al,\gamma_0}>0$ such that, for every
$\tau\ge\tau_{\al,\gamma_0}$,
  the adaptive rate of convergence on the scale of classes $\left\{ \Theta_{s}\times {\cal G}_{\al,\tau}^\otimes: s\in [d]\right\}$

is given by \(\sigma^2 \Phi_{\rm adp}(s; \eta)\).
\end{proposition}

The proof is given in Appendix B.4 of the Supplementary Material.

\section{General exponentially decaying noise distributions}\label{sec:exponential minimax}

In the previous sections, we have focused on symmetric noise distributions.
In this section, we show that the symmetry assumption plays a crucial role in establishing the preceding results; once this assumption is relaxed, the minimax rate may increase.

To formalize the problem, we consider the following larger class of noise distributions that are not necessarily symmetric (referred to as non-symmetric).
\begin{definition}[Distributions with exponentially decaying tails]\label{def: exponential}
  For some \( \al,\tau > 0 \), let \( \mathcal{H}_{\al,\tau} \) denote the class of distributions on \( \mathbb{R} \) such that for any \( P \in \mathcal{H}_{\al,\tau} \) and random variable \( W \sim P \), it holds that
  \[
  \mathbb{E}(W) = 0,
  \quad \mathbb{E}(W^2) = 1,
  \quad \text{and} \quad
  \forall t \ge 0, \quad
  \mathbb{P}(|W|\ge t)\le 2\exp\!\left\{-\left(\frac{t}{\tau}\right)^\al\right\}.
  \]
\end{definition}

We denote by $\mathcal{H}_{\al,\tau}^\otimes$ the class of product distributions on $\mathbb{R}^d$ whose marginal distributions belong to $\mathcal{H}_{\al,\tau}$.
In other words,
\(
\mathcal{H}_{\al,\tau}^\otimes
=
\left\{
    \bigotimes_{j=1}^{d} P_j :
    P_j \in \mathcal{H}_{\al,\tau}, \; j = 1, \ldots, d
\right\}
\).

For any $s \in [d]$, let $\lambda_{\cal H} \ge 0$ be the unique solution, if it exists, to
\begin{equation}\label{eq:solution exponential 0}
  \sum_{j\ge s^2} \exp\!\left(-\left|\frac{\lambda}{\eta_j}\right|^{\al}\right) = s.
\end{equation}
The left-hand side of \Cref{eq:solution exponential 0} is continuous and strictly decreasing in $\lambda$, tending to $d-s^2+1$ as $\lambda \to 0$ and to $0$ as $\lambda \to +\infty$.
Therefore, \Cref{eq:solution exponential 0} admits a unique solution whenever $s^2+s\le d+1$.

The following theorem provides a lower bound on the maximal risk of any estimator of $L(\theta)$ when the noise distribution belongs to $\mathcal{H}_{\al,\tau}^{\otimes}$.
\begin{theorem}\label{thm: exponential lower}
  For any $\al>0$, there exist $\tau_\al > 0$ and $c > 0$ such that for all $\tau\ge \tau_\al$ and $s\in [d]$, it holds that
  \begin{equation}\label{eq:exponential lower symmetric}
    \inf_{\hat{L}_s}\sup_{\theta\in \Theta_{s}}\sup_{\caP_\xi\in \mathcal{H}_{\al,\tau}^\otimes}\bbE_{\theta,\caP_\xi}\left(\hat{L}_s-L(\theta)\right)^2\ge c\sigma^2\Phi_{\rm o}(s;\eta).
  \end{equation}

  Furthermore, if (i) $s\ge3$ and $s^2+s\le d+1$, and (ii) there is some constant $\bar{C}$ such that
  \begin{equation}\label{eq:exponential lower condition}
    \sum_{j\ge s^2} \frac{\lambda_{\cal H}^4}{\eta_j^4}\exp\left(-2\left|\frac{\lambda_{\cal H}}{\eta_j}\right|^{\al}\right)\le \bar{C},
  \end{equation}
  where $\lambda_{\cal H}$ is defined in \Cref{eq:solution exponential 0},
  then there exists a constant $c^\prime>0$ depending on $(\al, \tau, \bar{C})$ such that
  \begin{equation}\label{eq:exponential lower asymmetric}
    \inf_{\hat{L}_s}\sup_{\theta\in \Theta_{s}}\sup_{\caP_\xi\in \mathcal{H}_{\al,\tau}^\otimes}\bbE_{\theta,\caP_\xi}\left(\hat{L}_s-L(\theta)\right)^2\ge c^\prime \sigma^2 s^2\lambda_{\cal H}^2.
  \end{equation}
\end{theorem}

Theorem~\ref{thm: exponential lower} establishes two complementary lower bounds.
The bound in \Cref{eq:exponential lower symmetric} follows directly from \Cref{thm: nonadaptive linear}, since ${\cal G}_{\al,\tau}^\otimes \subseteq {\cal H}_{\al,\tau}^\otimes$ according to \Cref{def: symmetric exponential,def: exponential}.
We refer to this as the \emph{symmetric lower bound}, as it is the minimax rate under symmetric noise distributions.

In contrast, the bound in \eqref{eq:exponential lower asymmetric} is obtained from a construction that explicitly exploits the asymmetry of the noise distribution; we therefore refer to it as the \emph{asymmetric lower bound}.
As demonstrated in \Cref{sec:example}, when the loading vector is homogeneous, the asymmetric lower bound can be strictly larger than its symmetric counterpart for $s = o(\sqrt{d})$, indicating that the results established in the symmetric setting no longer hold once the noise distribution loses symmetry.

The asymmetric lower bound is one of the main technical novelties of the paper. Unlike the symmetric case, it cannot be obtained from the generalized-Gaussian \(\chi^2\) argument. Instead, we construct a least favorable asymmetric noise distribution and develop a new lower-bound argument based on Hellinger control.
Related Hellinger-based methods were used in earlier works such as \cite{comminges2021adaptive,li2023robust}, mainly for estimating the quadratic functional \(\|\theta\|_2^2=\sum_{j=1}^d \theta_j^2\). These methods do not directly extend to the estimation of a linear functional with a heterogeneous loading vector. To handle this, we again use the random sparsity strategy and construct a new prior for the asymmetric setting. We omit the technical details here and defer the full proof to the supplementary material.

The lower bound highlights the challenge introduced by asymmetry in the noise distribution. To complement this result, we establish an upper bound that generally matches the lower bound up to a logarithmic factor and can coincide with it in certain regimes.
For any integer $s \in [d]$, define
\begin{equation}\label{eq:index unknown}
  j_3(s) = \big\lceil s^2 \log^{2/\al}(ed/s) \big\rceil \wedge d,
\end{equation}
and consider the estimator
\begin{equation}\label{eq:exponential estimator}
  \hat{L}_{\mathcal{H}}
  = \sum_{j \le j_3(s)} \eta_j y_j
  + \sum_{j > j_3(s)} \eta_j y_j \,
    \mathbbm{1}\!\left\{ |y_j| > c_{\mathcal{H}}  \sigma \log^{1/\al}\!\left(\frac{ed}{s}\right) \right\},
\end{equation}
where $c_{\mathcal{H}} > 0$ is a constant depending only on $\al$ and $\tau$.

The next theorem gives an upper bound on the risk of the estimator \( \hat{L}_{\cal H} \).
\begin{theorem}\label{thm: exponential upper}
  Let \(1 \le s \le d\). Then, for any \( \al > 0, \tau > 0 \), there exist constants $c_{\cal H}$ and $C$ depending on $\al$ and $\tau$ such that the following holds for the estimator $\hat{L}_{\cal H}$ in \Cref{eq:exponential estimator}:
  \begin{equation}\label{eq:exponential upper bound}
    \sup_{\theta\in \Theta_{s}}\sup_{\caP_\xi\in \mathcal{H}_{\al,\tau}^\otimes}\bbE_{\theta,\caP_\xi}\left(\hat{L}_{\cal H}-L(\theta)\right)^2\le C\sigma^2{\sum_{j \le j_3(s)} \eta_j^2}.
  \end{equation}
\end{theorem}
Similar to the estimator in \Cref{eq:minimax nonadaptive estimator} for symmetric noise, the estimator $\hat{L}_{\mathcal{H}}$ defined in \eqref{eq:exponential estimator} employs a loading-dependent shrinkage strategy: it uses the plug-in estimator for components with large loadings ($j \le j_3(s)$) and applies thresholding for components with small loadings ($j > j_3(s)$).
Despite this similarity, the analysis differs substantially from the symmetric noise case. The key observation is that under symmetric noise, if $\theta_j = 0$, the thresholding estimator is unbiased, i.e., $\mathbb{E}[\,y_j \mathbbm{1}\{|y_j| \ge \lambda\}\,] = 0$ for any threshold $\lambda > 0$. Under non-symmetric noise, however, this property does not hold in general. Consequently, we employ different techniques to analyze the non-symmetric case.

The upper bound established in \Cref{thm: exponential upper} is nearly sharp.
\Cref{sec:example-homogeneous} shows an example where this upper bound can match the lower bounds in \Cref{thm: exponential lower} in some regimes.
In general, we have the following proposition for comparing this upper bound with the minimax rate characterized in \Cref{thm: nonadaptive linear}.

\begin{proposition}\label{prop:exponential match}
Suppose $d\ge 2$.
For $\Phi_{\rm o}(s;\eta)$ defined in \eqref{eq:nonadaptive rate} and $j_3(s)$ defined in \eqref{eq:index unknown}, we have
$$
\Phi_{\rm o}(s;\eta)\gtrsim \sum_{j\le (s^2\wedge d)}\eta_j^2
\;\gtrsim\; \frac{1}{\log^{2/\al}d}\sum_{j\le j_3(s)}\eta_j^2.
$$
\end{proposition}

This proposition shows that the symmetric lower bound $\sigma^2 \Phi_{\mathrm{o}}(s; \eta)$ matches the upper bound in \Cref{thm: exponential upper} up to logarithmic factors, indicating that the cost of non-symmetry is at most logarithmic.
In particular, in the dense regime where $s \gtrsim \sqrt{d}$, both the symmetric lower bound and the upper bound in Theorem~\ref{thm: exponential upper} are of order $\sigma^2 \sum_{j=1}^d \eta_j^2$, yielding the exact minimax rate.
In this regime, the rate is independent of both the sparsity level $s$ and the symmetry of the noise distribution.

\section{Examples}\label{sec:example}
We now present several examples to illustrate the rates derived above, covering both nonadaptive and adaptive estimation settings. Throughout this section, we set $\sigma=1$, and  we use the notation $\Phi_{\rm o}(s;\eta)$ and $\Phi_{\rm adp}(s;\eta)$ to denote the optimal rates corresponding to
minimax estimation with known $s$ and adaptive estimation, respectively.

In addition, in the homogeneous loading vector example, we also examine the upper bound and asymmetric lower bound for the minimax rate $\Phi_{\rm ns}(s;\eta)$ under general (possibly non-symmetric) noise.

\subsection{Homogeneous loading vector}\label{sec:example-homogeneous}
Assume that $\eta_1=\cdots=\eta_d=1$. Up to multiplicative constants, the rates can be computed as follows:
\begin{equation}\label{eq:homogeneous case}
    \Phi_{\rm o}(s;\mathbf{1}_d)\asymp s^2 \log^{2/\al}\!\left(1+\frac{d^{\al/2}}{s^\al}\right),\quad \Phi_{\rm adp}(s;\mathbf{1}_d)\asymp s^2 \log^{2/\al}\!\left(1+\frac{d^{\al/2}\log^{\al/2}(e s)}{s^\al}\right)
\end{equation}
and
$$\Phi_{\rm ns}(s;\mathbf{1}_d)\asymp \left\{
\begin{array}{cc}
   \ds s^2\log^{2/\al} d ,  & \quad\ds  3\le s\lesssim \frac{\sqrt{d}}{\log^{2/\al} d};\\
   d,  & \quad s\gtrsim \sqrt{d}.
\end{array}\right.$$

\begin{itemize}
\item \textit{Adaptive rate:}
    We can verify that Condition~\ref{cond: loading vector} is satisfied (with $s_{\rm cut} \asymp d^{\gamma}$ for some $0 < \gamma < 1/2$) and $\gamma_0=1$, and thus our adaptive estimator is rate-optimal. Comparing $\Phi_{\rm o}(s;\mathbf{1}_d)$ and $\Phi_{\rm adp}(s;\mathbf{1}_d)$, the cost of adaptation is at most $\log d$, which is attained when $s \gtrsim \sqrt{d \log d}$.

    \item \textit{Comparison with existing results:}
    The homogeneous case has been studied previously under Gaussian noise, for both nonadaptive \cite{collier2017minimax} and adaptive estimation \cite{collier2018optimal}.
    For comparison, we set $\al = 2$, and consider the sub-Gaussian noise class ${\cal G}_{2,\tau}$ for any $\tau>0$;
   \Cref{eq:homogeneous case} shows that the optimal rates remain the same as under Gaussian noise.

    \item \textit{Comparison of symmetric and non-symmetric noise:}
    For general (possibly non-symmetric) noise, the \textit{asymmetric lower bound}
in \Cref{thm: exponential lower} matches the upper bound in
\Cref{thm: exponential upper} when
$3\le s\lesssim \sqrt{d}/\log^{2/\al}d$.
In the dense regime $s\gtrsim\sqrt{d}$, the symmetric lower bound and
the upper bound are both of order $d$.

    Therefore, our analysis is sharp for $s \lesssim \sqrt{d}/\log^{2/\al} d$ and $s\gtrsim \sqrt{d}$, but there remains a logarithmic gap between the upper and lower bounds in the narrow intermediate regime $\sqrt{d}/\log^{2/\al} d \lesssim s \lesssim \sqrt{d}$, which presents an interesting direction for future research. Comparing $\Phi_{\rm o}(s;\mathbf{1}_d)$ and $\Phi_{\rm ns}(s;\mathbf{1}_d)$, we observe that the absence of symmetry in the noise distribution increases the minimax rate from $s^2 \log^{2/\al}(1+d^{\al/2}/s^\al)$ to $s^2 \log^{2/\al} d$, which is at most a polylogarithmic factor larger.
\end{itemize}
\subsection{Two-phase loading vector}
We consider a loading vector \(\eta\) of the form
\begin{equation}\label{eq:two phase loading}
    \eta_j=
    \begin{cases}
        \lambda, & j \le d_0,\\
        1, & j > d_0.
    \end{cases}
\end{equation}
This loading vector \(\eta\) is divided into two parts: the first subvector consists of a small number of large components, while the second subvector consists of many small components.
This structure is related to the homogeneous case, but it accommodates unbounded ratios of coordinates, a setting less studied in the literature.
For convenience, we parameterize \(d_0\) and \(\lambda\) as \(d_0 = d^{\gamma_d}\) and \(\lambda = d^{\gamma_\lambda}\) for some constants \(\gamma_d > 0\) and \(\gamma_\lambda > 0\).
To uncover an interesting phase transition phenomenon as discussed below, we assume \(2\gamma_\lambda + \gamma_d < 1\).

The minimax rate is given by
\begin{align*}
\Phi_{\mathrm{o}}(s;\eta)
&\asymp
\begin{cases}
d^{2\gamma_\lambda}s^{2}\log^{2/\al}\!\left(1+\dfrac{d^{\gamma_d\al/2}}{s^{\al}}\right),
& s \lesssim \dfrac{d^{\gamma_\lambda+\gamma_d/2}}{\log^{1/\al} d},\\[6pt]
s^{2}\log^{2/\al}\!\left(1+\dfrac{d^{\al/2}}{s^{\al}}\right),
& s \gtrsim \dfrac{d^{\gamma_\lambda+\gamma_d/2}}{\log^{1/\al} d},
\end{cases}
\end{align*}
while the adaptive rate is given by
\begin{align*}
\Phi_{\mathrm{adp}}(s;\eta)
&\asymp
s^{2}\!\left[
d^{2\gamma_\lambda}\log^{2/\al}\!\left(1+\dfrac{d^{\gamma_d\al/2}\log^{\al/2}(es)}{s^{\al}}\right)
+\log^{2/\al}\!\left(1+\dfrac{d^{\al/2}\log^{\al/2}(es)}{s^{\al}}\right)
\right].
\end{align*}
\begin{itemize}
    \item \textit{Phase transition in the minimax rate:}
    Unlike the homogeneous case, here a phase transition arises in the minimax (nonadaptive) rate $\Phi_{\rm o}$. Specifically, when $s \lesssim d^{\gamma_\lambda+\gamma_d/2}/\log^{1/\al} d$, the error of estimating the first subvector (with $\eta_j=\lambda$) dominates, and the minimax rate equals $\Phi_{\rm o}(s; \lambda \mathbf{1}_{d_0})$ in \eqref{eq:homogeneous case}, which explicitly depends on $\gamma_\lambda$ and $\gamma_d$.
    In contrast, when $s \gtrsim d^{\gamma_\lambda+\gamma_d/2}/\log^{1/\al} d$, there are sufficiently many nonzero entries in the second subvector (with $\eta_j=1$) such that the overall error is dominated by the estimation error resulting from the second subvector, and the minimax rate equals $\Phi_{\rm o}(s; \mathbf{1}_{d})$ in \eqref{eq:homogeneous case}, which is independent of $\gamma_d$ and $\gamma_\lambda$.
    Such a phase transition occurs if and only if $\gamma_d + 2\gamma_\lambda < 1$ (i.e., $\lambda^2 d_0 \ll d$); otherwise, the first subvector always dominates and the second subvector never contributes to the minimax rate.

    \item \textit{Adaptive rate:}
    We can verify that Condition~\ref{cond: loading vector} is satisfied (with $s_{\rm cut} \asymp d^{\gamma}$ for some $\gamma \in (0, \gamma_d/2)$ and $\gamma_0=1$), and thus our adaptive estimator is rate-optimal.
\end{itemize}

\subsection{Exponentially decaying loading vector}
Following \cite{chhor2024sparse}, we consider the exponentially decaying loading vector: Let $\varphi:[0,\infty)\to \mathbb{R}$ be a non-decreasing convex function such that $\varphi(0)=0$,
and let
$$
\eta_j = \exp\!\left(-\varphi(j-1)\right),
\qquad j \in \{1,\ldots,d\}.
$$
We allow the function $\varphi$ to depend on $d$; for example,  $\varphi(x)=2(x/d)^\gamma$ for $\gamma\ge 1$.

Given $\varphi$, define $j_0 = \min\{\, j \mid \eta_j < 1/2 \,\}$
if the set is nonempty; otherwise set $j_0 = d+1$.
The optimal rates are given by
$$\Phi_{\rm o}(s;\eta)\asymp s^2\log^{2/\al}\left(1+\frac{j_0^{\al/2}}{s^\al}\right)$$
and $$\Phi_{\rm adp}(s;\eta)\asymp \left\{
\begin{array}{cc}
   \ds s^2\log^{2/\al}\left(1+\frac{j_0^{\al/2}\log^{\al/2}(es)}{s^\al}\right),  &  s\le \sqrt{j_0\log(ej_0)};\\
   j_0\log(ej_0),  &  \text{otherwise}.
\end{array}\right.$$

\begin{itemize}
    \item If $j_0$ is bounded, then $\Phi_{\rm adp}(s;\eta)\asymp \Phi_{\rm o}(s;\eta)$ and therefore is optimal.
    \item If $j_0$ goes to infinity, then we can verify \Cref{cond: loading vector} with $s_0\asymp\sqrt{j_0\log(e j_0)}$ and $\gamma_0=1$.
\end{itemize}
Comparing the above results with the homogeneous case \eqref{eq:homogeneous case}, we obtain
$$
\Phi_{\rm o}(s;\eta)\asymp \Phi_{\rm o}(s\wedge j_0;\mathbf{1}_{j_0}),
\qquad
\Phi_{\rm adp}(s;\eta)\asymp \Phi_{\rm adp}(s\wedge j_0;\mathbf{1}_{j_0}),
$$
which shows that an exponentially decaying loading vector behaves like a homogeneous vector with an effective dimension $j_0$.

\section{Extensions}\label{sec:discussion}

\subsection{Estimation with unknown noise variance}
The estimation problem with an unknown noise variance has been widely studied in the literature; see, for example, \cite{collier2017minimax,collier2018optimal,comminges2021adaptive,carpentier2022estimation}.
A common approach is to first construct an estimator of $\sigma^2$ and then substitute it into the original procedure in place of the true noise variance.
For example, \cite{comminges2021adaptive} proposed a median-of-means estimator:

Let $\gamma_{\al,\tau} \in (0,1/2]$ be a constant chosen sufficiently small, depending only on $\al$ and $\tau$.
Divide $\{1,\ldots,d\}$ into $m = \lfloor \gamma_{\al,\tau} d \rfloor$ disjoint subsets $B_1,\ldots,B_m$, each of cardinality
\[
|B_i| \;\ge\; k := \left\lfloor \frac{d}{m} \right\rfloor \;\ge\; \frac{1}{\gamma_{\al,\tau}} - 1.
\]
The estimator of $\sigma^2$ is then defined as
\begin{equation}\label{eq:mom estimator}
  \hat{\sigma}^2 = \operatorname{median}(\bar{\sigma}_1^2,\ldots,\bar{\sigma}_m^2),
  \qquad \text{ where }
  \bar{\sigma}_i^2 = \frac{1}{|B_i|}\sum_{j\in B_i} y_j^2,
  \quad i=1,\ldots,m.
\end{equation}

\cite{comminges2021adaptive} showed that $\hat{\sigma}^2$ is a consistent estimator of $\sigma^2$ if the noise distribution belongs to ${\cal G}_{\al,\tau}$.

Based on the estimator $\hat{\sigma}^2$ in \eqref{eq:mom estimator}, most of our results for estimating the linear functional $L(\theta)$ extend naturally to the setting with an unknown noise variance.
Due to space limitations, we restrict attention to the estimation problem with known sparsity.
Replacing $\sigma$ in the estimator $\hat{L}_s$ defined in \eqref{eq:minimax nonadaptive estimator} with $\hat{\sigma}$ defined in \eqref{eq:mom estimator} and adjusting the constant in front of the threshold,  we obtain the estimator
\[
\hat{L}_s^\prime
= \sum_{j \le j_1(s)} \eta_j y_{j}
  + \sum_{j > j_1(s)} \eta_j y_j \,\mathbbm{1}\!\left\{ |\eta_j y_j| > \factorone \sqrt{2}\hat{\sigma}\tau \lambda_{o} \right\}.
\]
The next theorem provides an upper bound on its risk.

\begin{theorem}\label{thm: unknown sigma}
For any $\al,\tau > 0$, there exist constants $\gamma_{\al,\tau}$ and $C$ depending only on $\al$ and $\tau$, such that if $1 \le s < \lfloor \gamma_{\al,\tau} d \rfloor / 4$, then
\[
\sup_{\theta \in \Theta_{s}}
\sup_{{\caP}_\xi \in {\cal G}_{\al,\tau}^\otimes}
\mathbb{E}_{\theta,{\caP}_\xi}\!\left(\hat{L}_s^\prime - L(\theta)\right)^2
\;\le\; C \sigma^2 \Phi_{\rm o}(s;\eta).
\]
\end{theorem}

Together with \Cref{thm: nonadaptive linear}, this theorem shows that the same minimax rate is attainable under an unknown noise level as in the known-noise case, provided the sparsity level is not too large.

\subsection{Linear hypothesis testing}
The results on estimating the linear functional $L(\theta)$ can be applied to the problem of testing linear hypotheses over $\Theta_{s}$. For simplicity, we focus on the case of Gaussian noise, i.e., $\caP_\xi = {\cal N}^\otimes$.
For any $t_0\in \mathbb{R}$, consider testing the null hypothesis
\begin{equation}\label{eq:null hypothesis}
  H_0: \theta \in \Theta_{s,0} = \left\{ \theta \in \Theta_{s} : L(\theta) = t_0 \right\}
\end{equation}
against the alternative
\begin{equation}\label{eq:alternative hypothesis}
  H_1: \theta \in \Theta_{s}(\rho) = \left\{ \theta \in \Theta_{s} : |L(\theta)-t_0| \ge \rho \right\},
\end{equation}
for some $\rho > 0$.
Our analysis does not depend on the value of $t_0$.

Let $\Delta$ be a test with values in $\{0,1\}$.
The risk of $\Delta$ is defined as the sum of the type I error and the maximal type II error:
\[
\sup_{\theta\in \Theta_{s,0}}\bbP_{\theta, {\caN}^\otimes}(\Delta=1) \;+\; \sup_{\theta \in \Theta_{s}(\rho)} \bbP_{\theta, {\caN}^\otimes}(\Delta=0).
\]

A benchmark quantity is the \emph{minimax risk of testing}, defined as
\[
\mathcal{R}_s(\rho)
= \inf_{\Delta} \left\{ \sup_{\theta\in \Theta_{s,0}}\bbP_{\theta, {\caN}^\otimes}(\Delta=1) \;+\; \sup_{\theta \in \Theta_{s}(\rho)} \bbP_{\theta, {\caN}^\otimes}(\Delta=0) \right\},
\]
where the infimum is taken over all $\{0,1\}$-valued tests.

The \emph{minimax separation rate for testing} $H_0:\theta\in\Theta_{s,0}$ versus $H_1: \theta \in \Theta_{s}(\rho)$ is the smallest $r_{\min} > 0$ such that the following two properties hold:
\begin{enumerate}
\item[(i)] For any $\varepsilon \in (0,1)$, there exists $A_\varepsilon > 0$ such that, for all $A > A_\varepsilon$,
$\mathcal{R}_s(A r_{\min}) \le \varepsilon$;
\item[(ii)] For any $\varepsilon \in (0,1)$, there exists $a_\varepsilon > 0$ such that, for all $0 < A < a_\varepsilon$,
$\mathcal{R}_s(A r_{\min}) \ge 1 - \varepsilon$.
\end{enumerate}
The next theorem provides the exact expression for the minimax separation rate.

\begin{theorem}\label{thm: testing}
  For any integers $s$ and $d$ such that $2 \le s \le d$, and any $\sigma > 0$,
  the minimax separation rate for testing $H_0:\theta\in\Theta_{s,0}$ versus $H_1: \theta \in \Theta_{s}(\rho)$ is
  $
  r_{\min} = \sigma \sqrt{\Phi_{\rm o}(s;\eta)}$.
\end{theorem}

The proof of the lower bound in \Cref{thm: testing} follows arguments similar to those in \Cref{thm: nonadaptive linear}; for completeness, the details are provided in
Section~C.5 in the supplementary material.

For the upper bound, recall that the estimator $\hat{L}_s$ defined in \eqref{eq:minimax nonadaptive estimator} achieves the optimal squared error rate $\sigma^2 \Phi_{\rm o}(s;\eta)$ over $\Theta_{s}$.
Using this estimator, we construct the test
\begin{equation}\label{eq:test}
  \Delta_s = \mathbbm{1}\!\left\{\, |\hat{L}_s-t_0| > B\sigma \sqrt{\Phi_{\rm o}(s;\eta)} \,\right\},
\end{equation}
for some constant $B > 0$, and we show in
Section~A.6 in the supplementary material

that $\Delta_s$ achieves the minimax rate for testing with an appropriate choice of $B$.
The proof also reveals that the results extend straightforwardly to the broader class of symmetric sub-Weibull noise distributions.

\section{Future Directions}\label{sec:future}

Because the sequence model in \Cref{eq:sequence model} is closely related to high-dimensional regression, a natural first direction is to investigate the implications of our results for inference on linear functionals in regression models. We briefly mention several questions in this direction.

\begin{enumerate}
    \item Applying our results to \Cref{eg: linear model prediction} shows that, under orthogonal designs, the minimax rate is determined by the induced loading vector \(\eta_j=b_j^{-1}x_{0,j}\). Thus, for a fixed target \(x_0\), the orthogonal-design reduction yields a  benchmark through the quantity \(\Phi_o(s;\eta)\), which reflects both sparsity level \(s\) and the induced loading structure. It would be interesting to understand whether this benchmark could help compare candidate design matrices in experimental design problems \cite{pukelsheim2006optimal}.

    \item More broadly, our sharp lower bounds may provide useful benchmarks beyond the orthogonal-design setting. Existing minimax theory for sparse regression often relies on structural assumptions on the loading vector \(x_0\), such as bounded coordinate ratios or polynomial decay \cite{cai2017confidence,cai2021optimal}. By contrast, our results show that, in the orthogonal reduction, a sharp minimax rate can be formulated for arbitrary induced loadings \(\eta\) without such assumptions. It would therefore be of interest to determine whether similar sharp characterizations remain possible for more general designs.

    \item Related regression procedures beyond the Gaussian setting often exploit additional structure in the error distribution, such as symmetry or exchangeability \cite{hartigan1970exact,meinshausen2015group,lei2021assumption}, as well as exact, permutation, or conformal ideas \cite{guan2024conformal,wen2025residual}. It would be interesting to understand whether the heterogeneous minimax benchmark identified here can be combined with such techniques to obtain sharp guarantees for inference on linear functionals in broader regression models.
\end{enumerate}

Another natural extension is to allow correlated noise. Recent work has obtained sharp results under equicorrelated noise \cite{kotekal2023minimax,kotekal2025sparsity}, while more general dependence structures have been studied in \cite{liu2021minimax,wang2024correlated}. Extending the present heterogeneous-loading theory beyond product noise would require new lower- and upper-bound arguments that explicitly incorporate the covariance structure.

A more specific analytical question concerns the remaining logarithmic gap under non-symmetric noise.
Unlike in the symmetric-noise case, thresholding need not center the null coordinates at zero.
This suggests that the current upper-bound analysis may be too loose, and that a sharper analysis may be needed to close the gap.

Finally, it is natural to extend the present framework beyond linear functionals. A first step would be to study weighted polynomial-type functionals such as \(\sum_{j=1}^d \eta_j \theta_j^m\). Existing work \cite{collier2020estimation} used polynomial approximation methods to analyze more general smooth even functionals via expansions in even-order terms. It would be of interest to determine whether similar approximation-theoretic ideas can be combined with the heterogeneous calibration and divergence bounds developed here to treat a broader class of functionals.

\section{Proof of the lower bounds}\label{sec:proof lower}
In this section, we prove the lower bounds in \Cref{thm: nonadaptive linear,thm: adaptive lower}.

Throughout, we set \(\sigma=1\) without loss of generality. For probability measures \(P\) and \(Q\) on the same measurable space \((\mathcal X,\mathcal U)\), define the total variation distance, Hellinger distance, and chi-square divergence by

\[
\mathrm{TV}(P,Q)=\sup_{B\in\mathcal U}|P(B)-Q(B)|, \mathrm{H}(P,Q)=\sqrt{\int\left(\sqrt{\frac{\diff P}{\diff Q}}-1\right)^2\,\diff Q}, \chi^2(P\,\|\,Q)=\int\left(\frac{\diff P}{\diff Q}-1\right)^2\,\diff Q.
\]

\subsection{General Tools}

For nonadaptive estimation, we prove the lower bounds using Le Cam's method \cite{lecam1973convergence}. Let $\mu$ be a probability measure on $\bbR^d$. Denote by $\bbP_{\mu,\caP_\xi}$ the mixture probability measure:
$$\bbP_{\mu,\caP_\xi}=\int_{\bbR^d}\bbP_{\theta,\caP_\xi}\, \mu(\diff \theta).$$
Specifically, we denote by $\bbP_{0,\caP_\xi}$ the mixture probability measure when $\mu$ is the Dirac measure at $0$, i.e., $\bbP_{0,\caP_\xi}=\bbP_{\theta=0,\caP_\xi}$.
The following lemma, as a special case of \cite[Theorem 2.15]{tsybakov2009introduction}, is the key tool we will use to prove the lower bounds.
\begin{lemma}\label{lem: general tool}
  Suppose that there exist two distributions $P_1=\bbP_{\mu_1,\caP_\xi^1}$ and $P_2=\bbP_{\mu_2,\caP_\xi^2}$ such that
  \[\mu_1(L(\theta)\le c,\theta\in \Theta)\ge 1-\beta_1,\quad \mu_2(L(\theta)\ge c+2t,\theta\in \Theta)\ge 1-\beta_2,\quad \mathrm{TV}(P_1,P_2)\le \beta_3\] for some $c\in \bbR$, $t>0$, and $\beta_1,\beta_2,\beta_3\in (0,1)$, and some subset $\Theta\subseteq \bbR^d$. Then, for any estimator $\hat{T}$ of $L(\theta)$, we have
$$\inf_{\hat{T}}\sup_{\theta\in \Theta}\max_{j\in \{1,2\}} \bbP_{\theta, \caP_{\xi}^j}(|\hat{T}-L(\theta)|\ge t)\ge \frac{1-\beta_1-\beta_2-\beta_3}{2}.$$
\end{lemma}

Since $\mathbb{E}(\hat{T}-L(\theta))^2 \ge t^2 \bbP(|\hat{T}-L(\theta)| \ge t)$, \Cref{lem: general tool} immediately yields a lower bound on the minimax mean squared error. To control the total variation distance, we introduce the following lemmas. Lemma~\ref{lem: chi square} is the main tool for bounding the $\chi^2$ divergence in the symmetric case, while Lemma~\ref{lem: hellinger control} is used to control the Hellinger divergence in the asymmetric case.

\begin{lemma}\label{lem: chi square}
For $\al> 0$, there exists a constant $\tau_\al>0$ and a distribution in ${\cal G}_{\al,\tau_\al}$ with density $f_\al^{(0)}$ such that, for
\( f_\al^{(1)}(\cdot) := f_\al^{(0)}(\cdot - \gamma) \) with some \( \gamma \in \mathbb{R} \), it holds that
\[
1+\chi^2(f_\al^{(1)} \,\|\, f_\al^{(0)}) \;\le\; C_{\al, 1}\exp\!\left( \left|\frac{\gamma}{C_{\al, 2}}\right|^\al\right),
\]
where $C_{\al, 1},C_{\al, 2}$ are positive constants depending only on $\al$. For $\al\le 1$ and $\al=2$, $C_{\al, 1}=1$.
\end{lemma}

When $f_0$ and $f_1$ are Gaussian distributions, the corresponding $\chi^2$-divergence admits a closed-form expression, which has been widely used in the minimax theory of Gaussian models \cite{collier2017minimax,collier2018optimal,cai2017confidence,cai2018accuracy,cai2021optimal,bradic2022testability,liu2021minimax}. Lemma~\ref{lem: chi square} extends this idea to the exponentially decaying noise setting, where no closed-form formula is available. Nevertheless, the optimality of our results highlights the effectiveness of Lemma~\ref{lem: chi square}.
\begin{lemma}[{\cite[Theorem 7.6]{ibragimov2013statistical}}]\label{lem: hellinger control}
Suppose $\mathcal{F}=\{p(\cdot; \theta):\theta\in \mathbb{R}^p\}$ is a family of density functions w.r.t. a $\sigma$-finite measure $\nu$ on $\mathbb{R}^d$ equipped with the Borel algebra $\mathcal{B}$.
If the experiment $(\mathbb{R}^d, \mathcal{B}, \mathcal{F})$ is regular (see \cite[Chapter 7.1]{ibragimov2013statistical} for a definition),

then
  \begin{align*}
    \int_{\mathbb{R}^d} \left[\sqrt{p(x;\theta)}-\sqrt{p(x;\theta+h)}\right]^2
    \nu(\diff x)
    \le \frac{|h|^2}{4}\int_0^1 {\rm tr}\left[ I(\theta+ sh)\right] \diff s, \quad \forall\theta, h\in \mathbb{R}^p,
  \end{align*}
where $I(\vartheta)$ denotes the Fisher information matrix of the family $\mathcal{F}$ at parameter value $\vartheta$.
\end{lemma}
\subsection{Proof of the lower bounds in Theorem~\ref{thm: nonadaptive linear}}\label{sec:proof nonadaptive linear}

Without loss of generality, assume $\sigma=1$.
From Lemma~\ref{lem: chi square}, for any $\al > 0$, there exists $\tau_\al > 0$ and a distribution $f_\al^0$ such that $f_\al^0 \in {\cal G}_{\al,\tau}$ for all $\tau \ge \tau_\al$.
For simplicity, let $F_\al^{\otimes}$ denote the product distribution on $\mathbb{R}^d$ with i.i.d. marginals having density $f_\al^0$.

We apply Lemma~\ref{lem: general tool} with $\Theta = \Theta_{s}$, taking
\( P_1 = \mathbb{P}_{0, F_\al^{\otimes}} \)
and
\( P_2 = \mathbb{P}_{\mu, F_\al^{\otimes}} \),
where $\mu$ is a distribution on $\mathbb{R}^d$ that will be specified later.
In what follows, we establish the lower bound for the case where $s > C_1$ and $\lambda_{o} s+\nu> C_2|\eta_1|$ with constants $C_1,C_2 > 0$ to be specified.
The proof for the complementary case is provided in the supplementary material.

Define the probability measure \(\mu\) on $\bbR^d$ as follows: for \(\theta\sim \mu\), its coordinates are independent and
  $$\theta_j=b_j\gamma_j,\quad \forall j\in [d],$$
where $b_j={\rm Ber}(\pi_j)$,
$$\pi_j=c_1\cdot\frac{ |\eta_j|\exp(-\beta_+/|\eta_j|^\al)}{\sqrt{\sum_{i=1}^d \eta_i^2\exp(-\beta_+/|\eta_i|^\al)}},\quad\text{and}\quad \gamma_j=\left\{\begin{array}{lc}
   C_{\al, 2} {\rm sign}(\eta_j)& j\le j_1(s)\\
   C_{\al, 2} {\lambda_{o}}/{\eta_j}& j> j_1(s)
\end{array}\right.$$
for some constant $c_1 > 0$ to be specified later, where the constant $C_{\al, 2}$ is as given in \Cref{lem: chi square}. It can be easily verified that $\pi_j\in (0,1)$.
From \Cref{eq:solution nonadaptive}, we have
$$\bbE_{\mu}\left\|\theta\right\|_0=\sum_{j=1}^d\pi_j\le \frac{c_1}{2}s,\quad\text{and}\quad{\rm Var}_{\mu}\left\|\theta\right\|_0\le \sum_{j=1}^d\pi_j(1-\pi_j)\le\sum_{j=1}^d\pi_j\le \frac{c_1}{2}s.$$

Since $\eta_j\gamma_j= C_{\al, 2}\max(|\eta_j|, \lambda_{o})$ is non-increasing in $j$, we have
$$
{\rm Var}_{\mu} [L(\theta)]= \sum_{j=1}^d\eta_j^2\gamma_j^2\cdot \pi_j[1-\pi_j]
  \le \eta_1\gamma_1\cdot\sum_{j=1}^d\eta_j\gamma_j\pi_j=  C_{\al, 2}\max(|\eta_1|, \lambda_{o}) \cdot\bbE_{\mu}L(\theta).
$$
Using $\eta_j\gamma_j\ge  C_{\al, 2} |\eta_j|$, we have $\bbE_{\mu} L(\theta)=\sum_{j=1}^d\eta_j\gamma_j\pi_j\ge  C_{\al, 2} \sum_{j=1}^{d}|\eta_j|\pi_j= c_1 C_{\al, 2}  \nu$. Using $\eta_j\gamma_j\ge C_{\al, 2} \lambda_{o}$, we have $\bbE_{\mu} L(\theta)=\sum_{j=1}^d\eta_j\gamma_j\pi_j\ge   C_{\al, 2} \lambda_{o} \sum_{j=1}^{d}\pi_j$.
If $\lambda_{o}>0$, we have $\sum_{j=1}^{d}\pi_j=\frac12 c_1s$ and thus $\bbE_{\mu} L(\theta)\ge \frac12 c_1  C_{\al, 2} \lambda_{o}s$.
This inequality also holds when $\lambda_{o}=0$.

Consequently, we have
$$\bbE_{\mu}L(\theta)\ge \frac{c_1}{4} C_{\al, 2}  [\lambda_{o} s+\nu]>0.$$
By Chebyshev's inequality, we have, for any $c_1\in (0,2)$,
$$
\mu\left(\|\theta\|_0>s\right)\le \mu\left(\|\theta\|_0-\mathbb{E}_{\mu}\|\theta\|_0> (1-\frac{c_1}{2}) s\right) \le \frac{{\rm Var}_{\mu}[\left\|\theta\right\|_0]}{(1-c_1/2)^2s^2}\le \frac{c_1/2}{(1-c_1/2)^2s},
$$
and
$$
\begin{aligned}
\mu\left(L(\theta)< \frac{c_1}{8}C_{\al, 2} [\lambda_{o} s+\nu]\right)& \le \mu\left(L(\theta)<\frac{1}{2}\bbE_{\mu}L(\theta)\right) \le \frac{4{\rm Var}_{\mu} L(\theta)}{\left[\bbE_{\mu}[L(\theta)]\right]^2}\le \frac{16 \max(|\eta_1|, \lambda_{o})}{c_1 (\lambda_{o}s+\nu)}.
\end{aligned}
$$
Therefore, $\forall c_1\in (0,2),\al_1\in (0,1)$, we can choose $C_1,C_2$ large enough such that
$$\forall s\ge C_1,\lambda_{o} s+\nu\ge C_2|\eta_1|,\quad\mu\left(\theta\in \Theta_{s}, L(\theta)\ge \frac{c_1}{8} C_{\al, 2}[\lambda_{o} s+\nu]\right)\ge 1-\al_1.$$

Let $P_2 = \mathbb{P}_{\mu, F_\al^{\otimes}}$ and let $f_\al^j$ be the distribution of the $j$-th component under $P_2$.
Lemma~\ref{lem: chi square} implies that
  \begin{align*}
    1+\chi^2(P_2\|P_1)&=\prod_{j=1}^d \left[1+\chi^2\left(f_\al^j\|f_\al^{(0)}\right)\right]\\
    &\stackrel{(1)}{=}\prod_{j=1}^d\left[1+\pi_j^2 \cdot \chi^2(f_\al^{(0)}(\cdot-\gamma_j)\|f_\al^{(0)})\right]\\
    &\stackrel{(2)}{\le} \exp\left[\sum_{j=1}^d \pi_j^2\chi^2(f_\al^{(0)}(\cdot-\gamma_j)\|f_\al^{(0)})\right]\\
    &\stackrel{(3)}{\le} \exp\left[c_1^2C_{\al, 1}\frac{e\sum_{j\le j_1(s)} \eta_j^2\exp(-2\beta_+/|\eta_j|^\al)+\sum_{j>j_1(s)}\eta_j^2\exp(-\beta_+/|\eta_j|^\al)}{\sum_{j=1}^d\eta_j^2\exp(-\beta_+/|\eta_j|^\al)}\right]\\
    &\le \exp(c_1^2C_{\al, 1}e),
  \end{align*}
  where we have applied Fubini's theorem in (1),
  $1+x\le e^x$ in (2),
  and Lemma~\ref{lem: chi square} in (3).
  Therefore, for any $\al_2\in (0,1-\al_1)$, we can choose some small $c_1$ such that $\mathrm{TV}(P_1,P_2)\le \sqrt{\chi^2(P_2\|P_1)}/2\le \al_2$.
Specifically, we can take $\al_1=\al_2=1/4$ and fix the constants $c_1$, $C_1$, and $C_2$.
Using Lemma~\ref{lem: general tool} with $t= 2^{-4}c_1 C_{\al, 2}(\lambda_{o} s+\nu)$, we have
$$
\inf_{\hat{L}_s} \sup_{\theta\in\Theta_{s}}\sup_{{\caP}_\xi\in  \mathcal{G}_{\al,\tau}^\otimes}
    \mathbb{P}_{\theta, {\caP}_\xi} \left(\! | \hat{L}_s - L(\theta) |
    \;\ge\; 2^{-4} c_1 C_{\al, 2}(\lambda_{o} s+\nu)  \right) \ge \frac14.
$$

\subsection{Proof of Theorem~\ref{thm: adaptive lower}}
  The lower bound for the adaptive estimation requires the following new technical tools.
  \begin{lemma}\label{lem: adaptive lower}
    Let \(P, Q_1\) and $Q_2$ be three probability measures defined over the same measurable space \((\mathcal{X}, \mathcal{U})\). Then for any $q>0$, we have
    \begin{equation}
      \inf_{{\cal A}\in {\cal U}}\left\{P({\cal A})q+Q_1({\cal A}^c)\right\}\ge \max_{0<\tau<1}\left[\frac{q\tau}{1+q\tau}\left(1-\tau(\chi^2(Q_2||P)+1)\right)-{\rm TV}(Q_1,Q_2)\right].
    \end{equation}
  \end{lemma}
  \Cref{lem: adaptive lower} follows directly from
  \cite[Lemma 8]{collier2018optimal} since for any $A\in \mathcal{U}$, we have $P({\cal A})q+Q_1({\cal A}^c)\ge P({\cal A})q+Q_2({\cal A}^c)-{\rm TV}(Q_1, Q_2)$.

For any integer $s\in[d\wedge s_0]$, denote
$A_s=s\lambda_*(s)+\nu_*(s)$.
We have
$\Phi_*(s;\eta)\le A_s^2\le2\Phi_*(s;\eta)$.

\begin{lemma}
\label{lem:adaptive-prior}
Fix $\al>0$ and $\gamma\in(0,2)$.
Let $F_\al^\otimes$ be defined as in
\Cref{sec:proof nonadaptive linear}.
There exist constants
$c_0,K_1,K_2>0$ (depending on
$\al$ and $\gamma$) such that if
$s_0\ge s\ge K_1$
and $A_s\ge K_2\log(es)\max\{|\eta_1|,\lambda_*(s)\}$,
then there exists a probability measure $\mu$ on $\mathbb R^d$ satisfying that
$\mu\{\theta\in\Theta_s:L(\theta)\ge c_0A_s\}\ge\frac34$,
and  $1+\chi^2\bigl(\bbP_{\mu,F_\al^\otimes}
\,\|\,\bbP_{0,F_\al^\otimes}\bigr)\le s^\gamma$.
\end{lemma}
The proof of \Cref{lem:adaptive-prior} is provided in the supplementary material.

Throughout, $\al$ and $\gamma$ are fixed, and the estimator $\widehat T$ satisfies
\eqref{eq:adaptive lower small}.
If $s>s_0$, we apply the argument below with $s_0$ in place of $s$.
Indeed, $\Theta_{s_0}\subseteq\Theta_s$,
$\Phi_{\rm adp}(s;\eta)=\Phi_{\rm adp}(s_0;\eta)$, and
$s\wedge s_0=s_0$. Hence, it suffices to consider $s\le s_0$.
First of all, the lower-bound part of \Cref{thm: nonadaptive linear}, applied with
sparsity one, implies that for some constant $c_\al>0$,
\begin{equation}\label{eq:adaptive-lower-one-sparse}
\sup_{\theta\in\Theta_1}
\sup_{{\caP}_\xi\in{\cal G}_{\al,\tau}^{\otimes}}
\bbE_{\theta,{\caP}_\xi}
\bigl(\widehat T-L(\theta)\bigr)^2
\ge c_\al\Phi_{\rm o}(1;\eta).
\end{equation}

The rest of the proof consists of three steps.

\textbf{Step 1:} \textit{Reduction to the nontrivial case.}
Let $c_0,K_1,K_2$ be the constants in \Cref{lem:adaptive-prior}. We may increase $K_1$ so that $s\ge K_1$ implies
$s\ge2K_2\log(es)$.
We aim to restrict attention to the regime where
\begin{equation}\label{eq:adaptive-lower-prior-regime-main}
s\ge K_1,
\qquad
\Phi_{\rm adp}(s;\eta)=\Phi_*(s;\eta),
\qquad
A_s\ge K_2\log(es)\max\{|\eta_1|,\lambda_*(s)\}.
\end{equation}
In the other cases listed below,
we will show
$\frac{\Phi_{\rm adp}(s;\eta)}{s^\gamma}\lesssim \Phi_{\rm o}(1;\eta)$ and  the result is proved using \eqref{eq:adaptive-lower-one-sparse}.

\textit{(1).}
If $s<K_1$, we have
$s^{2-\gamma}\log(es)\le K_1^{2-\gamma}\log(eK_1)$, which only depends on
$\al$ and $\gamma$.
The last comparison in \Cref{prop:linear decrease} gives
$\frac{\Phi_{\rm adp}(s;\eta)}{s^\gamma}
\lesssim
s^{2-\gamma}\log(es)\Phi_{\rm o}(1;\eta)
\lesssim
\Phi_{\rm o}(1;\eta)$.

\textit{(2).}
If $0<\al<2$, $K_1\le s\le s_0\wedge d$, and
$\Phi_*(s;\eta)\le\Phi_{*}(1;\eta)\log^2(es)$, then
$$\frac{\Phi_{\rm adp}(s;\eta)}{s^\gamma}
=\frac{\Phi_{*}(1;\eta)\log^2(es)}{s^\gamma}
\lesssim\Phi_{\rm o}(1;\eta).$$

\textit{(3).}
If it holds that $K_1\le s \le s_0\wedge d$,
$\Phi_{\rm adp}(s;\eta)=\Phi_*(s;\eta)$, and
$A_s<K_2\log(es)\max\{|\eta_1|,\lambda_*(s)\}$, then
we must have $\lambda_*(s)<|\eta_1|$; otherwise, $\lambda_*(s)\ge|\eta_1|$ implies that
$s\lambda_*(s)\le A_s<K_2\log(es)\lambda_*(s)$, which contradicts
$s\ge K_1\ge 2K_2\log(es)$. Therefore, $\lambda_*(s)<|\eta_1|$ implies that
$\Phi_{\rm o}(1;\eta)\gtrsim|\eta_1|^2$.
We have
$\frac{\Phi_{\rm adp}(s;\eta)}{s^\gamma}
\le\frac{A_s^2}{s^\gamma}
<K_2^2\frac{\log^2(es)}{s^\gamma}|\eta_1|^2
\lesssim \Phi_{\rm o}(1;\eta)$.

\smallskip
In the following, we focus on the regime in
\eqref{eq:adaptive-lower-prior-regime-main}.

\textbf{Step 2:}
\textit{Prior and constrained testing.}
By \Cref{lem:adaptive-prior}, we have a probability measure $\mu$ such that
\begin{equation}\label{eq:adaptive-lower-good-set-main}
\mu(\mathcal E_s)\ge\frac34, \quad \text{ where } \quad
\mathcal E_s
=\{\theta\in\Theta_s:L(\theta)\ge c_0A_s\},
\end{equation}
and
$1+\chi^2\bigl(\bbP_{\mu,F_\al^\otimes}
\,\|\,\bbP_{0,F_\al^\otimes}\bigr)\le s^\gamma$.
Let $\mu_1$ be the restriction of $\mu$ on  $\mathcal E_s$, i.e., $\mu_1(B)=\frac{\mu(B\cap \mathcal E_s)}{\mu(\mathcal E_s)}$ for every measurable $B\subseteq\mathbb R^d$.
Let $P=\bbP_{0,F_\al^\otimes}$,
$Q_2=\bbP_{\mu,F_\al^\otimes}$, and
$Q_1=\bbP_{\mu_1,F_\al^\otimes}$.
By data processing inequality,
we have $\mathrm{TV}(Q_1,Q_2)
\le\mathrm{TV}(\mu_1,\mu)
\le\mu(\mathcal E_s^c)
\le\frac14$.
Applying \Cref{lem: adaptive lower} with $q=4s^\gamma$ and evaluating
its right-hand side at $u=s^{-\gamma}/2$ gives
\begin{equation}\label{eq:adaptive-lower-testing-main}
\inf_{{\cal A}\in\mathcal U}
\{qP({\cal A})+Q_1({\cal A}^c)\}\ge \frac23\left(1-\frac12\right)-\frac14
=\frac1{12}.
\end{equation}

\textbf{Step 3:}
\textit{Conversion to estimation risk.}
Define
${\cal A}=\{\widehat T\ge c_0A_s/2\}$.
For $\mu_1$-almost every $\theta$, we have
$L(\theta)\ge c_0A_s$. For such $\theta$,
${\cal A}^c$ implies
$|\widehat T-L(\theta)|\ge\frac{c_0A_s}{2}$.
Using the definition of $Q_1$ and \eqref{eq:adaptive lower small}, we obtain
\begin{align*}
Q_1({\cal A}^c)
=\int
\bbP_{\theta,F_\al^\otimes}({\cal A}^c)\,
\mu_1(\diff\theta)
\le
\frac{4}{c_0^2A_s^2}
\int
\bbE_{\theta,F_\al^\otimes}
\bigl(\widehat T-L(\theta)\bigr)^2\,
\mu_1(\diff\theta)
\le
\frac{4}{c_0^2C_0}
\frac{\Phi_{\rm adp}(s;\eta)}{A_s^2}
\le
\frac{4}{c_0^2C_0},
\end{align*}
where the last inequality follows from
$\Phi_{\rm adp}(s;\eta)=\Phi_*(s;\eta)\le A_s^2$ in
\eqref{eq:adaptive-lower-prior-regime-main}.
Choosing $C_0\ge96/c_0^2$, we have
$Q_1({\cal A}^c)\le1/24$.
It follows from \eqref{eq:adaptive-lower-testing-main} that
$P({\cal A})\ge1/(96s^\gamma)$.
Since$P=\bbP_{0,F_\al^\otimes}$, we have
\begin{align*}
\sup_{\theta\in\Theta_1}
\sup_{{\caP}_\xi\in{\cal G}_{\al,\tau}^{\otimes}}
\bbE_{\theta,{\caP}_\xi}
\bigl(\widehat T-L(\theta)\bigr)^2
\ge
\bbE_{0,F_\al^\otimes}\left[\widehat T^2\right]
&\ge
\frac{c_0^2A_s^2}{4}P({\cal A})\\
&\ge
\frac{c_0^2A_s^2}{384s^\gamma}
\ge
\frac{c_0^2}{384}
\frac{\Phi_{\rm adp}(s;\eta)}{s^\gamma},
\end{align*}
where the last inequality uses
$\Phi_{\rm adp}(s;\eta)=\Phi_*(s;\eta)\le A_s^2$.
Combining this bound with \eqref{eq:adaptive-lower-one-sparse} proves
\eqref{eq:adaptive-lower}.

\section*{Supplementary material}
\label{SM}
The supplementary material contains complete technical details and proofs that are omitted from the main text for brevity, such as the proof of \Cref{thm: exponential lower}, derivations of the upper bounds, propositions, and the analytic rates in the example section.
\bibliographystyle{amsplain}
\bibliography{ref}
\newpage
\appendix
\setcounter{section}{0}
\pagenumbering{arabic}
\setcounter{page}{1}
\bigskip
\bigskip
\begin{center}
\textbf{\Large Supplement to ``Minimax and adaptive estimation of general linear functionals under sparsity''}
\end{center}
\numberwithin{equation}{section}
\numberwithin{theorem}{section}
\numberwithin{lemma}{section}
\numberwithin{proposition}{section}
\numberwithin{corollary}{section}
\numberwithin{assumption}{section}
\numberwithin{condition}{section}
\numberwithin{definition}{section}
\numberwithin{remark}{section}
\numberwithin{example}{section}

\Cref{sec:pf-upper} provides detailed proofs of all upper bound results in the paper.
\Cref{sec:pf-aux} provides the proof for the uniqueness of the solution to \Cref{eq:solution nonadaptive}
and includes the proofs of the propositions on the rate characteristics.
\Cref{sec:additional lower bound} includes the proofs of Theorem~\ref{thm: exponential lower} and the lemmas used to prove the lower bounds in \Cref{sec:proof lower}.
\Cref{sec:pf-example} contains the details for the examples in \Cref{sec:example}.
\Cref{app:comparison-chhor} compares our lower-bound construction with that of \cite{chhor2024sparse}.

In the appendix, we will denote by \( C_i, i = 1, 2, \ldots \), absolute positive constants, and by \( C, c \) absolute positive constants that may vary from line to line.
\section{Proof of the upper bounds}\label{sec:pf-upper}
We will use the following lemmas.
{Throughout this section, the proofs follow the same bias-variance decomposition as in the main text. The only additional work is to control the tail contributions uniformly over the noise classes $\mathcal{G}_{\al,\tau}^{\otimes}$ and $\mathcal{H}_{\al,\tau}^{\otimes}$.}
\begin{lemma}\label{lem: tail}
  For $X\sim P$ for some $P\in {\cal G}_{\al,\tau}$ (or ${\cal H}_{\al,\tau}$), we have
\begin{align}
    \mathbb{E}X^{2q}\le C_{q}^*, q\in \mathbb{N}, \label{eq:xn-tail}
\end{align}
where $C_{q}^*$ are positive constants depending on $\al$ and $\tau$.
\end{lemma}
\Cref{lem: tail} can be proved using the tail condition in \Cref{def: symmetric exponential} and integration by parts.
\begin{lemma}\label{lem: probability bound}
  For $\al,\tau>0$, if $\zeta$ is chosen sufficiently large, then for $\hat{s}$ defined in \eqref{eq:selected index}, we have
  \begin{eqnarray}
    \sup_{\theta\in \Theta_s}\sup_{\mathcal{P}_\xi\in {\cal G}_{\al,\tau}^\otimes} \bbP_{\theta, \mathcal{P}_\xi}\left(\hat{s}>s\right)\le c_{3,0} s^{-7}, \qquad \forall s\le s_*,  \label{eq:probability bound}
  \end{eqnarray}
  where $c_{3,0}$ is a constant that depends only on $\al$ and $\tau$.
\end{lemma}
\begin{lemma}\label{lem: fourth bound}
  For $\al,\tau>0$, there exists a constant $c_{4,0}>0$ such that
  $$\sup_{\theta\in \Theta_s}\sup_{{\caP}_\xi\in {\cal G}_{\al,\tau}^\otimes}\bbE_{\theta, {\caP}_\xi}\left(\hat{L}_{s}^*-L(\theta)\right)^4\le c_{4,0}\sigma^4\Phi_{\rm adp}^2(1;\eta)s^4\log^2(es),\quad \forall s\le s_0\wedge d.$$
\end{lemma}
\begin{lemma}[{\cite[Proposition 1]{comminges2021adaptive}}]\label{lem: mom}
For any $\al,\tau>0$,  there exist constants $\gamma_{\al,\tau}\in (0,1/2]$, $c_{\al,\tau}^* > 0$, and $C_{\al,\tau}^* > 0$ depending only on $\al$ and $\tau$ such that, for any integers $s$ and $d$ satisfying $1 \le s < \lfloor \gamma_{\al,\tau} d \rfloor/4$ and $\hat{\sigma}$ defined in \eqref{eq:mom estimator}, we have
\[
\inf_{\theta\in \Theta_s}\inf_{{\caP}_\xi\in {\cal G}_{\al,\tau}^\otimes}
\bbP_{\theta,\mathcal{P}_\xi}\!\left( \frac{1}{2} \le \frac{\hat{\sigma}^2}{\sigma^2} \le \frac{3}{2} \right)
\;\ge\; 1 - \exp(-c_{\al,\tau}^* d),
\]
and
\[
\sup_{\theta\in \Theta_s}\sup_{{\caP}_\xi\in {\cal G}_{\al,\tau}^\otimes}
\frac{\bbE_{\theta,\mathcal{P}_\xi}\!\bigl|\hat{\sigma}^2 - \sigma^2\bigr|}{\sigma^2}
\;\le\; C_{\al,\tau}^*.
\]
\end{lemma}
In this section, we will use the notation $S=\left\{j:\theta_j\neq 0\right\}$ for the true parameter vector $\theta$.

\subsection{Proof of the upper bound in Theorem~\ref{thm: nonadaptive linear}}
\begin{proof}
 Note that $y_j = \theta_j + \sigma \xi_j$, where $\xi_j$ are mean-zero, unit-variance random variables.
If $\lambda_{o} = 0$, then no thresholding is applied in $\hat{L}_s$, and
\[
\hat{L}_s - L(\theta) = \sigma \sum_{j=1}^d \eta_j \xi_j.
\]
In this case, it is immediate that
\[
\forall \, \theta\in\Theta_s,\; {\caP}_\xi \in {\cal G}_{\al,\tau}^{\otimes},
\quad
\mathbb{E}_{\theta, {\caP}_\xi}\big(\hat{L}_s-L(\theta)\big)^2
= \sigma^2 \sum_{j=1}^d \eta_j^2.
\]
Hence, in the sequel, we restrict attention to the case where $\lambda_{o} > 0$.
For all ${\caP}_\xi \in {\cal G}_{\al,\tau}^{\otimes}$ and $\theta\in\Theta_s$, we have
\begin{equation}\label{eq:nonadaptive error decomposition}
\begin{aligned}
\mathbb{E}_{\theta, {\caP}_\xi}(\hat{L}_s-L(\theta))^2
\;\le\; 3\Biggl[
&\underbrace{\mathbb{E}_{\theta, {\caP}_\xi}\!\left(\sum_{j \le j_1(s)} \eta_j y_j - \sum_{j \le j_1(s)} \eta_j \theta_j \right)^2}_{\text{denoted as I}} \\
&+ \underbrace{\mathbb{E}_{\theta, {\caP}_\xi}\!\left(\sum_{\substack{j> j_1(s) \\ j \in S^c}} \sigma \eta_j \xi_j \bbo\!\left\{ |\eta_j \xi_j| \ge  \factorone \tau \lambda_{o} \right\}\right)^2}_{\text{denoted as II}} \\
&+ \underbrace{\mathbb{E}_{\theta, {\caP}_\xi}\!\left(\sum_{\substack{j> j_1(s) \\ j \in S}} \eta_j y_j \bbo\!\left\{ |\eta_j y_j| \ge \factorone \sigma \tau \lambda_{o} \right\} - \sum_{\substack{j> j_1(s) \\ j \in S}} \eta_j \theta_j \right)^2}_{\text{denoted as III}}
\Biggr].
\end{aligned}
\end{equation}

\textit{Term I.}
\[
\mathbb{E}_{\theta, {\caP}_\xi}\!\left(\sum_{j \le j_1(s)} \eta_j y_j - \sum_{j \le j_1(s)} \eta_j \theta_j \right)^2
= \sigma^2 \sum_{j \le j_1(s)} \eta_j^2.
\]

\textit{Term II.}
\begin{align*}
\mathbb{E}_{\theta, {\caP}_\xi}\!\left(\sum_{\substack{j>j_1(s) \\ j \in S^c}} \sigma \eta_j \xi_j
\bbo\!\left\{|\eta_j \xi_j| \ge \factorone \tau \lambda_{o}\right\}\right)^2
&= \sigma^2 \sum_{\substack{j>j_1(s) \\ j \in S^c}}
\mathbb{E}\!\eta_j^2 \xi_j^2 \bbo\!\left\{ |\xi_j| \ge \frac{\factorone\tau \lambda_{o}}{|\eta_j|}\right\} \\
&\stackrel{(*)}{\le} \sigma^2 \sum_{\substack{j>j_1(s) \\ j \in S^c}}
\sqrt{\mathbb{E}(\eta_j^4 \xi_j^4)}\cdot
\sqrt{\mathbb{P}\!\left(|\xi_j| \ge \frac{\factorone\tau \lambda_{o}}{|\eta_j|}\right)} \\
&\stackrel{(**)}{\le}  2 \sigma^2 \sum_{\substack{j>j_1(s) \\ j \in S^c}}
\eta_j^2 \sqrt{C_2^*} \cdot \exp(-\beta/|\eta_j|^\al) \\
&\le 2\sqrt{C_2^*}\sigma^2 \nu^2,
\end{align*}
where inequality $(*)$ follows from H\"older's inequality, inequality $(**)$ from \Cref{lem: tail} and \Cref{def: symmetric exponential}, and the last inequality is due to the definition of $\nu$ in \Cref{eq:nonadaptive rate}.

\textit{Term III.}
\begin{align*}
&\mathbb{E}_{\theta, {\caP}_\xi}\!\left(
\sum_{\substack{j>j_1(s) \\ j \in S}} \eta_j y_j \bbo\!\left\{|\eta_j y_j| \ge \factorone\sigma \tau \lambda_{o} \right\}
- \sum_{j>j_1(s)} \eta_j \theta_j \right)^2 \\
&= \mathbb{E}_{\theta, {\caP}_\xi}\!\left(
\sigma \sum_{\substack{j>j_1(s) \\ j \in S}} \eta_j \xi_j
- \sum_{\substack{j>j_1(s) \\ j \in S}} \eta_j y_j \bbo\!\left\{|\eta_j y_j| < \factorone\sigma \tau \lambda_{o} \right\}
\right)^2 \\
&\le 2s \cdot \sum_{\substack{j>j_1(s) \\ j \in S}}
\left[\sigma^2 \mathbb{E}_{\theta, {\caP}_\xi}(\eta_j \xi_j)^2
+ \mathbb{E}\!\left(\eta_j y_j \bbo\!\left\{|\eta_j y_j| < \factorone\sigma \tau \lambda_{o}\right\}\right)^2\right] \\
&\le 2s \cdot \sum_{\substack{j>j_1(s) \\ j \in S}}
\left[\sigma^2 \eta_j^2 + \factorone\factorone \sigma^2 \lambda_{o}^2 \tau^2\right] \le 2(1+\factorone\factorone\tau^2)\sigma^2 \lambda_{o}^2 s^2.
\end{align*}

Combining the bounds for I, II, and III, we obtain the desired result \eqref{eq:nonadaptive upper bound} in Theorem~\ref{thm: nonadaptive linear}.
\end{proof}
\subsection{Proof of Theorem~\ref{thm: adaptive estimator}}\label{sec:adaptive estimator}
\begin{proof}
Let $\theta\in \Theta_s, {\caP}_\xi\in {\cal G}_{\al,\tau}^\otimes$. Similar to the proof of Theorem~\ref{thm: nonadaptive linear}, we can assume $\lambda_*(s)>0$ and decompose the error as follows:
\begin{align*}
\bbE_{\theta, {\caP}_\xi}(\hat{L}_{s}^*-L(\theta))^2\le& 3\biggl[\underbrace{\bbE_{\theta, {\caP}_\xi}\biggl(\sum_{j\le j_2(s)}\eta_j y_j-\sum_{j\le j_2(s)}\eta_j\theta_j\biggr)^2}_{\text{denoted as I}}\\
&\qquad+\underbrace{\bbE_{\theta, {\caP}_\xi}\biggl(\sum_{\substack{j>j_2(s)\\j\in S^c}}\sigma \eta_j\xi_j\bbo\left\{|\eta_j\xi_j|\ge  \factorone\tau\lambda_*(s)\right\}\biggr)^2}_{\text{denoted as II}}\\
&\qquad\qquad+\underbrace{\bbE_{\theta, {\caP}_\xi}\biggl(\sum_{\substack{j>j_2(s)\\j\in S}}\eta_jy_j\bbo\left\{|\eta_jy_j|\ge \factorone\sigma\tau\lambda_*(s)\right\} -\sum_{\substack{j>j_2(s)\\j\in S}}\eta_j\theta_j\biggr)^2}_{\text{denoted as III}}\biggr].
\end{align*}

\textit{Term I.}
$$\bbE_{\theta, {\caP}_\xi}\left(\sum_{j\le j_2(s)}\eta_j y_j-\sum_{j\le j_2(s)}\eta_j\theta_j\right)^2=\sigma^2\sum_{j\le j_2(s)}\eta_j^2\lesssim \sigma^2\sum_{j=1}^d\eta_j^2\exp\left(-\frac{(\beta_*(s))_+}{|\eta_j|^\al}\right).$$

\textit{Term II.}

Note that the symmetry of $\xi_j$ implies that $\bbE_{\theta, {\caP}_\xi}\bigl(\eta_j\xi_j\bbo\left\{|\eta_j\xi_j|\ge \factorone\tau\lambda_*(s)\right\}\bigr)=0$. Using the independence across indices, we have
\begin{align*}
\bbE_{\theta, {\caP}_\xi}\biggl(\sum_{\substack{j>j_2(s)\\j\in S^c}}\sigma \eta_j\xi_j\bbo\left\{|\eta_j\xi_j|\ge \factorone\tau\lambda_*(s)\right\}\biggr)^2&=\sigma^2\sum_{\substack{j>j_2(s)\\j\in S^c}}\bbE \biggl(\eta_j^2\xi_j^2\bbo\left\{|\xi_j|\ge \tau\lambda_*(s)/|\eta_j|\right\}\biggr)\\
&\le \sigma^2\sum_{\substack{j>j_2(s)\\j\in S^c}}\sqrt{\bbE(\eta_j^4\xi_j^4)}\sqrt{\bbP(|\xi_j|\ge \factorone\tau\lambda_*(s)/|\eta_j|)}\\
&\le 2\sqrt{C_2^*}\sigma^2\sum_{j=1}^d\eta_j^2\exp\left(-\frac{(\beta_*)_+}{|\eta_j|^\al}\right).
\end{align*}

\textit{Term III.}
\begin{align*}
&\bbE_{\theta, {\caP}_\xi}\biggl(\sum_{\substack{j>j_2(s)\\j\in S}}\eta_jy_j\bbo\left\{|\eta_jy_j|\ge \factorone\sigma\tau\lambda_*(s)\right\} -\sum_{\substack{j>j_2(s)\\j\in S}}\eta_j\theta_j\biggr)^2\\
=&\bbE_{\theta, {\caP}_\xi}\biggl(\sigma\sum_{\substack{j>j_2(s)\\j\in S}} \eta_j\xi_j-\sum_{\substack{j>j_2(s)\\j\in S}}\eta_jy_j\bbo\left\{|\eta_jy_j|<\factorone\sigma\tau\lambda_*(s)\right\}\biggr)^2\\
\le& 2s\cdot\sum_{\substack{j>j_2(s)\\j\in S}}\left[\sigma^2\bbE(\eta_j\xi_j)^2+\bbE(\eta_jy_j\bbo\left\{|\eta_jy_j|<\factorone\sigma\tau\lambda_*(s)\right\})^2\right]\\
\le& 2s\cdot \sum_{\substack{j>j_2(s)\\j\in S}}\left[\sigma^2 \eta_j^2+\factorone\factorone\sigma^2\tau^2\lambda^2_*(s)\right]\le 2[1+\factorone\factorone\tau^2]\sigma^2\lambda^2_*(s) s^2.
\end{align*}
Combining the bounds for I, II, and III, we obtain the desired result in \Cref{thm: adaptive estimator}.
\end{proof}
\subsection{Proof of Theorem~\ref{thm: adaptive upper}}
\begin{proof}
For $\theta\in \Theta_s, {\caP}_\xi\in {\cal G}_{\al,\tau}^\otimes$, we have
$$\bbE_{\theta, {\caP}_\xi}(\tilde{L}_*-L(\theta))^2=\underbrace{\bbE_{\theta, {\caP}_\xi}\left[(\tilde{L}_*-L(\theta))^2\bbo\left\{\hat{s}\le s\right\}\right]}_{\text{ denoted as I}}+\underbrace{\bbE_{\theta, {\caP}_\xi}\left[(\tilde{L}_*-L(\theta))^2\bbo\left\{\hat{s}> s\right\}\right]}_{\text{ denoted as II}}.$$

\textit{Term I.} Recall that $s_0=s_*+1$. By definition of $\hat{s}$, we have $\hat{s}\le s_0$ and on the event $\{\hat{s}\le s\}$,
$$(\hat{L}_{\hat{s}}^*-L(\theta))^2\le 2\left[\omega_s^2+(\hat{L}_s^*-L(\theta))^2\right]\qquad \text{if } s<s_0 \text{ or }s\ge s_0>\hat{s}.$$
Thus,
\begin{align}
\forall s < s_0: \bbE_{\theta, {\caP}_\xi} \left[ (\hat{L}_{\hat{s}}^* - L(\theta))^2\bbo_{\hat{s} \le s} \right]
&\le 2\left[\zeta\sigma^2 \Phi_{\rm adp}(s;\eta) + \bbE_{\theta, {\caP}_\xi} \left( \hat{L}_s^* - L(\theta) \right)^2\right], \label{eq:adaptive_bound 1} \\
\forall s \ge s_0: \bbE_{\theta, {\caP}_\xi} \left[ (\hat{L}_{\hat{s}}^* - L(\theta))^2\bbo_{\hat{s} \le s} \right]
&\le \bbE_{\theta, {\caP}_\xi} \left[ (\hat{L}_{\hat{s}}^* - L(\theta))^2 \left(\bbo_{\hat{s} \le s, \hat{s} < s_0} +\bbo_{\hat{s} = s_0} \right) \right] \nonumber \\
&\le 2\left[\zeta\sigma^2\Phi_{\rm adp}(s;\eta) + \bbE_{\theta, {\caP}_\xi} \left( \hat{L}_s^* - L(\theta) \right)^2 \right] \label{eq:adaptive_bound 2}\\
&\qquad + \bbE_{\theta, {\caP}_\xi} \left( \hat{L}_{s_0}^* - L(\theta) \right)^2.\nonumber
\end{align}

For $s = 1, \ldots, s_0-1$, by Theorem~\ref{thm: adaptive estimator}, we have
\[
\bbE_{\theta, {\caP}_\xi} \left( \hat{L}_s^* - L(\theta) \right)^2 \le c_{1,0}\sigma^2\Phi_{\rm adp}(s;\eta)
\]
for some absolute constant \(c_{1,0}>0\).

For \(s\in [s_0, d]\), $\lambda_*(s)=0$ and there is no thresholding in the estimator \(\hat{L}_s^*\), and we have
\[
\bbE_{\theta, {\caP}_\xi} \left( \hat{L}_s^* - L(\theta) \right)^2
=\bbE_{\theta, {\caP}_\xi} \left( \hat{L}_{s_0}^* - L(\theta) \right)^2
=\sigma^2\sum_{j=1}^d\eta_j^2\le \sigma^2\Phi_{\rm adp}(s;\eta).
\]

Combining with \eqref{eq:adaptive_bound 1} and \eqref{eq:adaptive_bound 2}, we have
\[\bbE_{\theta, {\caP}_\xi} \left[\left( \hat{L}_{\hat{s}}^* - L(\theta) \right)^2\bbo\left\{\hat{s}\le s\right\}\right] \lesssim \sigma^2\Phi_{\rm adp}(s;\eta), \qquad s = 1, \ldots, d.\]

\textit{Term II.} Since by definition $\hat{s}\le s_0$, we have $s<\hat{s}\le s_0$ when the indicator inside the expectation equal 1.
In this case, \Cref{lem: probability bound,lem: fourth bound} can be applied to obtain the following for any $\theta\in\Theta_s$:
\begin{equation}\label{eq:adaptive decompose}
    \begin{aligned}
      &\sup_{\theta\in\Theta_{s}}\bbE_{\theta, {\caP}_\xi} \left[ (\hat{L}_{\hat{s}}^* - L(\theta))^2\bbo_{\hat{s} > s} \right]\\
      =&\sup_{\theta\in\Theta_{s}}\sum_{s<s^\prime\le s_0\wedge d} \bbE_{\theta, {\caP}_\xi} \left[(\hat{L}_{s^\prime}^*-L(\theta))^2\bbo\left\{\hat{s}=s^\prime\right\} \right]   \\
      \le& \sum_{s<s^\prime\le s_0\wedge d}\sup_{\theta\in\Theta_{s}}\sqrt{\bbE_{\theta, {\caP}_\xi}(\hat{L}_{s^\prime}^*-L(\theta))^4}\cdot \sup_{\theta\in\Theta_{s}}\sqrt{ \bbP_{\theta, {\caP}_\xi}\left(\hat{s}=s^\prime\right)}  \\
      \le& \sum_{s<s^\prime\le s_0\wedge d}\sup_{\theta\in \Theta_{s^\prime}}\sup_{{\caP}_\xi\in {\cal G}_{\al,\tau}^\otimes}\sqrt{\bbE_{\theta, {\caP}_\xi}(\hat{L}_{s^\prime}^*-L(\theta))^4}\sup_{\theta\in\Theta_{s^{\prime}-1}}\sup_{{\caP}_\xi\in {\cal G}_{\al,\tau}^\otimes} \sqrt{\bbP_{\theta, {\caP}_\xi}\left(\hat{s}>s^\prime-1\right)}\\
      \le& \sqrt{c_{3,0}\cdot c_{4,0}}\sigma^2\cdot \Phi_{\rm adp}(1;\eta)\sum_{s<s^\prime\le s_0\wedge d} (s^\prime)^{(4-7)/2}\log(es^\prime),
    \end{aligned}
\end{equation}
where the first inequality is due to H\"older's inequality,
and we have applied \Cref{lem: probability bound,lem: fourth bound} in the last line.
Therefore, we have $\bbE_{\theta, {\caP}_\xi} \left[ (\hat{L}_{\hat{s}}^* - L(\theta))^2\bbo_{\hat{s} > s} \right]\le c\cdot \sigma^2 \Phi_{\rm adp}(1;\eta)$
for some constant $c>0$.
The proof is completed since we have $\Phi_{\rm adp}(1;\eta)\le c_{2,0}\Phi_{\rm adp}(s;\eta)$ and $\sum_{i=1}^\infty i^{-3/2}\log(ei)<\infty$.
\end{proof}
\subsection{Proof of Theorem~\ref{thm: exponential upper}}
\begin{proof}Let $\theta\in \Theta_s$ and ${\caP}_\xi\in {\cal H}_{\al,\tau}^\otimes.$

  If $j_3(s) = d$, then no thresholding is applied in the estimator $\hat{L}_{\cal H}$, and $\xi_j$ are independent noise with mean zero and unit variance. Clearly, we have
  $$
  \bbE_{\theta,{\caP_\xi}}(\hat{L}_{\cal H}-L(\theta))^2
= \bbE_{\theta,{\caP_\xi}}\biggl(\sigma\sum_{j=1}^d\eta_j\xi_j\biggr)^2
= \sigma^2\sum_{j=1}^d\eta_j^2,
  $$
  and the desired result is proved.

 When $j_3(s)<d$, we have $s\le \sqrt{d}$.
  Let \(\varepsilon=(\varepsilon_1,\ldots,\varepsilon_d)\in \bbR^d\) with
  \[\varepsilon_j=y_j\bbo\left\{|y_j|\ge c_{\cal H}\cdot\sigma\cdot\log^{1/\al}\left( \frac{e d}{s} \right)\right\}-\theta_j.\]
  Similar to \eqref{eq:nonadaptive error decomposition}, we have
  \begin{equation}\label{eq:exponential error decomposition}
    \bbE_{\theta, {\caP_\xi}}(\hat{L}_{\cal H}-L(\theta))^2\le 2\biggl[\underbrace{\bbE_{\theta, {\caP_\xi}}\biggl(\sum_{j\le j_3(s)}\eta_j y_j-\sum_{j\le j_3(s)}\eta_j\theta_j\biggr)^2}_{\text{denoted as I}}+\underbrace{\bbE_{\theta, {\caP_\xi}}\biggl(\sum_{j>j_3(s)}\eta_j\varepsilon_j\biggr)^2}_{\text{denoted as II}}\biggr].
  \end{equation}

\textit{Term I.} $$\bbE_{\theta, {\caP_\xi}}\biggl(\sum_{j\le j_3(s)}\eta_j y_j-\sum_{j\le j_3(s)}\eta_j\theta_j\biggr)^2=\sigma^2\sum_{j\le j_3(s)}\eta_j^2.$$

\textit{Term II.}
  \begin{align*}
    \bbE_{\theta, {\caP_\xi}}\biggl(\sum_{j>j_3(s)}\eta_j\varepsilon_j\biggr)^2\le |\eta_{j_3(s)+1}|^2\cdot \bbE_{\theta, {\caP_\xi}}\|\varepsilon\|_1^2.
  \end{align*}
  For \(j\in S^c\), we have
  \begin{align*}
    \varepsilon_j&=\sigma \xi_j\bbo\left\{|\xi_j|\ge c_{\cal H}\cdot\log^{1/\al}\left( \frac{e d}{s} \right)\right\}.
  \end{align*}
  Since \[\bbP(|\xi_j|\ge t)\le 2\exp\left\{-\left(\frac{t}{\tau}\right)^\al\right\}\qquad\text{for any } t\ge 0,\]
  we have
  \begin{align*}
    \bbE \varepsilon_j^2&= \sigma^2 \left(\bbE \xi_j^2\bbo\left\{|\xi_j|\ge c_{\cal H}\log^{1/\al}\left( \frac{e d}{s} \right)\right\}\right)\\
    &\le \sigma^2\sqrt{\bbE\xi_j^4}\cdot\sqrt{\bbP\left(|\xi_j|\ge c_{\cal H}\cdot\log^{1/\al}\left( \frac{e d}{s} \right)\right)}\\
    &\le \sigma^2\sqrt{2 C_2^* \left(\frac{s}{ed}\right)^{(c_{\cal H}/{\tau})^\al}},
  \end{align*}
  where the last inequality comes from \Cref{lem: tail}.
Let $C_1=\sqrt{2C_2^*}/e^2$ and $c_{\cal H}= \tau 4^{1/\al}$.
We have

  \[
  \begin{aligned}
\bbE_{\theta, \caP_\xi}\left\|\varepsilon_{S^c}\right\|_1^2 & \le |S^c|\biggl(\sum_{j\in S^c}\bbE\varepsilon_j^2\biggr) \\
&  \le  \sigma^2d^2\sqrt{2 C_2^* \left(\frac{s}{ed}\right)^{ 4 }}\\
& \le C_1 \sigma^2 s^2.
  \end{aligned}
  \]
  For \(j\in S\), we have
  \begin{align*}
    |\varepsilon_j|=\left|\sigma \xi_j- y_j\bbo\left\{|y_j|<c_{\cal H}\cdot\sigma\cdot\log^{1/\al}\left( \frac{e d}{s} \right)\right\}\right|\le \sigma |\xi_j|+c_{\cal H}\cdot\sigma\cdot\log^{1/\al}\left( \frac{e d}{s} \right).
  \end{align*}
  The first term on the right hand side has a bounded second moment, and the number of nonzero coordinates is at most \(s\). Therefore, we have $\bbE_{\theta, \caP_\xi}\left\|\varepsilon_S\right\|_1^2\le C_2 \sigma^2 s^2 \log^{2/\al}\left(ed/s\right)$ for some constant \(C_2>0\) depending on $\al$ and $\tau$.

  Combining the above analysis, we have
  \[\bbE_{\theta, \caP_\xi}\left(\hat{L}_{\cal H}-L(\theta)\right)^2\le \sigma^2 \left[\sum_{j\le j_3(s)}\eta_j^2+2(C_1+C_2)|\eta_{j_3(s)+1}|^2 s^2\log^{2/\al}\left(\frac{ed}{s}\right)\right] .\]
  Since $s^2\log^{2/\al}(ed/s)\le j_3(s)<d$ by
\eqref{eq:index unknown} and $|\eta_j|$ is nonincreasing, we have

  \[|\eta_{j_3(s)+1}|^2 s^2\log^{2/\al}\left(\frac{ed}{s}\right)\le \sum_{j\le j_3(s)}\eta_j^2.\]
  Therefore, the proof of Theorem~\ref{thm: exponential upper} is completed.

\end{proof}
\subsection{Proof of \Cref{thm: unknown sigma}}
\begin{proof}
Let $\theta \in \Theta_s$ and ${\caP}_\xi \in {\cal G}_{\al,\tau}^\otimes$.
As before, we only need to consider the case where $\lambda_{\rm o} > 0$.
Similar to \eqref{eq:nonadaptive error decomposition}, we decompose the risk as
\begin{align*}
    \mathbb{E}_{\theta,{\caP}_\xi}\!\left(\hat{L}_s^\prime - L(\theta)\right)^2
    \;\le\; 3\Biggl[
    &\underbrace{\mathbb{E}_{\theta,{\caP}_\xi}\!\left(\sum_{j \le j_1(s)} \eta_j y_j - \sum_{j \le j_1(s)} \eta_j \theta_j\right)^2}_{\text{denoted as I}} \\[1ex]
    &\quad+ \underbrace{\mathbb{E}_{\theta,{\caP}_\xi}\!\left(\sum_{\substack{j > j_1(s) \\ j \in S^c}} \sigma \eta_j \xi_j \,\mathbbm{1}\!\left\{\sigma|\eta_j \xi_j| \ge \factorone\sqrt{2}\hat{\sigma}\tau \lambda_o\right\}\right)^2}_{\text{denoted as II}} \\[1ex]
    &\quad+ \underbrace{\mathbb{E}_{\theta,{\caP}_\xi}\!\left(\sum_{\substack{j > j_1(s) \\ j \in S}} \eta_j y_j \,\mathbbm{1}\!\left\{|\eta_j y_j| \ge \factorone\sqrt{2}\hat{\sigma}\tau \lambda_o\right\}
    - \sum_{\substack{j > j_1(s) \\ j \in S}} \eta_j \theta_j\right)^2}_{\text{denoted as III}}
    \Biggr].
\end{align*}

\textit{Term I.}
\[
\mathbb{E}_{\theta,{\caP}_\xi}\!\left(\sum_{j \le j_1(s)} \eta_j y_j - \sum_{j \le j_1(s)} \eta_j \theta_j\right)^2
= \sigma^2 \sum_{j \le j_1(s)} \eta_j^2
\;\lesssim\; \sigma^2 \Phi_{\rm o}(s;\eta).
\]

\textit{Term II.}

For $j\in S^c$, let
$$
H_j=\sigma\eta_j\xi_j
\mathbbm{1}\!\left\{
\sigma|\eta_j\xi_j|>
\factorone\sqrt{2}\hat{\sigma}\tau\lambda_o
\right\}.
$$
For distinct $j,k\in S^c$, changing $\xi_j$ to $-\xi_j$
leaves both the joint distribution of the noise variables and
$\hat{\sigma}$ unchanged, because $\theta_j=0$ and $\hat{\sigma}$
is a function of $(y_1^2,\ldots,y_d^2)$. Since this change of sign also changes the sign of the product $H_jH_k$, this product is symmetric around 0 and hence $\mathbb EH_jH_k=0$. Therefore, we have

\begin{align*}
    &\mathbb{E}_{\theta,{\caP}_\xi}\!\left(\sum_{\substack{j > j_1(s) \\ j \in S^c}} \sigma \eta_j \xi_j \,\mathbbm{1}\!\left\{\sigma|\eta_j \xi_j| \ge \factorone\sqrt{2}\hat{\sigma}\tau \lambda_o\right\}\right)^2 \\
    &= \sum_{\substack{j > j_1(s) \\ j \in S^c}}
    \mathbb{E}\!\left[ \sigma^2 \eta_j^2 \xi_j^2 \,\mathbbm{1}\!\left\{\sigma|\eta_j \xi_j| \ge \factorone\sqrt{2}\hat{\sigma}\tau \lambda_o\right\}\right] \\
    &= \sum_{\substack{j > j_1(s) \\ j \in S^c}}
    \mathbb{E}\!\left[\sigma^2 \eta_j^2 \xi_j^2 \,\mathbbm{1}\!\left\{\sigma|\eta_j \xi_j|\ge \factorone\sqrt{2}\hat{\sigma}\tau \lambda_o,\, \hat{\sigma} < \frac{\sigma}{\sqrt{2}}\right\}\right] \\
    &\quad+ \sum_{\substack{j > j_1(s) \\ j \in S^c}}
    \mathbb{E}\!\left[\sigma^2 \eta_j^2 \xi_j^2 \,\mathbbm{1}\!\left\{\sigma|\eta_j \xi_j| \ge \factorone\sqrt{2}\hat{\sigma}\tau \lambda_o,\, \hat{\sigma} \ge \frac{\sigma}{\sqrt{2}}\right\}\right].
\end{align*}
For the first term, we apply H\"older's inequality and \Cref{lem: mom} to obtain
\begin{align*}
    \sum_{\substack{j > j_1(s) \\ j \in S^c}}
    \sigma^2 \mathbb{E}\!\left[\eta_j^2 \xi_j^2 \,\mathbbm{1}\!\left\{\hat{\sigma} < \frac{\sigma}{\sqrt{2}}\right\}\right]
    &\le \sigma^2 \sum_{j=1}^d \sqrt{ \mathbb{E}(\eta_j^4 \xi_j^4) \cdot \mathbb{P}(\hat{\sigma} < \frac{\sigma}{\sqrt{2}})} \\
    &\le \sigma^2 \exp\!\left(-\frac{c_{\al,\tau}^*}{2}d\right) \sum_{j=1}^d \eta_j^2\\
    & \le \sigma^2 d \exp\!\left(-\frac{c_{\al,\tau}^*}{2}d\right) \eta_1^2 \\
& \lesssim\; \sigma^2 \Phi_{\rm o}(s;\eta),
\end{align*}
where the last inequality is because $d \exp\!\left(-\tfrac{c_{\al,\tau}^*}{2}d\right)\lesssim 1$  and $\eta_1^2 \lesssim \Phi_{\rm o}(s;\eta)$ by the definition of $\Phi_{\rm o}(s;\eta)$ (either $\eta_1^2\le \lambda_o^2$ or  $\lambda_o(s)/|\eta_1|\le 1$, which implies that $\nu^2\ge \eta_1^2 e^{-1}$).

For the second term, we first note that
$\mathbbm{1}\!\left\{\sigma|\eta_j \xi_j| \ge \factorone\sqrt{2}\hat{\sigma}\tau \lambda_o,\, \hat{\sigma} \ge \frac{\sigma}{\sqrt{2}}\right\}\le \mathbbm{1}\!\left\{\sigma|\eta_j \xi_j| \ge \factorone \sigma\tau \lambda_o\right\}$.
Then we can bound the sum as in the proof of \Cref{thm: nonadaptive linear}:
\begin{align*}
    \sum_{\substack{j > j_1(s) \\ j \in S^c}}
    \mathbb{E}\!\left[\sigma^2 \eta_j^2 \xi_j^2 \,\mathbbm{1}\!\left\{|\eta_j \xi_j| \ge \factorone\tau \lambda_o\right\}\right]
    &\le \sigma^2 \sum_{j=1}^d \eta_j^2 \exp\!\left(-\frac{\lambda_{\rm o}^\al}{|\eta_j|^\al}\right) \lesssim \sigma^2 \Phi_{\rm o}(s;\eta).
\end{align*}

\textit{Term III.}
\Cref{lem: mom} implies that
$\mathbb{E}\hat{\sigma}^2\le \sigma^2(1+C_{\al, \tau}^*)$.
Let
$$
A=\{j>j_1(s):j\in S\}.
$$
We have $|A|\le s$ and
\begin{align*}
    &\mathbb{E}_{\theta,{\caP}_\xi}\!\left(\sum_{j\in A} \eta_j y_j \,\mathbbm{1}\!\left\{|\eta_j y_j| \ge \factorone\sqrt{2}\hat{\sigma}\tau \lambda_o\right\}
    - \sum_{j\in A} \eta_j \theta_j\right)^2 \\
    =& \mathbb{E}_{\theta,{\caP}_\xi}\!\left(\sum_{j\in A} \eta_j y_j \,\mathbbm{1}\!\left\{|\eta_j y_j| < \factorone\sqrt{2}\hat{\sigma}\tau \lambda_o\right\}
    - \sum_{j\in A} \sigma \eta_j \xi_j\right)^2 \\
    \le &2s \sum_{j\in A} \mathbb{E}\!\left[\eta_j^2 y_j^2 \,\mathbbm{1}\!\left\{|\eta_j y_j| < \factorone\sqrt{2}\hat{\sigma}\tau \lambda_o\right\} + \sigma^2 \eta_j^2 \xi_j^2 \right] \\
     \le &2s \sum_{j\in A} \left( \factorone\factorone2^{}\tau^2 \lambda_o^2\mathbb{E}\hat{\sigma}^2+ \sigma^2 \eta_j^2  \right) \\
         \le & \factorone\factorone2^{2}s^2 \tau^2 \lambda_o^2 \sigma^2(1+C_{\al, \tau}^*) + 2 s \sigma^2 \sum_{j\in A} \lambda_o^2  \\
\lesssim &
     \sigma^2 \Phi_{\rm o}(s;\eta).
\end{align*}

Combining the bounds for I, II, and III, we obtain the desired result in \Cref{thm: unknown sigma}.
\end{proof}
\subsection{Proof of the upper bound in \Cref{thm: testing}}\label{sec:test upper}
Note that for the test $\Delta_s$ defined in \eqref{eq:test}, we have
$$\forall \theta\in \Theta_{s,0},\,\bbP_{\theta, {\caN}^\otimes}(\Delta_s=1)=\bbP_{\theta,{\caN}^\otimes}\left(|\hat{L}_s-t_0|\ge B\sigma \sqrt{\Phi_{\rm o}(s;\eta)}\right)\le \frac{\bbE_{\theta, {\caN}^\otimes}(\hat{L}_s-L(\theta))^2}{B^2\sigma^2\Phi_{\rm o}(s;\eta)}.$$
For any $A>B$ and any $ \theta\in \Theta_s(A\sigma \sqrt{\Phi_{\rm o}(s;\eta)})$, we have $|L(\theta)-t_0|\ge A\sigma \sqrt{\Phi_{\rm o}(s;\eta)}$. Therefore,
\begin{align*}
\bbP_{\theta, {\caN}^\otimes}(\Delta_s=0)&=\bbP_{\theta,{\caN}^\otimes}(|\hat{L}_s-t_0|< B\sigma \sqrt{\Phi_{\rm o}(s;\eta)})\\
&\le \bbP_{\theta,{\caN}^\otimes}(|\hat{L}_s-L(\theta)| > (A-B)\sigma \sqrt{\Phi_{\rm o}(s;\eta)})\\
&\le \frac{\bbE_{\theta, {\caN}^\otimes}(\hat{L}_s-L(\theta))^2}{(A-B)^2\sigma^2\Phi_{\rm o}(s;\eta)}.
\end{align*}
By \Cref{thm: nonadaptive linear}, there is some constant $C>0$ such that
$$\sup_{\theta\in \Theta_s}\bbE_{\theta, {\caN}^\otimes}(\hat{L}_s-L(\theta))^2\le C \sigma^2\Phi_{\rm o}(s;\eta).$$
For any $\varepsilon>0$, we can choose $B$ and $A_\varepsilon$ large enough such that
$C/B^2<\varepsilon/2$ and $C/(A-B)^2<\varepsilon/2$, which implies that
$$\forall A\ge A_\varepsilon, \,\sup_{\theta\in \Theta_{s,0}}\bbP_{\theta, {\caN}^\otimes}(\Delta_s=1) \;+\; \sup_{\theta \in \Theta_s(A\sigma\sqrt{\Phi_{\rm o}(s;\eta)})} \bbP_{\theta, {\caN}^\otimes}(\Delta_s=0)\le \varepsilon.
$$
This completes the proof.

\subsection{Proofs of the Lemmas}\label{sec:technical lemmas upper}
\begin{proof}[Proof of Lemma~\ref{lem: probability bound}]
  Let $\theta\in \Theta_s$ and ${\caP}_\xi\in {\cal G}_{\al,\tau}^\otimes$.
  Below, we drop the subscript in $\mathbb{P}_{\theta,\caP_\xi}$.

  Let $m=s_0\wedge d$.
  By the definition of the estimator $\hat L_s^*$ in \eqref{eq:nonadaptive estimator}, the following identity holds for any pair $s,s^\prime$ satisfying $s<s^\prime\le m$:

  \begin{align*}
    \hat{L}_{s^\prime}^*- \hat{L}_s^*=&\sum_{\substack{j_2(s)<j\le j_2(s^\prime)\\j\in S}}\eta_jy_j\bbo\left\{|\eta_jy_j|\le \factorone\sigma\tau\lambda_*(s)\right\}\\
    &\quad +\sum_{\substack{j>j_2(s^\prime)\\j\in S}}\eta_jy_j\bbo\left\{\factorone\sigma\tau\lambda_*(s^\prime)<|\eta_jy_j|\le \factorone\sigma\tau \lambda_*(s)\right\}\\
     &\quad\quad+\sigma\sum_{\substack{j_2(s)<j\le j_2(s^\prime)\\j\notin S}}\eta_j\xi_j\bbo\left\{|\eta_j\xi_j|\le  \factorone\tau\lambda_*(s)\right\}\\
    &\quad\quad\quad +\sigma\sum_{\substack{j>j_2(s^\prime)\\j\notin S}}\eta_j\xi_j\bbo\left\{ \factorone\tau\lambda_*(s^\prime)<|\eta_j\xi_j|\le \factorone\tau\lambda_*(s)\right\}.
  \end{align*}
  The sum of the first two terms on the right hand side is bounded by
  $$ \factorone\tau\sigma s\lambda_*(s)\le \factorone\tau\sigma\sqrt{\Phi_{\rm adp}(s;\eta)}\le \sqrt{c_{2,0}}\tau\sigma\sqrt{\Phi_{\rm adp}(s^\prime;\eta)}, $$
  where the second inequality is due to \Cref{prop:linear decrease}.
  Therefore, once $\zeta>9 c_{2,0}\tau^2$, we have
  \begin{align*}
    \biggl|&\sum_{\substack{j_2(s)<j\le j_2(s^\prime)\\j\in S}}\eta_jy_j\bbo\left\{|\eta_jy_j|\le \factorone\sigma\tau\lambda_*(s)\right\}\\
    &\quad+\sum_{\substack{j>j_2(s^\prime)\\j\in S}}\eta_jy_j\bbo\left\{\sigma\tau\lambda_*(s^\prime)<|\eta_jy_j|\le \factorone\sigma\tau \lambda_*(s)\right\}\biggr|\le \frac{1}{3}\omega_{s^\prime},
  \end{align*}
  which implies that
  \begin{align*}
    &\bbP_{\theta}\left(|\hat{L}_{s^\prime}^*-\hat{L}_s^*|\ge \omega_{s^\prime}\right)\\
    \le &\underbrace{\bbP\biggl(\biggl|\sum_{\substack{j_2(s)<j\le j_2(s^\prime)\\j\notin S}}\eta_j\xi_j\bbo\left\{|\eta_j\xi_j|\le  \factorone\tau\lambda_*(s)\right\}\biggr|\ge \frac{\omega_{s^\prime}}{3\sigma}\biggr)}_{\text{denoted as I}}\\
    &\,+\underbrace{\bbP_{\theta}\biggl(\biggl|\sum_{\substack{j>j_2(s^\prime)\\j\notin S}}\eta_j\xi_j\bbo\left\{ \factorone\tau\lambda_*(s^\prime)<|\eta_j\xi_j|\le  \factorone\tau\lambda_*(s)\right\}\biggr|\ge \frac{\omega_{s^\prime}}{3\sigma}\biggr)}_{\text{denoted as II}}.
  \end{align*}

\medskip
  \textit{Term I.}
  The random variables $\eta_j\xi_j\bbo\left\{|\eta_j\xi_j|\le \factorone\tau \lambda_*(s)\right\}$ are independent, mean zero, and bounded by $K=\factorone\tau\lambda_*(s)$. Let $U_1=\{j_2(s)<j\le j_2(s^\prime): j\notin S\}$. By Bernstein's inequality for bounded distributions \cite[Theorem 2.8.4]{vershynin2018high}, we have
  $$\bbP\left\{\biggl|\sum_{j\in U_1}\eta_j\xi_j\bbo\left\{|\eta_j\xi_j|\le \tau \lambda_*(s)\right\}\biggr|\ge t\right\}\le 2\exp\left(-\frac{t^2}{\tilde{\sigma}^2+Kt/3}\right),\forall t>0,$$
  where
  \[\tilde{\sigma}^2=\sum_{j\in U_1}\bbE \eta_j^2\xi_j^2\bbo\left\{|\eta_j\xi_j|\le \factorone\tau \lambda_*(s)\right\}\le \sum_{j\le j_2(s^\prime)}\eta_j^2.\]

  For $t = \omega_{s^\prime}/(3\sigma)$, we have
\[
t^2 = \frac{\zeta}{9}\Phi_{\rm adp}(s^\prime;\eta)
\;\ge\; \frac{\log(es^\prime)}{9e}\sum_{j \le j_2(s^\prime)} \eta_j^2,
\]
and
\[
\frac{t}{K} \;\ge\; \frac{\sqrt{\zeta}\,\sqrt{\Phi_{\rm adp}(1;\eta)}\,\log(es^\prime)}{\factorone\tau\lambda_*(s)}.
\]
Since $\sqrt{\Phi_{\rm adp}(1;\eta)} \gtrsim \lambda_*(1) \gtrsim \lambda_*(s)$, we may choose $\zeta$ sufficiently large such that
\[
\mathbb{P}\!\left(\left|\sum_{\substack{j_2(s)<j\le j_2(s^\prime)\\ j\notin S}}
   \eta_j \xi_j \bbo\{|\eta_j\xi_j|\le \factorone\tau\lambda_*(s)\}\right|
   \;\ge\; \frac{\sqrt{\zeta}}{3}\sqrt{\Phi_{\rm adp}(s^\prime;\eta)}\right)
   \;\le\; 2(s^\prime)^{-8}.
\]
Moreover, when $\al \ge 2$, the truncated variables
$\xi_j \bbo\{|\eta_j\xi_j|\le \tau\lambda_*(s)\}$ are sub-Gaussian, and thus Hoeffding's inequality \cite[Theorem 2.6.3]{vershynin2018high} yields the same probability bound.

\medskip
\textit{Term II.} Let $U_2=\{j>j_2(s^\prime):j\notin S\}$. If
$\lambda_*(s^\prime)=0$, then $j_2(s^\prime)=d$ and $U_2$ is empty, so
$\mathrm{II}=0$. We therefore assume that $\lambda_*(s^\prime)>0$. Set
$$
a=\factorone\tau\lambda_*(s^\prime),\qquad
b=\factorone\tau\lambda_*(s).
$$
Since $s<s^\prime$ and $s\mapsto\lambda_*(s)$ is nonincreasing, we have
$0<a\le b$.
By the symmetry and independence of the noise variables, we may write $\xi_j=\varepsilon_jR_j$, where $R_j$ has the
same distribution as $|\xi_j|$ and the variables $(\varepsilon_j:j\in U_2)$
are independent Rademacher variables independent of $(R_j:j\in U_2)$. Define
$$
\overline Z_j=|\eta_j|R_j\bbo\{a<|\eta_j|R_j\le b\},\qquad
Q=\sum_{j\in U_2}\overline Z_j^2, \qquad
\mathcal R=\sigma(R_j:j\in U_2).
$$
Conditional on $\mathcal R$, Hoeffding's inequality for Rademacher sums gives
\begin{align*}
q_{\mathcal R}
&:=\bbP\left(\left.\left|\sum_{j\in U_2}\eta_j\varepsilon_jR_j
\bbo\{a<|\eta_j|R_j\le b\}\right|
\ge \frac{\omega_{s^\prime}}{3\sigma}\,\right|\,\mathcal R\right)\\
&\le 2\exp\left(-\frac{\zeta\Phi_{\rm adp}(s^\prime;\eta)}{18Q}\right),
\end{align*}
where the right hand side is interpreted as zero when $Q=0$.
For every $\Delta>0$, splitting according to the events $\{Q\le\Delta\}$
and $\{Q>\Delta\}$ gives
\begin{equation}\label{eq:probability bound 4}
\mathrm{II}=\bbE q_{\mathcal R}=\bbE q_{\mathcal R}\bbo\{Q\le\Delta\} + \bbE q_{\mathcal R}\bbo\{Q>\Delta\}
\le 2\exp\left(-\frac{\zeta\Phi_{\rm adp}(s^\prime;\eta)}{18\Delta}\right)
+\bbP(Q>\Delta).
\end{equation}
Choose
$$
\Delta=\frac{\zeta\Phi_{\rm adp}(s^\prime;\eta)}{144\log(es^\prime)}.
$$
The first term on the right hand side of \eqref{eq:probability bound 4} is then
equal to $2(es^\prime)^{-8}$.
We next bound $\bbP(Q>\Delta)$, where we will use the following Fuk-Nagaev inequality \citep[page 78]{petrov1995limit}.

\begin{lemma}[Fuk--Nagaev inequality]\label{lem: fuk-nagaev}
Let $p>2$ and $v>0$. Assume that $X_1,\ldots,X_n$ are independent random
variables such that $\bbE X_i=0$ and $\bbE|X_i|^p<\infty$ for
$i=1,\ldots,n$. Then
$$
\bbP\left(\sum_{i=1}^nX_i>v\right)
\le \left(1+\frac{2}{p}\right)^p v^{-p}
\sum_{i=1}^n\bbE|X_i|^p
+\exp\left(-\frac{2v^2}{(p+2)^2e^p\sum_{i=1}^n\bbE X_i^2}\right).
$$
\end{lemma}
For every integer $p\ge1$,
Lemma~\ref{lem: tail} and the Cauchy--Schwarz inequality give, for $j\in U_2$ and
$\eta_j\neq0$,
\begin{equation}\label{eq:lem-probability-bound-II-CS}
    \begin{aligned}
\bbE\overline Z_j^{2p}
&\le |\eta_j|^{2p}\bbE\left[|\xi_j|^{2p}
\bbo\left\{|\xi_j|>\frac{\factorone\tau\lambda_*(s^\prime)}{|\eta_j|}\right\}\right]\\
&\le |\eta_j|^{2p}\sqrt{\bbE|\xi_j|^{4p}}
\sqrt{\bbP\left(|\xi_j|>\frac{\factorone\tau\lambda_*(s^\prime)}{|\eta_j|}\right)}\\
&\le \sqrt{2C_{2p}^*}\,|\eta_j|^{2p}
\exp\left(-\frac{(\beta_*(s^\prime))_+}{|\eta_j|^\al}\right).
\end{aligned}
\end{equation}

Let
$$
D_{s^\prime}^2:=\sum_{j=1}^d\eta_j^2
\exp\left(-\frac{(\beta_*(s^\prime))_+}{|\eta_j|^\al}\right).
$$
For $j\in U_2$, let $X_j=\overline Z_j^2-\bbE\overline Z_j^2$. These random
variables are independent and centered.
Since $\bbE|X_j|^p\le2^p\bbE\overline Z_j^{2p}$ and
$|\eta_j|<\lambda_*(s^\prime)$ for $j\in U_2$, we have
\begin{align*}
\sum_{j\in U_2}\bbE|X_j|^p
&\le C_{p,\al,\tau}[\lambda_*(s^\prime)]^{2p-2}
D_{s^\prime}^2,
\end{align*}
where $C_{p,\al,\tau}=2^p\sqrt{2C_{2p}^*}$.
Moreover,
$$
\Phi_{\rm adp}(s^\prime;\eta)\ge\Phi_*(s^\prime;\eta)
=(s^\prime)^2\lambda_*^2(s^\prime)
+\log(es^\prime)D_{s^\prime}^2,
$$
and hence
\begin{align*}
[\lambda_*(s^\prime)]^{2p-2}D_{s^\prime}^2
&=\frac{[(s^\prime)^2\lambda_*^2(s^\prime)]^{p-1}
[\log(es^\prime)D_{s^\prime}^2]}
{(s^\prime)^{2p-2}\log(es^\prime)}\\
&\le\frac{[\Phi_{\rm adp}(s^\prime;\eta)]^p}
{(s^\prime)^{2p-2}\log(es^\prime)}.
\end{align*}
Consequently,
\begin{equation}\label{eq:annulus pth moments}
\sum_{j\in U_2}\bbE|X_j|^p
\le \frac{C_{p,\al,\tau}}{(s^\prime)^{2p-2}\log(es^\prime)}
[\Phi_{\rm adp}(s^\prime;\eta)]^p.
\end{equation}
Taking $p=2$ gives
\begin{equation}\label{eq:annulus centered moments}
\sum_{j\in U_2}\bbE X_j^2
\le \frac{C_{2,\al,\tau}}{(s^\prime)^2\log(es^\prime)}
[\Phi_{\rm adp}(s^\prime;\eta)]^2.
\end{equation}
\Cref{eq:lem-probability-bound-II-CS} with $p=1$ gives
\begin{equation}\label{eq:annulus second moments}
\sum_{j\in U_2}\bbE\overline Z_j^2
\le \sqrt{2C_{2}^*} D_{s^\prime}^2
\le \frac{\sqrt{2C_{2}^*} }{\log(es^\prime)}
\Phi_{\rm adp}(s^\prime;\eta).
\end{equation}
Choose the constant $\zeta$ sufficiently large, depending only on $\al$ and $\tau$, such that $\zeta/288\ge \sqrt{2C_{2}^*}$.
By definition of $\Delta$ and
\eqref{eq:annulus second moments},
$$
v:=\Delta-\sum_{j\in U_2}\bbE\overline Z_j^2
\ge\frac{\Delta}{2}
=\frac{\zeta\Phi_{\rm adp}(s^\prime;\eta)}{288\log(es^\prime)}.
$$

We now apply Lemma~\ref{lem: fuk-nagaev}. First observe that
$$
Q=\sum_{j\in U_2}X_j
+\sum_{j\in U_2}\bbE\overline Z_j^2.
$$
Thus, by the definition of $v$,
$$
\{Q>\Delta\}
=
\left\{\sum_{j\in U_2}X_j>v\right\}.
$$
If $\sum_{j\in U_2}\bbE X_j^2=0$, then
$\sum_{j\in U_2}X_j=0$ almost surely and the desired probability is zero.
We may therefore assume that
$\sum_{j\in U_2}\bbE X_j^2>0$.

For notational convenience, let
$L=\log(es^\prime)$ and
$\Phi=\Phi_{\rm adp}(s^\prime;\eta)$.
Applying Lemma~\ref{lem: fuk-nagaev} with $p=6$ gives
$$
\bbP(Q>\Delta)\le T_1+T_2,
$$
where
\begin{align*}
T_1
&=
\left(1+\frac{2}{6}\right)^6v^{-6}
\sum_{j\in U_2}\bbE|X_j|^6,\\
T_2
&=
\exp\left(
-\frac{2v^2}
{8^2e^6\sum_{j\in U_2}\bbE X_j^2}
\right).
\end{align*}

\textit{Bounding $T_1$.}
By the lower bound
$v\ge\zeta\Phi/(288L)$ and
\eqref{eq:annulus pth moments} with $p=6$,
\begin{align*}
T_1
&\le
\left(\frac{4}{3}\right)^6
\left(\frac{288L}{\zeta\Phi}\right)^6
\frac{C_{6,\al,\tau}\Phi^6}
{(s^\prime)^{10}L}\\
&=
\left(\frac{4}{3}\right)^6
288^6C_{6,\al,\tau}\zeta^{-6}
\frac{L^5}{(s^\prime)^{10}}.
\end{align*}
Since $s<s^\prime$ are positive integers,
$s^\prime\ge2$. The function
$[\log(ex)]^5/x^2$ is continuous on $[2,\infty)$
and converges to zero as $x\to\infty$, so $\frac{L^5}{(s^\prime)^{2}}$ is bounded.
Therefore,
$$
T_1
\le
C_{\al,\tau}\zeta^{-6}(s^\prime)^{-8}
$$
for some constant $C_{\al,\tau}$.

\textit{Bounding $T_2$.}
By the same lower bound on $v$ and
\eqref{eq:annulus centered moments},
\begin{align*}
T_2
&\le
\exp\left\{
-\frac{2[\zeta\Phi/(288L)]^2}
{8^2e^6C_{2,\al,\tau}\Phi^2/[(s^\prime)^2L]}
\right\}\\
&=
\exp\left(
-c_{\al,\tau}\zeta^2
\frac{(s^\prime)^2}{L}
\right),
\end{align*}
where
$
c_{\al,\tau}
=(32\cdot288^2e^6C_{2,\al,\tau})^{-1}$.
For $s^\prime\ge2$, the inequalities
$\log(s^\prime)\le s^\prime/2$ and
$L\le s^\prime$ imply that
$$
\frac{(s^\prime)^2}{L}
\ge
2\log(s^\prime).
$$
We can therefore increase $\zeta$, depending only on
$\al$ and $\tau$, so that
$2c_{\al,\tau}\zeta^2\ge8$. With this choice,
$$
T_2\le(s^\prime)^{-8}.
$$

Combining the bounds for $T_1$ and $T_2$ gives
$$
\bbP(Q>\Delta)
\le
(C_{\al,\tau}+1)(s^\prime)^{-8}.
$$
Substituting this estimate into
\eqref{eq:probability bound 4}, and recalling that its first term equals
$2(es^\prime)^{-8}$, we conclude that for some constant $C_{\al,\tau}^\prime$,
$$
\mathrm{II}
\le
C_{\al,\tau}^\prime(s^\prime)^{-8}.
$$

\bigskip

Combining the bounds for Terms I and II, there exists a constant
$c_{\al,\tau}>0$, depending only on $\al$ and $\tau$, such that
\begin{equation}\label{eq:lem-probability-bound-pairwise}
    \bbP\left(
|\hat L_s^*-\hat L_{s^\prime}^*|\ge\omega_{s^\prime}
\right)
\le c_{\al,\tau}(s^\prime)^{-8}
\end{equation}

for all integers $s<s^\prime\le m$, uniformly over
$\theta\in\Theta_s$ and
$\mathcal P_\xi\in\mathcal G_{\al,\tau}^{\otimes}$.

We next reduce the comparisons in the definition of $\hat s$ to this range.

If $s_0\le d$, then by definition of $s_0$, we have $\lambda_*(r)=0$ for every $r\in[d]$ with
$r\ge s_0$. Hence $j_2(r)=d$ and
$$
\hat L_r^*=\sum_{j=1}^d\eta_jy_j=\hat L_{s_0}^*,
\qquad r\in[d],\quad r\ge s_0
$$
Moreover, the definitions of $\Phi_*$ and $\Phi_{\rm adp}$ give
$\Phi_{\rm adp}(r;\eta)=\Phi_{\rm adp}(s_0;\eta)$ and therefore
$\omega_r=\omega_{s_0}$ for every $r\in[d]$ with $r\ge s_0$.
Thus all comparisons indexed by $r\in[d]$ with $r>s_0$ coincide with
the comparison at $r=s_0=m$.

If $s_0>d$, then $m=d$ and there are no
indices in $[d]$ larger than $m$.

In both cases, we have
\begin{equation}\label{eq:comparison-event-saturation}
\bigcup_{\substack{r\in[d]\\r>s}}
\left\{|\hat L_s^*-\hat L_r^*|>\omega_r\right\}
=
\bigcup_{s<r\le m}
\left\{|\hat L_s^*-\hat L_r^*|>\omega_r\right\}.
\end{equation}

Recall the definition of $\hat s$ in \eqref{eq:selected index}. If
$|\hat L_s^*-\hat L_r^*|\le\omega_r$ for every $r\in[d]$ with $r>s$,
then $s$ is an admissible index in the defining minimum and hence
$\hat s\le s$. Taking the contrapositive and using
\eqref{eq:comparison-event-saturation}, we obtain
$$
\{\hat s>s\}
\subseteq
\bigcup_{s<r\le m}
\left\{|\hat L_s^*-\hat L_r^*|>\omega_r\right\}.
$$
Using the union bound and \Cref{eq:lem-probability-bound-pairwise}, we have
\begin{align*}
\bbP(\hat s>s)
&\le
\sum_{s<r\le m}
\bbP\left(|\hat L_s^*-\hat L_r^*|>\omega_r\right)\\
&\le c_{\al,\tau}\sum_{r=s+1}^{m}r^{-8}\\
&\le c_{\al,\tau}\int_s^\infty x^{-8}\,\mathrm dx
=\frac{c_{\al,\tau}}{7}s^{-7}.
\end{align*}
This proves \eqref{eq:probability bound} with
$c_{3,0}=c_{\al,\tau}/7$.

\end{proof}
\begin{proof}[Proof of Lemma~\ref{lem: fourth bound}]
  Let $\theta\in \Theta_s$ and ${\caP}_\xi\in {\cal G}_{\al,\tau}^\otimes$.

If $s_0=d+1$, then every admissible $s\in[d]$ satisfies
$s\le s_*=d$, which is the second case below.

If $s_0\le d$, we first consider $s=s_0$.
We have $\lambda_*(s)=0$ and $\hat{L}_s^*-L(\theta)=\sigma\sum_{j=1}^d\eta_j\xi_j.$
  Therefore
  \begin{align*}
    \bbE_{\theta,{\caP}_\xi}\left[\hat{L}_s^*-L(\theta)\right]^4\le &\sum_{j=1}^d \sigma^4\eta_j^4\bbE\xi_j^4+3\sum_{i\neq j}\sigma^4\eta_i^2\eta_j^2\bbE\xi_i^2\bbE\xi_j^2\\
    \le &\sigma^4(C_2^*\sum_{j=1}^d\eta_j^4+3(\sum_{j=1}^d\eta_j^2)^2).
  \end{align*}
  Since $\lambda_*(s_0)=0$, we have
\[
\Phi_*(s_0;\eta)
=
\log(es_0)\sum_{j=1}^d\eta_j^2
\le
\Phi_{\rm adp}(s_0;\eta).
\]
Therefore,
  $$\bbE_{\theta,{\caP}_\xi}\left[\hat{L}_s^*-L(\theta)\right]^4\lesssim \sigma^4\Phi_{\rm adp}^2(s;\eta).$$

  We next consider the case $s \le s_*$. Following a similar decomposition as in \Cref{sec:adaptive estimator}, we obtain
  \begin{align*}
    &\bbE_{\theta,{\caP}_\xi}\left(\hat{L}_s^*-L(\theta)\right)^4\\
    \le&27\,\biggl[\underbrace{\bbE_{\theta,{\caP}_\xi}\biggl(\sigma\sum_{j\le j_2(s)}\eta_j\xi_j\biggr)^4}_{\text{denoted as I}}+\underbrace{\bbE_{\theta,{\caP}_\xi}\biggl(\sum_{\substack{j>j_2(s)\\j\in S^c}}\sigma \eta_j\xi_j\bbo\left\{|\eta_j\xi_j|\ge  \factorone\tau\lambda_*(s)\right\}\biggr)^4}_{\text{denoted as II}}\\
    &\quad+\underbrace{\bbE_{\theta,{\caP}_\xi}\biggl(\sum_{\substack{j>j_2(s)\\j\in S}}\eta_jy_j\bbo\left\{|\eta_jy_j|\ge \factorone\sigma\tau\lambda_*(s)\right\} -\sum_{\substack{j>j_2(s)\\j\in S}}\eta_j\theta_j\biggr)^4}_{\text{denoted as III}}\biggr].
  \end{align*}

  \textit{Term I.} Similar to the previous analysis, we have
  $$\bbE_{\theta,{\caP}_\xi}\biggl(\sigma\sum_{j\le j_2(s)}\eta_j\xi_j\biggr)^4\lesssim\sigma^4\biggl(\sum_{j\le j_2(s)}\eta_j^2\biggr)^2\lesssim \sigma^4\Phi_{\rm adp}^2(s;\eta).$$

  \textit{Term II.} Let $U=\left\{j>j_2(s)\mid j\notin S\right\}$. We have
  \begin{align*}
    &\bbE_{\theta,{\caP}_\xi}\biggl(\sum_{j\in U}\sigma \eta_j\xi_j\bbo\left\{|\eta_j\xi_j|\ge \factorone\tau\lambda_*(s)\right\}\biggr)^4\\
    =&\sigma^4\left[\sum_{j\in U}\eta_{j}^4\bbE\biggl(\xi_{j}^4\bbo\left\{|\xi_{j}|\ge \frac{\factorone\tau\lambda_*(s)}{|\eta_{j}|}\right\}\biggr)\right.\\
    &\qquad+\left.3\sum_{\substack{j_1,j_2\in U\\j_1\neq j_2}}\eta_{j_1}^2\eta_{j_2}^2\bbE\left(\xi_{j_1}^2\xi_{j_2}^2\bbo\left\{|\xi_{j_1}|\ge  \frac{\factorone\tau\lambda_*(s)}{|\eta_{j_1}|}\right\}\bbo\left\{|\xi_{j_2}|\ge  \frac{\factorone\tau\lambda_*(s)}{|\eta_{j_2}|}\right\}\right)\right]\\
    \le &\sigma^4\left[\sqrt{C_4^*}\sum_{j\in U}|\eta_j|^4\exp\left(-\frac{(\beta_*)_+}{|\eta_j|^\al}\right)+3C_2^*\left(\sum_{j\in U}|\eta_j|^2\exp\left(-\frac{(\beta_*)_+}{|\eta_j|^\al}\right)\right)^2\right]\\
    \le&\sigma^4\left[\sqrt{C_4^*}\lambda_*^2(s)\sum_{j=1}^d|\eta_j|^2\exp\left(-\frac{(\beta_*)_+}{|\eta_j|^\al}\right)+3C_2^*\left(\sum_{j=1}^d|\eta_j|^2\exp\left(-\frac{(\beta_*)_+}{|\eta_j|^\al}\right)\right)^2\right]\\
    \lesssim&\sigma^4\Phi_{\rm adp}^2(s;\eta).
  \end{align*}

  \textit{Term III.}
  \begin{align*}
    &\bbE_{\theta,{\caP}_\xi}\biggl(\sum_{\substack{j>j_2(s)\\j\in S}}\eta_jy_j\bbo\left\{|\eta_jy_j|\ge \factorone\sigma\tau\lambda_*(s)\right\} -\sum_{\substack{j>j_2(s)\\j\in S}}\eta_j\theta_j\biggr)^4\\
    \le&8s^3\sum_{\substack{j>j_2(s)\\j\in S}}\left[\bbE_{\theta,{\caP}_\xi}\left(\eta_jy_j\bbo\left\{|\eta_jy_j|< \factorone\sigma\tau\lambda_*(s)\right\}\right)^4+\sigma^4\eta_j^4\cdot\bbE\xi_j^4\right]\\
    \le&8s^3\sum_{\substack{j>j_2(s)\\j\in S}}\left[2^{4/\al}\sigma^4\tau^4\lambda_*^4(s)+C_2^*\sigma^4\eta_j^4\right]\lesssim \sigma^4  \Phi_{\rm adp}^2(s;\eta).
  \end{align*}
  Therefore, we have obtained that
  $$\sup_{\theta\in \Theta_s}\bbE_{\theta,{\caP}_\xi}\left(\hat{L}_s^*-L(\theta)\right)^4\lesssim \sigma^4\Phi_{\rm adp}^2(s;\eta).$$
  Combined with \Cref{prop:linear decrease}, we have completed the proof.
\end{proof}

\section{Auxiliary Results}\label{sec:pf-aux}
\subsection{Existence and uniqueness of the solution to \Cref{eq:solution nonadaptive}}\label{sec:existence uniqueness linear}
In this subsection, we show that \Cref{eq:solution nonadaptive} has a unique solution.
As pointed out in \Cref{sec:nonadaptive minimax}, it is sufficient to show that the left-hand side of \Cref{eq:solution nonadaptive} is a continuous and strictly decreasing function in $\beta$ on $\bbR$.
\begin{lemma}\label{lem: decrease}
  The following function $\phi:\bbR\to\bbR$ is continuous and strictly decreasing:
  $$\phi(\beta)=\frac{\sum_{j=1}^d|\eta_j| \exp(-\beta/|\eta_j|^\al)}{\sqrt{\sum_{j=1}^d \eta_j^2\exp(-\beta/|\eta_j|^\al)}}.$$
\end{lemma}
\begin{proof}[Proof of Lemma~\ref{lem: decrease}]
  The function $\phi$ is continuous and differentiable. Its derivative can be written as
\begin{align*}
& \phi'(\beta)\\
=&\frac12 \Big(\sum_{j=1}^d \eta_j^2 e^{-\beta/|\eta_j|^\al}\Big)^{-3/2}\times  \\
& \qquad \left\{
\Big(\sum_{j=1}^d |\eta_j| e^{-\beta/|\eta_j|^\al}\Big)
\Big(\sum_{j=1}^d |\eta_j|^{2-\al} e^{-\beta/|\eta_j|^\al}\Big)
- 2 \Big(\sum_{j=1}^d |\eta_j|^{1-\al} e^{-\beta/|\eta_j|^\al}\Big)
\Big(\sum_{j=1}^d \eta_j^2 e^{-\beta/|\eta_j|^\al}\Big)\right\}\\
=& \frac{1}{2}
\Big(\sum_{j=1}^d \eta_j^2 e^{-\beta/|\eta_j|^\al}\Big)^{-3/2}
\Big(\sum_{j=1}^d {|\eta_j|^{1-\al}} e^{-\beta/|\eta_j|^\al}\Big)^2 \\[1ex]
&\quad \times \Bigg\{
\frac{\sum_{j=1}^d |\eta_j| e^{-\beta/|\eta_j|^\al}}
     {\sum_{j=1}^d {|\eta_j|}^{1-\al} e^{-\beta/|\eta_j|^\al}}
\cdot
\frac{\sum_{j=1}^d |\eta_j|^{2-\al}e^{-\beta/|\eta_j|^\al}}
     {\sum_{j=1}^d {|\eta_j|}^{1-\al} e^{-\beta/|\eta_j|^\al}}
- 2 \cdot
\frac{\sum_{j=1}^d \eta_j^2 e^{-\beta/|\eta_j|^\al}}
     {\sum_{j=1}^d {|\eta_j|}^{1-\al} e^{-\beta/|\eta_j|^\al}}
\Bigg\}.
\end{align*}
  Define a random variable $Y$ following a probability measure $\mu=\sum_{j=1}^d w_j\delta_{|\eta_j|}$, where
  $$w_j=\frac{{|\eta_j|}^{1-\al}\exp(-\beta/|\eta_j|^\al)}{\sum_{j=1}^d {|\eta_j|}^{1-\al}\exp(-\beta/|\eta_j|^\al)}.$$
  Since $|\eta_j|>0,\forall j\in [d]$, we have $\bbE Y^{\al+1}>0$  and $\bbE Y \cdot \bbE Y^{\al}\le \bbE Y^{1+\al}$. Therefore,
  $$
\begin{aligned}
& \phi'(\beta)\\
= & \frac{1}{2}\left(\sum_{j=1}^d \eta_j^2\exp(-\beta/|\eta_j|^\al)\right)^{-3/2}\left(\sum_{j=1}^d |\eta_j|^{1-\al}\exp(-\beta/|\eta_j|^\al)\right)^2\left[\bbE Y\cdot \bbE Y^\al-2\bbE Y^{1+\al}\right]\\
< & 0.
\end{aligned}
$$

\end{proof}

\subsection{Proof of \Cref{prop:nonadaptive linear}}
\begin{proof}
  By the definition of $j_1$, the following holds up to absolute constants:
  \begin{align*}
    [\lambda_{o} s+\nu]^2\asymp&\lambda_{o}^2 s^2+\nu^2\\
    \asymp &\lambda_{o}^2 s^2+\sum_{j=1}^d \eta_j^2\exp(-\beta_+/|\eta_j|^\al)\\
    \asymp &\lambda_{o}^2 s^2+\sum_{j\le j_1} \eta_j^2\exp(-\beta_+/|\eta_j|^\al)+\sum_{j>j_1} \eta_j^2\exp(-\beta_+/|\eta_j|^\al)\\
    \asymp &\lambda_{o}^2 s^2+\sum_{j\le j_1} \eta_j^2+\sum_{j>j_1} \eta_j^2\exp(-\beta_+/|\eta_j|^\al).
  \end{align*}

  If $\sum_{j>j_1} \eta_j^2\exp(-\beta_+/|\eta_j|^\al)\le \sum_{j\le j_1} \eta_j^2\exp(-\beta_+/|\eta_j|^\al)$, then we have $[\lambda_{o} s+\nu]^2\asymp \lambda_{o}^2 s^2+\sum_{j\le j_1} \eta_j^2$. Otherwise, we have
  \begin{align*}
    \sum_{j>j_1}\eta_j^2\exp(-\beta_+/|\eta_j|^\al)&\le \lambda_{o} \sum_{j>j_1}|\eta_j|\exp(-\beta_+/|\eta_j|^\al)\\
    &\le \frac{1}{2}\lambda_{o} s\cdot\sqrt{\sum_{j=1}^d\eta_j^2\exp(-\beta_+/|\eta_j|^\al)}\\
    &\le \frac{1}{2}\lambda_{o} s\cdot\sqrt{2\sum_{j>j_1} \eta_j^2\exp(-\beta_+/|\eta_j|^\al)},
  \end{align*}
  where the first inequality holds since $|\eta_j|\le \lambda_{o}$ for $j> j_1$, the second inequality holds since the definition of $\beta$ implies that $\sum_{j=1}^d |\eta_j|\exp(-\beta_+/|\eta_j|^\al)\le s/2\cdot \sqrt{\sum_{j=1}^d \eta_j^2\exp(-\beta_+/|\eta_j|^\al)}$. Therefore, we have $\sum_{j>j_1} \eta_j^2\exp(-\beta_+/|\eta_j|^\al)\lesssim \lambda_{o}^2 s^2$, which completes the proof.
\end{proof}
\subsection{Proof of \Cref{prop:linear decrease}}
\begin{proof}

Note that
 \begin{align*}
  \Phi_*(s;\eta)=&\log (es)\left[\sum_{j=1}^d\eta_j^2\exp\left(-\frac{(\beta_*)_+}{|\eta_j|^\al}\right)+\lambda^2_*(s)\left(\frac{s}{\sqrt{\log (es)}}\right)^2\right]\\
  \asymp &\log (es) \left[\sqrt{\sum_{j=1}^d\eta_j^2\exp\left(-\frac{(\beta_*)_+}{|\eta_j|^\al}\right)}+\lambda_*(s)\cdot \frac{\sum_{j=1}^d|\eta_j| \exp\left(-{\beta_*}/{|\eta_j|^\al}\right)}{\sqrt{\sum_{j=1}^d \eta_j^2\exp\left(-{\beta_*}/{|\eta_j|^\al}\right)}}\right]^2
 \end{align*}
 where the second equation uses the definition of $\lambda_*(s)$ in \eqref{eq:solution linear adaptive}.
 Note that $\beta_*$ is a non-increasing function of \(s\).
  Let $$\phi_1(\beta)=\sqrt{\sum_{j=1}^d \eta_j^2\exp(-\beta/|\eta_j|^\al)}+\beta^{1/\al}\cdot \frac{\sum_{j=1}^d|\eta_j| \exp(-\beta/|\eta_j|^\al)}{\sqrt{\sum_{j=1}^d \eta_j^2\exp(-\beta/|\eta_j|^\al)}}$$
  and
  $$\phi_2(\beta)=\frac{\sum_{j=1}^d \eta_j^2\exp(-\beta/|\eta_j|^\al)}{\sum_{j=1}^d |\eta_j|\exp(-\beta/|\eta_j|^\al)}.$$

We first consider the case
$\beta_*(s)\ge\beta_*(s^\prime)>0$.

Since $\log(e(s \wedge s_0))$ is non-decreasing, whereas $\log(es)/{s^2}$ is nonincreasing,
  we only need to show that there exists a constant \(c>0\) such that for all \(\beta_1\ge \beta_2>0\), $$\phi_1(\beta_1)\le c\cdot \phi_1(\beta_2)\quad\text{and}\quad\phi_2(\beta_1)\ge c\cdot \phi_2(\beta_2).$$
  The proof for $\phi_2(\beta)$ is simple, as we have

  \begin{align*}
\phi_2'(\beta)
= \Biggl(\sum_{j=1}^d |\eta_j| \exp\!\bigl(-\beta/|\eta_j|^\al\bigr)\Biggr)^{-2}
\Biggl[&
\sum_{j=1}^d \eta_j^2 \exp\!\bigl(-\beta/|\eta_j|^\al\bigr)\!\sum_{j=1}^d |\eta_j|^{1-\al}\exp\!\bigl(-\beta/|\eta_j|^\al\bigr) \Biggr.
\\
& \Biggl. - \sum_{j=1}^d |\eta_j|^{2-\al}\exp\!\bigl(-\beta/|\eta_j|^\al\bigr)\!\sum_{j=1}^d |\eta_j|\exp\!\bigl(-\beta/|\eta_j|^\al\bigr)
\Biggr].
\end{align*}

  The above term is non-negative following the proof of Lemma~\ref{lem: decrease}. Therefore, $\phi_2(\beta)$ is non-decreasing in $\beta$.
  The remaining analysis concerns $\phi_1(\beta)$.
  Following the notation in \Cref{sec:existence uniqueness linear}, we have
  \begin{align*}
    \phi_1(\beta)=\sqrt{\sum_{j=1}^d \eta_j^2\exp(-\beta/|\eta_j|^\al)}+\beta^{1/\al}\phi(\beta).
  \end{align*}
  There are two components in $\phi_1(\beta)$, and the first term $$\sqrt{\sum_{j=1}^d \eta_j^2\exp(-\beta/|\eta_j|^\al)}$$
  is clearly strictly decreasing in $\beta$. Therefore, for $\beta=\beta_1$, we can assume
  \begin{equation*}
    \sqrt{\sum_{j=1}^d\eta_j^2\exp(-\beta/|\eta_j|^\al)}< \left(\frac{\al\beta}{2}\right)^{1/\al} \phi(\beta).
  \end{equation*}
  Otherwise, we have
  \begin{align*}
    \phi_1(\beta_1)=&\sqrt{\sum_{j=1}^d \eta_j^2\exp(-\beta_1/|\eta_j|^\al)}+\beta_1^{1/\al}\phi(\beta_1)\\
    \le &\left[1+\left(\frac{2}{\al}\right)^{1/\al}\right]\sqrt{\sum_{j=1}^d \eta_j^2\exp(-\beta_1/|\eta_j|^\al)}\\
    \le &\left[1+\left(\frac{2}{\al}\right)^{1/\al}\right]\sqrt{\sum_{j=1}^d \eta_j^2\exp(-\beta_2/|\eta_j|^\al)}\\
    \le &\left[1+\left(\frac{2}{\al}\right)^{1/\al}\right]\phi_1(\beta_2).
  \end{align*}
  Define
  $$\beta_{\min}=\min\left\{\beta^\prime\ge \beta_2: \sqrt{\sum_{j=1}^d\eta_j^2\exp(-\beta/|\eta_j|^\al)}\le \left(\frac{\al\beta}{2}\right)^{1/\al}\phi(\beta)\,\,\text{holds}\,\, \forall\beta\in [\beta^\prime,\beta_1]\right\}.$$
  The existence of such $\beta_{\min}$ is guaranteed by the continuity of functions involved in the above definition.

  Below we are going to show that $\phi_1(\beta)$ is strictly decreasing in $\beta\in [\beta_{\min},\beta_1]$.

  Following the calculations in Section \ref{sec:existence uniqueness linear}, we have

  \scalebox{0.75}{
  \parbox{\linewidth}{
  \begin{align*}
    &\phi'_1(\beta)=-\frac{\sum_{j=1}^d |\eta_j|^{2-\al}\exp(-\beta/|\eta_j|^\al)}{2\sqrt{\sum_{j=1}^d \eta_j^2\exp(-\beta/|\eta_j|^\al)}}+\frac{\beta^{1/\al-1}\sum_{j=1}^d|\eta_j|\exp(-\beta/|\eta_j|^\al)}{\al\sqrt{\sum_{j=1}^d \eta_j^2\exp(-\beta/|\eta_j|^\al)}}+\beta^{1/\al}\phi'(\beta)\\
    =&{\tiny -\frac{\sum_{j=1}^d |\eta_j|^{2-\al}\exp(-\beta/|\eta_j|^\al)}{2\sqrt{\sum_{j=1}^d \eta_j^2\exp(-\beta/|\eta_j|^\al)}}+\frac{\beta^{1/\al-1}\sum_{j=1}^d|\eta_j|\exp(-\beta/|\eta_j|^\al)}{\al\sqrt{\sum_{j=1}^d \eta_j^2\exp(-\beta/|\eta_j|^\al)}}-\beta^{1/\al}\frac{\sum_{j=1}^d|\eta_j|^{1-\al}\exp(-\beta/|\eta_j|^\al)}{2\sqrt{\sum_{j=1}^d \eta_j^2\exp(-\beta/|\eta_j|^\al)}}}\\
+&\beta^{1/\al}\left[\frac{{\tiny\left(\sum_{j=1}^d |\eta_j| \exp(-\beta/|\eta_j|^\al)\right)\left(\sum_{j=1}^d |\eta_j|^{2-\al}\exp(-\beta/|\eta_j|^\al)\right)-\left(\sum_{j=1}^d |\eta_j|^{1-\al}\exp(-\beta/|\eta_j|^\al)\right)\left(\sum_{j=1}^d \eta_j^2\exp(-\beta/|\eta_j|^\al)\right)}}{2\cdot\left(\sum_{j=1}^d \eta_j^2\exp(-\beta/|\eta_j|^\al)\right)^{3/2}}\right].
  \end{align*}
  }
  }

  The last term is non-positive following the proof of Lemma~\ref{lem: decrease}. And the first three terms can be rearranged as follows:

    \scalebox{0.85}{
  \parbox{\linewidth}{
  \begin{align*}
    &-\frac{\sum_{j=1}^d |\eta_j|^{2-\al}\exp(-\beta/|\eta_j|^\al)}{2\sqrt{\sum_{j=1}^d \eta_j^2\exp(-\beta/|\eta_j|^\al)}}+\frac{\beta^{1/\al-1}\sum_{j=1}^d|\eta_j|\exp(-\beta/|\eta_j|^\al)}{\al\sqrt{\sum_{j=1}^d \eta_j^2\exp(-\beta/|\eta_j|^\al)}}-\frac{\beta^{1/\al}\sum_{j=1}^d|\eta_j|^{1-\al}\exp(-\beta/|\eta_j|^\al)}{2\sqrt{\sum_{j=1}^d \eta_j^2\exp(-\beta/|\eta_j|^\al)}}\\
    =&-\frac{\beta^{1/\al}}{2\sqrt{\sum_{j=1}^d \eta_j^2\exp(-\beta/|\eta_j|^\al)}}\sum_{j=1}^d \left[\left(\beta^{-1/\al}|\eta_j|-\frac{2}{\al\beta}|\eta_j|^\al+1\right)|\eta_j|^{1-\al}\exp(-\beta/|\eta_j|^\al)\right].
  \end{align*}
  }
  }

  \begin{proposition}\label{prop:inequality}
For any $\al > 0$, we have
\[
f(t) = \al t^{2} - (\al+1)t + t^{1-\al}\ge 0,
\quad \forall t > 0.\]
\end{proposition}

Applying Proposition~\ref{prop:inequality} with
\(
t = {|\eta_j|}(\frac{2}{\al\beta})^{1/\al}
\),
we obtain
\[
\al |\eta_j|^2 \left(\frac{2}{\al\beta}\right)^{2/\al}
-(\al+1)|\eta_j|\left(\frac{2}{\al\beta}\right)^{1/\al}
+|\eta_j|^{1-\al}\left(\frac{2}{\al\beta}\right)^{1/\al-1} \;\ge\; 0.
\]
Rearranging yields
\[
|\eta_j|^{1-\al} - \frac{2}{\al\beta}\,|\eta_j|
\;\ge\; \frac{2}{\beta}\left[\,|\eta_j| - \left(\frac{2}{\al\beta}\right)^{1/\al} |\eta_j|^2 \right].
\]

From the definition of $\beta_{\min}$, we have
  $$\sum_{j=1}^d\left[\left(\frac{\al\beta}{2}\right)^{1/\al}|\eta_j|-\eta_j^2\right]\exp(-\beta/|\eta_j|^\al)\ge 0,\forall \beta\in [\beta_{\min},\beta_1].$$
  Therefore, we have
  \begin{align*}
    &\sum_{j=1}^d\left(|\eta_j|^{1-\al}-\frac{2}{\al\beta}|\eta_j|\right)\exp(-\beta/|\eta_j|^\al)\\
    \ge &\frac{2}{\beta}\left(\frac{2}{\al\beta}\right)^{1/\al}\sum_{j=1}^d \left[\left(\frac{\al\beta}{2}\right)^{1/\al}|\eta_j|-\eta_j^2\right]\exp(-\beta/|\eta_j|^\al)\ge 0.
  \end{align*}
  which implies that $\phi_1'(\beta) \le 0$ for all $\beta \in [\beta_{\min}, \beta_1]$, and hence $\phi_1(\beta)$ is non-increasing on the interval $[\beta_{\min}, \beta_1]$.

  If $\beta_{\min}=\beta_2$, then the proof is completed. Otherwise, from the definition of $\beta_{\min}$, we have
  $$
\sqrt{
\sum_{j=1}^d\eta_j^2
\exp\left(-\frac{\beta_{\min}}{|\eta_j|^\al}\right)
}
=
\left(\frac{\al\beta_{\min}}{2}\right)^{1/\al}
\phi(\beta_{\min}).
$$
Therefore,
  \begin{align*}
    \phi_1(\beta_1)\le \phi_1(\beta_{\min})=&\left[1+\left(\frac{2}{\al}\right)^{1/\al}\right]\sqrt{\sum_{j=1}^d\eta_j^2\exp(-\beta_{\min}/|\eta_j|^\al)}\\\le&\left[1+\left(\frac{2}{\al}\right)^{1/\al}\right] \sqrt{\sum_{j=1}^d\eta_j^2\exp(-\beta_{2}/|\eta_j|^\al)}\\
    \le &\left[1+\left(\frac{2}{\al}\right)^{1/\al}\right] \phi_1(\beta_2).
  \end{align*}

\bigskip

We next extend the two comparisons to the case
$\beta_*(s^\prime)\le 0$. For notational ease, define
$$
V=\sum_{j=1}^d\eta_j^2.
$$

First, if $\beta_*(s)\le0$, then the monotonicity of $\beta_*$ implies that
$\beta_*(r)\le0$ for every $r\in[s,s^\prime]$. Hence
$\Phi_*(r;\eta)=V\log(er)$ on this interval, and therefore
\begin{align*}
\frac{\Phi_*(s;\eta)}{\log(es)}
&=
\frac{\Phi_*(s^\prime;\eta)}{\log(es^\prime)}
=V,\\
\frac{\Phi_*(s^\prime;\eta)}
{(s^\prime)^2\log(es^\prime)}
&=
\frac{V}{(s^\prime)^2}
\le
\frac{V}{s^2}
=
\frac{\Phi_*(s;\eta)}{s^2\log(es)}.
\end{align*}

Second, suppose
$\beta_*(s)>0\ge\beta_*(s^\prime)$.
The defining equation for $\beta_*$ and the strict monotonicity of
$\phi$ imply that $\beta_*$ is continuous and non-increasing.
Thus, there exists $s_\circ\in(s,s^\prime]$ such that
$\beta_*(s_\circ)=0$.
The comparisons established for positive values of $\beta_*$
extend by continuity to the pair $(s,s_\circ)$.
Since $\Phi_*(r;\eta)=V\log(er)$ for
$r\in[s_\circ,s^\prime]$, it follows that
\begin{align*}
\frac{\Phi_*(s;\eta)}{\log(es)}
&\lesssim
\frac{\Phi_*(s_\circ;\eta)}{\log(es_\circ)}
=
\frac{\Phi_*(s^\prime;\eta)}{\log(es^\prime)}=V,\\
\frac{\Phi_*(s^\prime;\eta)}
{(s^\prime)^2\log(es^\prime)}
&=
\frac{V}{(s^\prime)^2}
\le
\frac{V}{s_\circ^2}
=
\frac{\Phi_*(s_\circ;\eta)}
{s_\circ^2\log(es_\circ)}
\lesssim
\frac{\Phi_*(s;\eta)}{s^2\log(es)}.
\end{align*}
Therefore, the two comparisons are proved.

We finally derive the remaining two consequences. By the first
comparison and the monotonicity of $\log(er)$, for
$1\le s\le s^\prime\le s_0\wedge d$,
$$
\Phi_*(s;\eta)
\lesssim
\frac{\log(es)}{\log(es^\prime)}
\Phi_*(s^\prime;\eta)
\le
\Phi_*(s^\prime;\eta).
$$
Thus, $\Phi_*(r;\eta)$ is almost non-decreasing on
$[1,s_0\wedge d]$. If $s_0=d+1$, this interval is $[1,d]$.
If $s_0\le d$, then $\Phi_*(r;\eta)$ is constant for
$r\ge s_0$, so the same conclusion again holds on $[1,d]$.

For $0<\al<2$, the function
$\Phi_*(1;\eta)\log^2(e(r\wedge s_0))$ is non-decreasing in $r$.
Together with the definition of $\Phi_{\rm adp}$, this shows that
$$
\Phi_{\rm adp}(s;\eta)
\lesssim
\Phi_{\rm adp}(s^\prime;\eta),
\qquad
1\le s\le s^\prime\le d.
$$
For $\al\ge2$, the same conclusion follows directly from
$\Phi_{\rm adp}(r;\eta)=\Phi_*(r;\eta)$. This proves the first
consequence.

Next, applying the second comparison with $s=1$ and
$s^\prime=r$ gives
$$
\Phi_*(r;\eta)
\lesssim
r^2\log(er)\Phi_*(1;\eta),
\qquad
1\le r\le s_0\wedge d.
$$
When $0<\al<2$, for $1\le r\le s_0\wedge d$,
$$
\Phi_*(1;\eta)\log^2(e(r\wedge s_0))
=
\Phi_*(1;\eta)\log^2(er)
\le
r^2\log(er)\Phi_*(1;\eta).
$$
For $\al\ge2$, there is no additional term in
$\Phi_{\rm adp}$. Since
$\Phi_*(1;\eta)=\Phi_{\rm o}(1;\eta)$, we conclude that
$$
\Phi_{\rm adp}(r;\eta)
\lesssim
r^2\log(er)\Phi_{\rm o}(1;\eta),
\qquad
1\le r\le s_0\wedge d.
$$
This proves the second consequence.

\end{proof}

\subsection{Proof of \Cref{prop:adaptive rate}}

\begin{proof}
Choose $\tau_{\al,\gamma_0}$ large enough that the lower bound in \Cref{thm: nonadaptive linear} and \Cref{thm: adaptive lower} with $\gamma=\gamma_0/2$ hold.

By \Cref{thm: adaptive upper},
$\psi_*(s;\eta)=\sigma^2\Phi_{\rm adp}(s;\eta)$ satisfies
\eqref{eq:adaptive_upper_bound}. It remains to verify the second
requirement in \Cref{def: adaptive rate}.

Consider any $\psi_*^\prime$ satisfying \eqref{eq:alternative_rate} and
\eqref{eq:improvement_condition}. For $r\in[d]$, define
$$
q_d(r)=\frac{\psi_*^\prime(r;\eta)}
{\sigma^2\Phi_{\rm adp}(r;\eta)},
$$
and denote the minimum as
$$
m_d=\min_{r\in[d]}q_d(r)=q_d(\widetilde s),
$$
for some $\widetilde s\in[d]$. By \eqref{eq:improvement_condition},
$m_d\to0$.

Since the lower bound in \Cref{thm: nonadaptive linear} holds for every estimator, \eqref{eq:alternative_rate} implies
\begin{equation}\label{eq:lower-q_d_r}
    q_d(r)\gtrsim\Phi_{\rm o}(r;\eta)/\Phi_{\rm adp}(r;\eta).
\end{equation}
Hence, by \eqref{eq:adaptive lower condition 1}, $q_d(r)$ is bounded away from zero uniformly over $r\le s_{\rm cut}$. Since $m_d\to0$, we have $\widetilde s>s_{\rm cut}$ for all sufficiently large $d$.

We next show that $\widetilde s\wedge s_0\to\infty$.
Since $\sigma^2\Phi_{\rm o}(s;\eta)$ is the minimax rate by \Cref{thm: nonadaptive linear} and
$\Theta_1\subseteq\Theta_r$, we have
$\Phi_*(1;\eta)\asymp\Phi_{\rm o}(1;\eta)
\lesssim\Phi_{\rm o}(r;\eta)$ for every $r\in[d]$.
Together with \eqref{eq:adaptive-oracle-comparison} and
\eqref{eq:threshold adaptive}, this yields
$\Phi_{\rm adp}(r;\eta)
\lesssim
\log^2\!\bigl(e(r\wedge s_0)\bigr)\Phi_{\rm o}(r;\eta)$ for every $r\in[d]$.
This, together with \eqref{eq:lower-q_d_r} at $r=\widetilde s$, gives
$$
m_d\gtrsim
\log^{-2}\!\bigl(e(\widetilde s\wedge s_0)\bigr).
$$
Since $m_d\to0$, we conclude that
$\widetilde s\wedge s_0\to\infty$.

For each $d\ge d_0$, choose $\widehat T_d$ so that the maximum in
\eqref{eq:alternative_rate} is at most $2C^\prime$.
Its maximal risk at $\widetilde s$ is at most
$2C^\prime\psi_*^\prime(\widetilde s;\eta)
=2C^\prime m_d\sigma^2
\Phi_{\rm adp}(\widetilde s;\eta)$.
Since $m_d\to0$, this bound is at most
$\sigma^2\Phi_{\rm adp}(\widetilde s;\eta)/C_0$
for all sufficiently large $d$.
Hence \eqref{eq:adaptive lower small} holds with
$s=\widetilde s$, and \Cref{thm: adaptive lower} applies.

At sparsity one, $\widehat T_d$ has maximal risk at most
$2C^\prime\psi_*^\prime(1;\eta)
=2C^\prime q_d(1)\sigma^2\Phi_{\rm adp}(1;\eta)$.
On the other hand, \eqref{eq:adaptive-lower} in \Cref{thm: adaptive lower} with
$s=\widetilde s$ gives a lower bound for the same risk.
Comparing the two bounds yields

\begin{equation}\label{eq:adaptive-rate-qd-one}
\begin{aligned}
q_d(1)
&\gtrsim
\frac{\Phi_{\rm adp}(\widetilde s;\eta)}
{(\widetilde s\wedge s_0)^{\gamma_0/2}
 \Phi_{\rm adp}(1;\eta)}\\
&\asymp
\frac{\Phi_{\rm adp}(\widetilde s;\eta)}
{(\widetilde s\wedge s_0)^{\gamma_0/2}
 \Phi_{\rm o}(1;\eta)} \\
&\gtrsim
(\widetilde s\wedge s_0)^{\gamma_0/2},
\end{aligned}
\end{equation}
where the second line follows from
$\Phi_{\rm adp}(1;\eta)\asymp\Phi_{\rm o}(1;\eta)$,
and the last line follows from
\eqref{eq:adaptive lower condition 2}, since
$\widetilde s>s_{\rm cut}$ for all sufficiently large $d$.

Combining \eqref{eq:adaptive-rate-qd-one} with
$m_d\gtrsim
\log^{-2}\!\bigl(e(\widetilde s\wedge s_0)\bigr)$ gives
$$
q_d(1)m_d
\gtrsim
\frac{(\widetilde s\wedge s_0)^{\gamma_0/2}}
{\log^2\!\bigl(e(\widetilde s\wedge s_0)\bigr)}
\longrightarrow\infty.
$$
Therefore, \eqref{eq:loss_condition} holds with $\bar{s}=1$,
which completes the proof.

\end{proof}

\subsection{Proof of \Cref{prop:exponential match}}
\begin{proof}

To ease notation, define
$$
t=\log^{1/\al}(ed/s),
\qquad
m=s^2\wedge d.
$$

By the definition of $j_1(s)$, we have
$|\eta_j|<\lambda_{o}$ for every $j>j_1(s)$. It follows that
\begin{align*}
\sum_{j\le s^2\wedge d}\eta_j^2 =\sum_{j\le m}\eta_j^2
&\le
\sum_{j\le j_1(s)}\eta_j^2
+m\lambda_{o}^2\\
&\le
\sum_{j\le j_1(s)}\eta_j^2
+s^2\lambda_{o}^2\\
&\lesssim
\Phi_{\rm o}(s;\eta),
\end{align*}
where the last inequality follows from
\Cref{prop:nonadaptive linear}.

It remains to compare the two partial sums.
Define $q=j_3(s)=\lceil s^2t^2\rceil\wedge d$.
Since $t\ge 1$, we have $m\le q$.

If $m=d$, then $q=d$, and the
desired comparison follows from $\log d\ge\log 2$.

Now suppose that
$m=s^2<d$.
We have $q\le\lceil mt^2\rceil
\le mt^2+1
\le2mt^2$.
Because $|\eta_j|$ is nonincreasing in $j$, we further have
\begin{align*}
\sum_{j\le q}\eta_j^2
&=
\sum_{j\le m}\eta_j^2
+\sum_{m<j\le q}\eta_j^2\\
&\le
\sum_{j\le m}\eta_j^2
+(q-m)\eta_m^2\\
&\le
\sum_{j\le m}\eta_j^2
+2mt^2\eta_m^2\\
&\le
(1+2t^2)\sum_{j\le m}\eta_j^2\\
&\le
3t^2\sum_{j\le m}\eta_j^2.
\end{align*}
Since
$t^2
=
\log^{2/\al}(ed/s)
\le
\log^{2/\al}(ed)
\lesssim
\log^{2/\al}d$, we have
$$
\sum_{j\le m}\eta_j^2
\gtrsim
\frac{1}{\log^{2/\al}d}
\sum_{j\le q}\eta_j^2.
$$
  This completes the proof.

\end{proof}

\section{Additional proof of the lower bounds}\label{sec:additional lower bound}
\subsection{Proof of Theorem~\ref{thm: exponential lower}}\label{subsec:proof exponential lower}
{This subsection contains the genuinely new ingredient for the non-symmetric model. The construction couples a sparse prior on the signal with a least-favorable asymmetric perturbation of the noise, and the proof separates into two parts: first we show that the alternative puts enough mass on large values of $L(\theta)$, and then we control the total variation distance through a Hellinger bound.}

  Let $f_0:\bbR\to[0,\infty)$ be a probability density with the following properties: $f_0$ is continuously differentiable, symmetric about $0$, supported on $[-3/2,3/2]$, with variance $1$ and finite Fisher information $I_{f_0}=\int (f_0^\prime(x))^2 [f_0(x)]^{-1}{\diff} x<\infty$. The existence of such a function is guaranteed by \cite[Lemma 7]{comminges2021adaptive}. Without loss of generality, we assume that $\sigma = 1$.

  Note that based on \eqref{eq:exponential lower symmetric} and \Cref{prop:exponential match}, we only consider the case where
  $$s^2\lambda^2_{\cal H}\ge C\sum_{j\le s^2}\eta_j^2$$
  for some sufficiently large constant $C>0$.
  The above inequality implies that
  $$\lambda_{\cal H}\ge \sqrt{C}\cdot \sqrt{\frac{1}{s^2}\sum_{j\le s^2}\eta_j^2}\ge \sqrt{C}\cdot |\eta_{s^2}|, $$
  which implies $\lambda_{\cal H}/ |\eta_j| \ge \sqrt{C}$, for any $j \ge {s^2}$.

  Let $b_1,\ldots,b_d$ be independently Bernoulli random variables with the probability of success $\bbP(b_j=1)=\pi_j,\forall j=1,\ldots,d$, where $\pi_j$ is defined as follows:
  $$\pi_j=\left\{\begin{array}{cc}
    0,& j< s^2\\
    \ds c\cdot \exp\left(-\left|\frac{\lambda_{\cal H}}{\eta_j}\right|^\al\right),& j\ge s^2
  \end{array}\right.,$$
  where $c>0$ is a constant specified later. Then by \eqref{eq:solution exponential 0}, we have
  $\sum_{j=1}^d\pi_j=c\cdot s.$
  Let $\mu$ denote the distribution of $\theta = (\gamma_1b_1,\ldots,\gamma_db_d)$ with $\gamma_j=c_2\cdot \lambda_{\cal H}/\eta_j,\forall j\in [d]$, where $c_2$ is some sufficiently small constant.

  Note that for $j\ge s^2$, we have
  $$\gamma_j\cdot \pi_j=
      \ds  c \cdot c_2\cdot \frac{\lambda_{\cal H}}{\eta_j}\cdot \exp\left(-\left|\frac{\lambda_{\cal H}}{\eta_j}\right|^\al\right),
    $$
    and
    $$\gamma_j^2\cdot \pi_j(1-\pi_j) \le
      \ds  c \cdot c_2^2 \cdot (\frac{\lambda_{\cal H}}{\eta_j})^2\cdot \exp\left(-\left|\frac{\lambda_{\cal H}}{\eta_j}\right|^\al\right).
    $$

Note that $\sup_{x}|x|^k \exp(-|x|^\al)=(e\al/k)^{-k/\al}$ only depends on $\al$.
For $k=1,2$, one can find $c_{2,1}$ such that for any $c_2\le c_{2,1}$, the following inequalities hold for all $j$:
\begin{equation}\label{eq:lower-exp-gamma-bound}
    |\gamma_j\pi_j| < 1/2, \qquad \gamma_j^2\cdot \pi_j(1-\pi_j)<1/4.
\end{equation}

  Then, we have
  \begin{align*}
  \|\theta\|_0 &=\sum_{j=1}^d b_j, \\
    L(\theta)&=\sum_{j=1}^d \eta_j\cdot \theta_j=c_2 \cdot \lambda_{\cal H}\cdot \|\theta\|_0, \\
    \bbE_\mu\left[\exp( t \|\theta\|_0)\right]&=\prod_{j=1}^d \left[1+\pi_j\cdot (\exp(t)-1)\right], \quad t\in \mathbb{R}.
  \end{align*}

  Using the basic inequality that  $1+x\le \exp(x)$ for $\forall x\in \mathbb{R}$, we have
  \begin{equation}\label{eq: pro 1}
  \begin{aligned}
    \bbP_\mu(\left\|\theta\right\|_0> s)
    \le& \exp(-s)\cdot\bbE_\mu\left[\exp(\|\theta\|_0)\right]\\
    = & \exp(-s)\cdot \prod_{j=1}^d\bbE\left[\exp(b_j)\right]\\
    =&\exp(-s)\cdot\prod_{j=1}^d \left[1+\pi_j\cdot (e-1)\right]\\
    \le& \exp(-s)\exp\left[(e-1)\cdot\sum_{j=1}^d \pi_j\right]=\exp(s\cdot ((e-1)c-1)),
  \end{aligned}
  \end{equation}
  and
  \begin{equation}\label{eq: pro 2}
  \begin{aligned}
    \bbP_\mu(\left\|\theta\right\|_0<c_1\cdot s)\le &\exp(c_1\cdot s)\bbE_\mu\left[\exp(-\|\theta\|_0)\right]\\
    =  &\exp(c_1\cdot s)\prod_{j=1}^d\bbE\left[\exp(-b_j)\right]\\
    = &\exp(c_1\cdot s)\cdot \prod_{j=1}^d \left[1+\pi_j\cdot (e^{-1}-1)\right]\\
    \le &\exp(c_1\cdot s)\cdot \exp\left[(e^{-1}-1)\sum_{j=1}^d \pi_j\right]=\exp(s\cdot (c_1-(1-e^{-1})\cdot c)).
  \end{aligned}
  \end{equation}

  Define
\[
h_0(x)
=
\exp\left\{3\bigl((e-1)x-1\bigr)\right\}
+\exp\left\{-3(1-e^{-1})x\right\},
\qquad x\ge 0.
\]
A direct calculation shows that \(h_0\) attains its minimum at
\[
x_0=\frac{2}{3(e-e^{-1})}\in(0,1),
\]
with
\[
h_0(x_0)
=
(1+e^{-1})\exp\left(-\frac{2}{e+1}\right)
<1.
\]
Set
\[
c=x_0=\frac{2}{3(e-e^{-1})}.
\]
Since \(h_0(c)<1\), we may choose a sufficiently small constant
\(c_1>0\) such that
\[
c_1<(1-e^{-1})c
\quad\text{and}\quad
h_0(c)+e^{3c_1}-1<1.
\]
For example, \(c_1=1/100\) satisfies these conditions.

Since both \((e-1)c-1\) and \(c_1-(1-e^{-1})c\) are negative, for
every \(s\ge 3\), the sum of the bounds in
\eqref{eq: pro 1} and \eqref{eq: pro 2} satisfies
\begin{align*}
  &\exp\left\{s\bigl((e-1)c-1\bigr)\right\}
  +\exp\left\{s\bigl(c_1-(1-e^{-1})c\bigr)\right\} \\
  &\quad\le
  \exp\left\{3\bigl((e-1)c-1\bigr)\right\}
  +\exp\left\{3\bigl(c_1-(1-e^{-1})c\bigr)\right\} \\
  &\quad\le h_0(c)+e^{3c_1}-1
  =:\delta_0<1.
\end{align*}
Consequently,
\begin{align*}
  \bbP_\mu\left(
    \theta\in\Theta_s,\,
    L(\theta)\ge c_1c_2s\lambda_{\mathcal H}
  \right)
  &=\bbP_\mu\left(c_1s\le\|\theta\|_0\le s\right) \ge 1-\delta_0>0.
\end{align*}

With a little abuse of notation, we let $f_j$ be the density of the random variable $\tilde{\delta}_j=\sigma_j \delta_j+\gamma_j\cdot(b_j - \pi_j)$ with $\delta_j$ being i.i.d. random variables with density $f_0$ and $$\sigma_j=\sqrt{1-\gamma_j^2\cdot \pi_j(1-\pi_j)}.$$
  $\sigma_j$ is well-defined in view of \Cref{eq:lower-exp-gamma-bound}.
  By construction, we have $$\bbE\, \tilde{\delta}_j=0\quad \text{and}\quad \bbE\, \tilde{\delta}_j^2=1,\quad  \forall  j=1,\ldots,d.$$

We establish a tail bound for \(\tilde{\delta}_j\) by considering two cases.

\textbf{(i)}.
  For any $t\in (0,2]$,
  if $\tau\ge 2/\log^{1/\al} 2$, then $(t/\tau)^\al \le \log(2)$. Thus,
  $$
\bbP(|\tilde{\delta}_j|\ge t)\le 1 \le  2\exp\left\{-\left(\frac{t}{\tau}\right)^\al\right\}.
  $$
Set $\tau_\al=2/\log^{1/\al} 2$.

\textbf{(ii)}.
  For any $t > 2$, we have $t-\frac32 \ge t/4$.
  \Cref{eq:lower-exp-gamma-bound} implies that when $b_j=0$, the absolute value of $\tilde{\delta}_j=\sigma_j\delta_j-\gamma_j\cdot \pi_j$ is smaller than $2$ since $\sigma_j \delta_j\in [-3/2,3/2]$ and \(|\gamma_j\pi_j|<1/2\).
  Therefore,
  \[
  \{|\widetilde\delta_j|\ge t\}
  \subseteq
\left\{
b_j=1,\,
|\gamma_j|\ge t - \frac32
\right\}
\subseteq
\left\{
b_j=1,\,
|\gamma_j|\ge\frac{t}{4}
\right\}.
\]
Therefore, for $t>2$, we have
  \begin{align*}
 \bbP(|\tilde{\delta}_j|\ge t)
 & \le \bbP\left(b_j=1, |\gamma_{j}|\ge \frac{t}{4} \right) \\
 & \le c\cdot \exp\left(-\left|\frac{\lambda_{\cal H}}{\eta_j}\right|^\al\right)\bbo\left\{c_2\cdot \frac{\lambda_{\cal H}}{|\eta_j|}\ge \frac{t}{4}\right\} \\
 & \le c\cdot \exp\left\{-\left(\frac{t}{4 c_2}\right)^\al\right\}.
  \end{align*}

If we pick a $c_2$ sufficiently small such that $c_2<\tau_\al/4$,
then for any $\tau\ge\tau_\al$, we have

  $$\bbP(|\tilde{\delta}_j|\ge t)\le c\cdot \exp\left\{-\left(\frac{t}{4 c_2}\right)^\al\right\}\le 2\exp\left\{-\left(\frac{t}{\tau}\right)^\al\right\}.$$

Define $F_0 = f_0^{\otimes d}$ and $F_1 = \bigotimes_{j=1}^d f_j$.
The tail bound for $\tilde{\delta}_j$ implies that  $F_1\in {\cal H}_{\al,\tau}^{\otimes}$.
Define $P_1=\bbP_{0,F_1}$ and $P_2=\bbP_{\mu,F_0}$.
To complete the proof using Lemma~\ref{lem: general tool}, it remains to control the total variation distance between $P_1$ and $P_2$ by $c c_2$ (up to a constant).
This is achieved by using the condition in \eqref{eq:exponential lower condition} and the following lemma.
{The next lemma is the point where the appendix condition \eqref{eq:exponential lower condition} enters: it bounds the aggregate Hellinger cost of the scale and location perturbations across coordinates.}

\begin{lemma}\label{lem: tv bound for exponential lower sparse}
In the proof of \Cref{eq:exponential lower asymmetric}, if $c_2\le1$, then
$$
{\rm TV}(P_1,P_2)
\le
\bar C^\prime c c_2,
$$
where $\bar C^\prime>0$ depends only on $\bar C$ in
\eqref{eq:exponential lower condition} and on the fixed density $f_0$.
\end{lemma}

The remaining step is to combine \Cref{lem: general tool} with \Cref{lem: tv bound for exponential lower sparse}.
Let $p_0:=1-\delta_0>0$. Choose $c_2$ such that
$$
0<c_2\le
\min\left\{
c_{2,1},\,1,\,\frac{\tau_\al}{8},\,
\frac{p_0}{4\bar C^\prime c}
\right\}.
$$
This choice satisfies all restrictions on $c_2$ used in the construction.
Furthermore, by \Cref{lem: tv bound for exponential lower sparse}, we have
\begin{equation}\label{eq:nonsym-tv-bound}
    {\rm TV}(P_1,P_2)
\le \bar C^\prime c c_2
\le \frac{p_0}{4}.
\end{equation}

If $\lambda_{\cal H}=0$, then \eqref{eq:exponential lower asymmetric} is
immediate.
Suppose that $\lambda_{\cal H}>0$ and set
$$
t:=\frac{c_1c_2}{2}s\lambda_{\cal H}>0.
$$
We will apply \Cref{lem: general tool} with $\Theta=\Theta_s$ with $t$ and
$$
\beta_1=\frac{p_0}{4},
\quad
\beta_2=\delta_0,
\quad
\beta_3=\frac{p_0}{4}, \quad\text{ and }\quad c=0.
$$
For $P_1=\bbP_{0,F_1}$, the prior is the point mass at zero, so
$$
\bbP\left(L(\theta)\le0,\ \theta\in\Theta_s\right)
=1\ge1-\beta_1.
$$
For $P_2=\bbP_{\mu,F_0}$, the alternative-event bound in
\eqref{eq: pro 1}--\eqref{eq: pro 2} gives
$$
\mu\left(L(\theta)\ge2t,\ \theta\in\Theta_s\right)
\ge 1-\delta_0
=1-\beta_2.
$$
The total variation bound in \eqref{eq:nonsym-tv-bound} is at most $\beta_3$. Since
$F_0,F_1\in\mathcal H_{\al,\tau}^{\otimes}$, \Cref{lem: general tool}
therefore yields
$$
\inf_{\hat L_s}
\sup_{\theta\in\Theta_s}
\sup_{\caP_\xi\in\mathcal H_{\al,\tau}^{\otimes}}
\bbP_{\theta,\caP_\xi}
\left(|\hat L_s-L(\theta)|\ge t\right)
\ge
\frac{1-\beta_1-\beta_2-\beta_3}{2}
=\frac{p_0}{4}.
$$
Consequently,
$$
\begin{aligned}
\inf_{\hat L_s}
\sup_{\theta\in\Theta_s}
\sup_{\caP_\xi\in\mathcal H_{\al,\tau}^{\otimes}}
\bbE_{\theta,\caP_\xi}
\left(\hat L_s-L(\theta)\right)^2
\,\ge\, t^2\frac{p_0}{4} \,=\, \frac{p_0c_1^2c_2^2}{16}
s^2\lambda_{\cal H}^2.
\end{aligned}
$$
This completes the proof for \Cref{eq:exponential lower asymmetric} with
$c^\prime=p_0c_1^2c_2^2/16>0$.

Finally, by increasing $\tau_\al$ if necessary, we may take the same
$\tau_\al$ to satisfy the requirement of
\Cref{thm: nonadaptive linear}.
This enlargement does not affect the
argument above, which only requires
$\tau_\al\ge 2/\log^{1/\al}2$.
Consequently, \Cref{eq:exponential lower symmetric} also holds.

\subsection{Proof of Lemma~\ref{lem: chi square}}
\begin{proof}
{We work with the generalized Gaussian density $f_\al(x;\sigma)\propto \exp(-|x|^\al/\sigma^\al)$ because it gives an explicit one-parameter family whose shift behavior can still be controlled sharply enough for the lower bound. The proof first normalizes the variance, then checks the exponential-tail condition, and finally bounds the shifted $\chi^2$ divergence.}
  Consider the scale family $\{F_\al(\sigma)\}$ on $\mathbb{R}$ with density
\[
f_\al(x;\sigma)=\frac{1}{c_\sigma}\exp\!\left(-\frac{|x|^\al}{\sigma^\al}\right),\qquad \al>0,\ \sigma>0,
\]
where the normalizing constant is
\[
c_\sigma=\int_{\mathbb{R}} \exp\!\left(-\frac{|x|^\al}{\sigma^\al}\right)\,\diff x
=\sigma\!\int_{\mathbb{R}} e^{-|y|^\al}\,\diff y
=\sigma\cdot \frac{2}{\al}\,\Gamma\!\left(\frac{1}{\al}\right).
\]

We choose $\sigma=\sigma_\al$ so that $\mathbb{E}_{F_\al(\sigma_\al)}[X^2]=1$. Using the change of
variables $x=\sigma y$,
\[
\int_{\mathbb{R}} x^2 e^{-|x|^\al/\sigma^\al}\,\diff x
=\sigma^3 \int_{\mathbb{R}} y^2 e^{-|y|^\al}\,\diff y
=\sigma^3\cdot \frac{2}{\al}\,\Gamma\!\left(\frac{3}{\al}\right),
\]
and
\[
\int_{\mathbb{R}} e^{-|x|^\al/\sigma^\al}\,\diff x
=\sigma \int_{\mathbb{R}} e^{-|y|^\al}\,\diff y
=\sigma\cdot \frac{2}{\al}\,\Gamma\!\left(\frac{1}{\al}\right).
\]
The condition $\mathbb{E}[X^2]=1$ is equivalent to
\[
\frac{1}{c_\sigma}\int_{\mathbb{R}} x^2 e^{-|x|^\al/\sigma^\al}\,\diff x =1
\quad\Longleftrightarrow\quad
\int_{\mathbb{R}} x^2 e^{-|x|^\al/\sigma^\al}\,\diff x
=\int_{\mathbb{R}} e^{-|x|^\al/\sigma^\al}\,\diff x,
\]
which yields
\[
\sigma_\al^2
=\frac{\displaystyle \int_{\mathbb{R}} e^{-|y|^\al}\,\diff y}
       {\displaystyle \int_{\mathbb{R}} y^2 e^{-|y|^\al}\,\diff y}
=\frac{\frac{2}{\al}\Gamma(1/\al)}{\frac{2}{\al}\Gamma(3/\al)}
=\frac{\Gamma(1/\al)}{\Gamma(3/\al)}.
\]

Hence, taking
\[
\sigma_\al=\sqrt{\frac{\Gamma(1/\al)}{\Gamma(3/\al)}}
\]
makes {$F_\al(\sigma_\al)$} mean-zero (by symmetry) with unit variance.
In order to satisfy the tail decay condition, $\tau$ is required to satisfy
\begin{equation}\label{eq:lem chi sq mgf bound}
\mathbb{E}_{F_\al(\sigma_\al)}\exp\left(\frac{|X|^\al}{\tau^\al}\right)=\int_{\mathbb{R}}\frac{1}{c_{\sigma_\al}}\exp\left[\left(\frac{1}{\tau^\al}-\frac{1}{\sigma_\al^\al}\right) |x|^\al\right] \diff x \le 2.
\end{equation}
The existence of such $\tau$ is guaranteed by the fact that the above integral is {decreasing} in $\tau$ and when $\tau\to\infty$, the integral converges to 1.
In fact, if we choose $\bar{\tau}_\al=\sigma_\al(1-2^{-\al})^{-1/\al}$ and define $\tau_\al=2^{1/\al}\bar{\tau}_\al$, then for any $\tau\ge \bar{\tau}_\al$, \Cref{eq:lem chi sq mgf bound} holds; we defer the derivation of $\bar{\tau}_\al$ to the end of this proof.
Using Markov's inequality, we see that $F_\al(\sigma_\al)\in {\cal G}_{\al,\tau}$ for all $\tau\ge \tau_\al$.

 {The remaining task is to bound} the $\chi^2$ divergence between {$f_\al^{(1)}(\cdot)=f_\al(\cdot-\gamma;\sigma_\al)$ and $f_\al^{(0)}(\cdot)=f_\al(\cdot;\sigma_\al)$.}
We have
\begin{align*}
  {1+\chi^2(f_\al^{(1)}\|f_\al^{(0)})=\int_{\mathbb{R}} \frac{f_\al^{(1)}(x)^2}{f_\al^{(0)}(x)}\,\diff x}&=\frac{1}{c_{\sigma_\al}}\int_{\mathbb{R}} \exp\left(\frac{|x|^\al}{\sigma_\al^\al}-\frac{2|x-\gamma|^\al}{\sigma_\al^\al}\right)\,\diff x.
\end{align*}
\begin{itemize}
  \item When $0<\al\le 1$, we have $|x|^\al\le |x-\gamma|^\al+|\gamma|^\al$ and therefore, we have
  \begin{align*}
    {1+\chi^2(f_\al^{(1)}\|f_\al^{(0)})}&=\frac{1}{c_{\sigma_\al}}\int_{\mathbb{R}} \exp\left(\frac{2|x|^\al}{\sigma_\al^\al}-\frac{2|x-\gamma|^\al}{\sigma_\al^\al}\right)\exp\left(-\frac{|x|^\al}{\sigma_\al^\al}\right)\,\diff x.\\
    &\le \frac{1}{c_{\sigma_\al}}\int_{\mathbb{R}} \exp\left(\frac{2|\gamma|^\al}{\sigma_\al^\al}\right)\exp\left(-\frac{|x|^\al}{\sigma_\al^\al}\right)\,\diff x\le \exp\left(\frac{2|\gamma|^\al}{\sigma_\al^\al}\right)\\
  \end{align*}
  \item {When $\al>1$, the map $x\mapsto |x|^\al$ is convex. Hence, for any $s\in(0,1)$,
\[
|x|^\al=\Big|(1-s)\frac{x-\gamma}{1-s}+s\frac{\gamma}{s}\Big|^\al
\le (1-s)\Big|\frac{x-\gamma}{1-s}\Big|^\al + s\Big|\frac{\gamma}{s}\Big|^\al
= (1-s)^{1-\al}|x-\gamma|^\al + s^{1-\al}|\gamma|^\al.
\]}
{
For any $s>0$ such that $(1-s)^{1-\al}<2$, we have
}

  \begin{align*}
    {1+\chi^2(f_\al^{(1)}\|f_\al^{(0)})}&=\frac{1}{c_{\sigma_\al}}\int_{\mathbb{R}} \exp\left(\frac{|x|^\al}{\sigma_\al^\al}-\frac{2|x-\gamma|^\al}{\sigma_\al^\al}\right)\,\diff x.\\
    &\le \frac{1}{c_{\sigma_\al}}\int_{\mathbb{R}} \exp\left(\frac{s^{1-\al}|\gamma|^\al}{\sigma_\al^\al}\right)\exp\left(-\frac{(2-(1-s)^{1-\al})|x|^\al}{\sigma_\al^\al}\right)\,\diff x\\
    &\le \exp\left(\frac{s^{1-\al}|\gamma|^\al}{\sigma_\al^\al}\right)\cdot \frac{1}{c_{\sigma_\al}}\int_{\mathbb{R}}\exp\left(-\frac{(2-(1-s)^{1-\al})|x|^\al}{\sigma_\al^\al}\right)\,\diff x,
  \end{align*}
where we have rewritten the variate $x-\gamma$ as $x$ after the first inequality.
  We can choose $s=s_\al\in (0,1)$ such that {$2-(1-s_\al)^{1-\al}>0$ and thus}

  $$\int_{\mathbb{R}}\exp\left(-\frac{(2-(1-s_\al)^{1-\al})|x|^\al}{\sigma_\al^\al}\right)\,\diff x<\infty.$$
  And therefore, we have
  $$1+\chi^2(f_\al^{(1)}\|f_\al^{(0)})\le \frac{1}{c_{\sigma_\al}}\int_{\mathbb{R}}\exp\left(-\frac{(2-(1-s_\al)^{1-\al})|x|^\al}{\sigma_\al^\al}\right)\,\diff x\cdot \exp\left(\frac{s_\al^{1-\al}|\gamma|^\al}{\sigma_\al^\al}\right).$$
\end{itemize}
In particular, when $\al=2$, we can obtain the sharper constant in the lemma.
Indeed, $\sigma_2^2=\Gamma(1/2)/\Gamma(3/2)=2$, so
$f_2^{(0)}$ is the standard Gaussian density. If
$X\sim f_2^{(0)}$, then
$$
\frac{f_2^{(1)}(X)}{f_2^{(0)}(X)}
=\exp\left(\gamma X-\frac{\gamma^2}{2}\right),
$$
and hence
$1+\chi^2(f_2^{(1)}\|f_2^{(0)})
=\bbE\exp(2\gamma X-\gamma^2)
=\exp(\gamma^2)$.
We can take $C_{2,1}=C_{2,2}=1$.
\end{proof}

\begin{proof}[Necessary and sufficient condition for \Cref{eq:lem chi sq mgf bound}]
We prove that \Cref{eq:lem chi sq mgf bound} holds if and only if $\tau \ge \sigma_\al(1-2^{-\al})^{-1/\al}$.
In the following, the expectation on the left hand side of \Cref{eq:lem chi sq mgf bound} is denoted as $I(\tau)$.

It is clear that $\tau>\sigma_\al$ is needed for \Cref{eq:lem chi sq mgf bound} to hold.
For such a $\tau$, define $a(\tau)=\frac{1}{\sigma_\al^\al}-\frac{1}{\tau^\al}>0$. We have
$$
\begin{aligned}
I(\tau)=&\int_{\mathbb{R}} \frac{1}{c_{\sigma_\al}} \exp \left[\left(\frac{1}{\tau^\al}-\frac{1}{\sigma_\al^\al}\right)|x|^\al\right] \diff x \\
=& \frac{2}{c_{\sigma_\al}} \int_0^{\infty} e^{-a(\tau) x^\al} \diff x \\
=& \frac{2}{c_{\sigma_\al}}\frac{1}{\al} \int_0^{\infty} y^{1/\al -1} e^{-a(\tau) y} \diff y \\
=& \frac{2}{c_{\sigma_\al}} \frac{1}{\al} \Gamma(1/\al) [a(\tau)]^{-1/\al}.
\end{aligned}
$$
Since $c_{\sigma_\al}=\sigma_\al \cdot \frac{2}{\al} \Gamma(1 / \al)$, we have
$$
\begin{aligned}
I(\tau)&=\frac{1}{\sigma_\al}\left(\frac{1}{\sigma_\al^\al}-\frac{1}{\tau^\al}\right)^{-1 / \al} \\
&=\left(1-\frac{\sigma_\al^\al}{\tau^\al}\right)^{-1 / \al}.
\end{aligned}
$$
Therefore, it is straightforward to see that
$$I(\tau)\le 2 \Leftrightarrow 1-\frac{\sigma_\al^\al}{\tau^\al} \ge 2^{-\al} \Leftrightarrow \tau\ge \sigma_\al (1-2^{-\al})^{-1/\al}.$$
This completes the proof.

\end{proof}

\subsection{Proof of \Cref{lem:adaptive-prior}}

The proof consists of four parts. We first construct a random sparsity
prior with capped Bernoulli probabilities and bound its uncapped weights.
We then control the random support size and
the value of $L(\theta)$ separately.
Finally, we analyze the chi-square divergence of the resulting product mixture.

\textit{Part 1: Prior construction.}
Write $\ell_s=\log(es)$ and $\bar\beta=(\beta_*(s))_+$. By the definitions
of $\lambda_*(s)$ and $\nu_*(s)$ in \Cref{sec:adaptive minimax},
\begin{equation}\label{eq:adaptive-prior-normalizer}
\bar\beta=\lambda_*^\al(s),
\qquad
D^2:=\sum_{j=1}^d\eta_j^2
\exp\left(-\frac{\bar\beta}{|\eta_j|^\al}\right)
=\frac{\nu_*^2(s)}{\ell_s}>0.
\end{equation}

Let $C_{\al,\chi}=C_{\al,1}(e+1)$, where $C_{\al,1}$ is the constant in
\Cref{lem: chi square}, and choose $c_1\in(0,1)$ such that
\begin{equation}\label{eq:adaptive-prior-c1-choice}
C_{\al,\chi}c_1^2\le\frac{\gamma}{2}.
\end{equation}
Define
\begin{equation}\label{eq:adaptive-prior-definition}
\begin{aligned}
r_j
&=c_1\sqrt{\ell_s}\,
\frac{|\eta_j|\exp(-\bar\beta/|\eta_j|^\al)}{D},
\qquad
\pi_j=1\wedge r_j,\\
\gamma_j
&=\begin{cases}
C_{\al,2}{\rm sign}(\eta_j),&j\le j_2(s),\\
C_{\al,2}\lambda_*(s)/\eta_j,&j>j_2(s),
\end{cases}
\end{aligned}
\end{equation}
where $C_{\al,2}$ is the other constant in \Cref{lem: chi square}.
Let $b_1,\ldots,b_d$ be independent random variables with
$b_j\sim{\rm Ber}(\pi_j)$, set $\theta_j=b_j\gamma_j$, and denote the
law of $\theta$ by $\mu$. Since $r_j\ge0$, the definition
$\pi_j=1\wedge r_j$ ensures that $0\le\pi_j\le1$.
Moreover, $e^{-2\bar\beta/|\eta_j|^\al}
\le e^{-\bar\beta/|\eta_j|^\al}$ since $\bar\beta\ge0$. By definition of $D^2$, we have
\begin{equation}\label{eq:adaptive-prior-square-sum}
\begin{aligned}
\sum_{j=1}^dr_j^2
&=\frac{c_1^2\ell_s}{D^2}
\sum_{j=1}^d\eta_j^2
\exp\left(-\frac{2\bar\beta}{|\eta_j|^\al}\right)
&\le\frac{c_1^2\ell_s}{D^2}
\sum_{j=1}^d\eta_j^2
\exp\left(-\frac{\bar\beta}{|\eta_j|^\al}\right)
=c_1^2\ell_s.
\end{aligned}
\end{equation}

\textit{Part 2: Control the support size.}
The support size $N:=\|\theta\|_0$ is random under $\mu$, so the prior is
not supported on $\Theta_s$ deterministically. We show instead that
$N\le s$ with probability at least $7/8$.

Recall $\phi(t)$ in \Cref{lem: decrease}.
By definition in \eqref{eq:solution linear adaptive}, we have $
\phi(\bar\beta)
\le\phi(\beta_*(s))
=\frac{s}{2\sqrt{\ell_s}}$.
It follows
from \eqref{eq:adaptive-prior-definition} that
\begin{equation}\label{eq:adaptive-prior-support-mean}
\begin{aligned}
\sum_{j=1}^d\pi_j
\le\sum_{j=1}^dr_j
&=c_1\sqrt{\ell_s}\,
\frac{\sum_j|\eta_j|e^{-\bar\beta/|\eta_j|^\al}}{D}\\
&=c_1\sqrt{\ell_s}\,\phi(\bar\beta)\\
&\le\frac{c_1s}{2}.
\end{aligned}
\end{equation}
Every $\gamma_j$ in \eqref{eq:adaptive-prior-definition} is nonzero:
this is immediate when $\lambda_*(s)>0$, while $\lambda_*(s)=0$ implies
$j_2(s)=d$.
Therefore, we have $N=\|\theta\|_0=\sum_{j=1}^db_j$. Independence of the Bernoulli variables
and \eqref{eq:adaptive-prior-support-mean} give
\begin{align*}
\bbE_\mu N
&=\sum_{j=1}^d\pi_j\le\frac{c_1s}{2},\\
{\rm Var}_\mu(N)
&=\sum_{j=1}^d\pi_j(1-\pi_j)
\le\sum_{j=1}^d\pi_j
\le\frac{c_1s}{2}.
\end{align*}
Since $c_1<1$, we also have
$s-\bbE_\mu N\ge s(1-c_1/2)\ge s/2$. Chebyshev's inequality therefore
gives
\begin{equation}\label{eq:adaptive-prior-support-tail}
\begin{aligned}
\mu(N>s)
&\le\frac{{\rm Var}_\mu(N)}{(s-\bbE_\mu N)^2}\\
&\le\frac{c_1s/2}{(s/2)^2}
=\frac{2c_1}{s}
\le\frac18,
\end{aligned}
\end{equation}
where the last inequality follows from $c_1<1$ and $s\ge K_1$ if we choose $K_1\ge16$.

\textit{Part 3: Control the separation.}
We next lower-bound $L(\theta)$ under the capped prior. Since
\[
\bbE_\mu L(\theta)=\sum_{j=1}^d w_j\pi_j,
\qquad \pi_j=1\wedge r_j,
\]
we first lower-bound the more explicit uncapped sum
$S_r:=\sum_jw_jr_j$ and then control the loss
$S_r-\bbE_\mu L(\theta)$ caused by capping the probabilities.

Define $M_s=\max\{|\eta_1|,\lambda_*(s)\}$ and for each $j$,
$w_j:=\eta_j\gamma_j
=C_{\al,2}\max\{|\eta_j|,\lambda_*(s)\}$.
We have
\begin{equation}\label{eq:adaptive-prior-bound-wj}
    0<w_j\le C_{\al,2}M_s.
\end{equation}

Since $w_j\ge C_{\al,2}|\eta_j|$, the definitions of $r_j$ and
$D$ yield
\begin{equation}\label{eq:adaptive-prior-uncapped-nu}
\begin{aligned}
S_r
&=\frac{c_1\sqrt{\ell_s}}{D}
\sum_{j=1}^dw_j|\eta_j|
\exp\left(-\frac{\bar\beta}{|\eta_j|^\al}\right)\\
&\ge\frac{C_{\al,2}c_1\sqrt{\ell_s}}{D}
\sum_{j=1}^d\eta_j^2
\exp\left(-\frac{\bar\beta}{|\eta_j|^\al}\right)\\
&=C_{\al,2}c_1\sqrt{\ell_s}D
=C_{\al,2}c_1\nu_*(s).
\end{aligned}
\end{equation}

If $\lambda_*(s)>0$, then $\bar\beta=\beta_*(s)$ and $\phi(\bar\beta)=\frac{s}{2\sqrt{\ell_s}}$, so
\eqref{eq:solution linear adaptive} and
\eqref{eq:adaptive-prior-definition} yield
$\sum_jr_j=c_1s/2$.
Since $w_j\ge C_{\al,2}\lambda_*(s)$ for every $j$, we have
\begin{equation}\label{eq:adaptive-prior-uncapped-lambda}
\sum_{j=1}^dw_jr_j
\ge C_{\al,2}\lambda_*(s)\sum_{j=1}^dr_j
=\frac{C_{\al,2}c_1}{2}s\lambda_*(s).
\end{equation}
When $\lambda_*(s)=0$,
\eqref{eq:adaptive-prior-uncapped-lambda} remains valid since the left hand side is nonnegative.

Combining
\eqref{eq:adaptive-prior-uncapped-nu} and
\eqref{eq:adaptive-prior-uncapped-lambda}, and using
$\max\{a,b\}\ge(a+b)/2$ for $a,b\ge0$, we obtain
\begin{equation}\label{eq:adaptive-prior-uncapped-separation}
\begin{aligned}
S_r
&\ge\frac12\left\{C_{\al,2}c_1\nu_*(s)
+\frac{C_{\al,2}c_1}{2}s\lambda_*(s)\right\}\\
&\ge\frac{C_{\al,2}c_1}{4}
\{\nu_*(s)+s\lambda_*(s)\}
=\frac{C_{\al,2}c_1}{4}A_s.
\end{aligned}
\end{equation}

Let $H=\{j:r_j>1\}$. By
\eqref{eq:adaptive-prior-square-sum}, we have
$|H|
\le\sum_{j\in H}r_j^2
\le c_1^2\ell_s$.
Then by the Cauchy--Schwarz inequality and
and \eqref{eq:adaptive-prior-square-sum}, we have
\begin{equation}\label{eq:adaptive-prior-capped-mass}
\sum_{j\in H}r_j
\le |H|^{1/2}
\left(\sum_{j\in H}r_j^2\right)^{1/2}
\le c_1^2\ell_s.
\end{equation}
For $j\notin H$, we have $\pi_j=r_j$, whereas for $j\in H$, we have
$\pi_j=1$.
The difference between the uncapped and the capped sums is
$$
S_r-\bbE_\mu L(\theta)
=\sum_{j\in H}w_j(r_j-1)
\leq C_{\al,2}M_s\sum_{j\in H}r_j
\leq C_{\al,2}M_s c_1^2\ell_s,
$$
where the first inequality follows from \eqref{eq:adaptive-prior-bound-wj} and the second from
\eqref{eq:adaptive-prior-capped-mass}.
By \eqref{eq:adaptive-prior-uncapped-separation}, we obtain
\begin{equation}\label{eq:adaptive-prior-capped-mean}
\begin{aligned}
\bbE_\mu L(\theta)
&\ge S_r-C_{\al,2}M_s c_1^2\ell_s\\
&\ge\frac{C_{\al,2}c_1}{4}A_s
-C_{\al,2}c_1^2\ell_sM_s.
\end{aligned}
\end{equation}
Choose
\begin{equation}\label{eq:adaptive-prior-K2-choice}
K_2\ge 8c_1,\quad \text{ and }\quad K_2 \ge \frac{256}{c_1}.
\end{equation}
Under the hypothesis $A_s\ge K_2\ell_sM_s$ of the lemma, the first
condition in \eqref{eq:adaptive-prior-K2-choice} gives
$$
c_1^2\ell_sM_s
\le\frac{c_1^2}{K_2}A_s
\le\frac{c_1}{8}A_s.
$$
Thus
\eqref{eq:adaptive-prior-capped-mean} yields
\begin{equation}\label{eq:adaptive-prior-mean-lower}
\bbE_\mu L(\theta)
\ge\frac{C_{\al,2}c_1}{8}A_s.
\end{equation}
Using \eqref{eq:adaptive-prior-bound-wj} and the independence of the $b_j$'s, we have
\begin{equation}\label{eq:adaptive-prior-L-variance}
\begin{aligned}
{\rm Var}_\mu L(\theta)
&=\sum_{j=1}^dw_j^2\pi_j(1-\pi_j)\\
&\le C_{\al,2}M_s\sum_{j=1}^dw_j\pi_j
=C_{\al,2}M_s\bbE_\mu L(\theta).
\end{aligned}
\end{equation}
Letting $c_0=C_{\al,2}c_1/16$, \eqref{eq:adaptive-prior-mean-lower} becomes $c_0A_s\le\bbE_\mu L(\theta)/2$.
Therefore,
$$
\{L(\theta)<c_0A_s\}
\subseteq
\left\{\bbE_\mu L(\theta)-L(\theta)
>\frac12\bbE_\mu L(\theta)\right\}.
$$
By Chebyshev's inequality, we have
\begin{equation}\label{eq:adaptive-prior-separation-tail}
\begin{aligned}
\mu\{L(\theta)<c_0A_s\}
&\le\frac{4{\rm Var}_\mu L(\theta)}
{[\bbE_\mu L(\theta)]^2}\\
(\because \text{by \eqref{eq:adaptive-prior-L-variance}})\qquad &\le\frac{4C_{\al,2}M_s}{\bbE_\mu L(\theta)}\\
(\because \text{by \eqref{eq:adaptive-prior-mean-lower}})\qquad &\le\frac{32M_s}{c_1A_s}\\
(\because A_s\ge K_2\ell_sM_s)\qquad &\le\frac{32}{c_1K_2\ell_s}\\
&\le\frac{32}{c_1K_2}\\
&\le\frac18,
\end{aligned}
\end{equation}
where the second last inequality uses $\ell_s=\log(es)\ge1$, and the last inequality follows from the
second condition in \eqref{eq:adaptive-prior-K2-choice}.

\medskip
We now verify the \textbf{first
conclusion} of the lemma by
combining the bounds on the support size and the functional separation.
The bounds
\eqref{eq:adaptive-prior-support-tail},
\eqref{eq:adaptive-prior-separation-tail}, and the union bound give
\begin{align*}
\mu\{\theta\in\Theta_s:L(\theta)\ge c_0A_s\}
&=\mu\{N\le s,L(\theta)\ge c_0A_s\}\\
&\ge1-\mu(N>s)-\mu\{L(\theta)<c_0A_s\}\\
&\ge1-\frac18-\frac18=\frac34.
\end{align*}

\textit{Part 4: Control the chi-square divergence.}
It remains to control the chi-square divergence of the observation
distribution induced by $\mu$.
We first calculate the divergence for each
coordinate and then combine the coordinatewise bounds.

Let $f_\al^{(0)}$ be the marginal density defining
$F_\al^\otimes$ in \Cref{sec:proof nonadaptive linear}. Under
$\bbP_{\mu,F_\al^\otimes}$, the $j$th observation has density
$$
g_j=(1-\pi_j)f_\al^{(0)}
+\pi_jf_\al^{(0)}(\cdot-\gamma_j).
$$
Define
$\Delta_j=\chi^2(f_\al^{(0)}(\cdot-\gamma_j)
\|f_\al^{(0)})$. A direct calculation gives
\begin{equation}\label{eq:adaptive-prior-coordinate-mixture}
\begin{aligned}
\chi^2(g_j\|f_\al^{(0)})
&=\int_{\mathbb R}
\left(\frac{g_j(x)}{f_\al^{(0)}(x)}-1\right)^2
f_\al^{(0)}(x)\,\diff x\\
&=\pi_j^2\int_{\mathbb R}
\left(\frac{f_\al^{(0)}(x-\gamma_j)}{f_\al^{(0)}(x)}-1\right)^2
f_\al^{(0)}(x)\,\diff x\\
&=\pi_j^2\Delta_j.
\end{aligned}
\end{equation}

We next bound $\Delta_j$.
By \Cref{lem: chi square},
$\Delta_j
\le 1+\Delta_j
\le C_{\al,1}\exp\left(\left|\frac{\gamma_j}{C_{\al,2}}\right|^\al\right)$. Recall the definition of $\gamma_j$ in \eqref{eq:adaptive-prior-definition}.
For $j\le j_2(s)$, $|\gamma_j/C_{\al,2}|^\al=1$. For $j>j_2(s)$,
\eqref{eq:adaptive-prior-definition} and
\eqref{eq:adaptive-prior-normalizer} give
$|\gamma_j/C_{\al,2}|^\al
=\bar\beta/|\eta_j|^\al$. Consequently,
\begin{equation}\label{eq:adaptive-prior-coordinate-divergence}
\Delta_j\le
\begin{cases}
C_{\al,1}e,&j\le j_2(s),\\
C_{\al,1}\exp\left(\bar\beta/|\eta_j|^\al\right),&j>j_2(s).
\end{cases}
\end{equation}

Using the definition of $r_j$ in \eqref{eq:adaptive-prior-definition}, we have
$$
r_j^2\Delta_j
\le\frac{C_{\al,1}c_1^2\ell_s}{D^2}
\begin{cases}
e\eta_j^2e^{-2\bar\beta/|\eta_j|^\al},&j\le j_2(s),\\
\eta_j^2e^{-\bar\beta/|\eta_j|^\al},&j>j_2(s).
\end{cases}
$$
Since $\bar\beta\ge0$,
$e^{-2\bar\beta/|\eta_j|^\al}
\le e^{-\bar\beta/|\eta_j|^\al}$.
Summing up these
coordinatewise bounds, we obtain
\begin{align*}
&e\sum_{j\le j_2(s)}\eta_j^2
e^{-2\bar\beta/|\eta_j|^\al}
+\sum_{j>j_2(s)}\eta_j^2
e^{-\bar\beta/|\eta_j|^\al}\\
&\qquad\le(e+1)
\sum_{j=1}^d\eta_j^2e^{-\bar\beta/|\eta_j|^\al}
=(e+1)D^2,
\end{align*}
where the equality follows from the definition of $D^2$ in
\eqref{eq:adaptive-prior-normalizer}.

Since $\pi_j\le r_j$, we have
\begin{equation}\label{eq:adaptive-prior-summed-divergence}
\begin{aligned}
\sum_{j=1}^d\pi_j^2\Delta_j
&\le\sum_{j=1}^dr_j^2\Delta_j\\
&\le\frac{C_{\al,1}c_1^2\ell_s}{D^2}
\left[
e\sum_{j\le j_2(s)}\eta_j^2
e^{-2\bar\beta/|\eta_j|^\al}
+\sum_{j>j_2(s)}\eta_j^2
e^{-\bar\beta/|\eta_j|^\al}
\right]\\
&\le C_{\al,\chi}c_1^2\ell_s,
\end{aligned}
\end{equation}
where $C_{\al,\chi}=C_{\al,1}(e+1)$.

The densities of
$\bbP_{\mu,F_\al^\otimes}$ and $\bbP_{0,F_\al^\otimes}$ are
$\prod_jg_j$ and $\prod_jf_\al^{(0)}$, respectively.
We combine the coordinatewise divergences into the product divergence as follows:
\begin{equation}\label{eq:adaptive-prior-product-divergence}
\begin{aligned}
1+\chi^2\bigl(\bbP_{\mu,F_\al^\otimes}
\,\|\,\bbP_{0,F_\al^\otimes}\bigr)
&=\int_{\mathbb R^d}
\prod_{j=1}^d\frac{g_j^2(x_j)}{f_\al^{(0)}(x_j)}
\,\diff x_1\cdots\diff x_d\\
(\because \text{Tonelli's theorem})\qquad &=\prod_{j=1}^d
\int_{\mathbb R}\frac{g_j^2(x)}{f_\al^{(0)}(x)}\,\diff x\\
(\because \text{\Cref{eq:adaptive-prior-coordinate-mixture}})\qquad &=\prod_{j=1}^d(1+\pi_j^2\Delta_j)\\
(\because ~ 1+x\le e^x\text{ for }x\ge0)\qquad &\le\exp\left\{\sum_{j=1}^d\pi_j^2\Delta_j\right\}\\
(\because \text{inequality in \eqref{eq:adaptive-prior-summed-divergence}})\qquad &\le\exp\{C_{\al,\chi}c_1^2\ell_s\}.
\end{aligned}
\end{equation}

It remains to show that $\exp\{C_{\al,\chi}c_1^2\ell_s\}$ is at most $s^\gamma$.
The choice in \eqref{eq:adaptive-prior-c1-choice} ensures that
$2C_{\al,\chi}c_1^2\le\gamma$.
Since $\ell_s=\log(es)$ and $s\ge K_1\ge16>e$, we have
\begin{align*}
\exp\{C_{\al,\chi}c_1^2\ell_s\}
&=e^{C_{\al,\chi}c_1^2}s^{C_{\al,\chi}c_1^2}\\
&\le s^{2C_{\al,\chi}c_1^2}\\
&\le s^\gamma.
\end{align*}
Combining this bound with
\eqref{eq:adaptive-prior-product-divergence} proves the \textbf{second conclusion}
of the lemma.

This completes the proof.

\subsection{Additional proof of the lower bounds in ~\Cref{thm: nonadaptive linear}}\label{sec: add proof lower}
\begin{proof}
We now prove the lower bound in \Cref{thm: nonadaptive linear} when
\[
s \le C_1
\quad \text{or} \quad
\lambda_{o} s + \nu \le C_2 |\eta_1|,
\]
for some constants $C_1, C_2 > 0$ and for any $\al > 0$, $\tau \ge \tau_\al$.

\textit{Case 1: $\lambda_{o} s + \nu \le C_2 |\eta_1|$.}
We need to show that
\[
\inf_{\hat{T}}\sup_{\theta\in\Theta_s}\sup_{{\caP}_\xi \in {\cal G}_{\al,\tau}^\otimes}
\mathbb{E}_{\theta,{\caP}_\xi}\!\left(|\hat{T}-L(\theta)|^2\right)
\;\gtrsim\; \eta_1^2.
\]
Let $c > 0$ be a small constant to be specified later, and let $\mu$ be the point mass at
$(c \,\mathrm{sgn}(\eta_1), 0, \ldots, 0)$.
Recall the definition of $F_\al^\otimes$ in \Cref{sec:proof nonadaptive linear}, and set
\[
P_1 = \mathbb{P}_{0, F_\al^\otimes},
\qquad
P_2 = \mathbb{P}_{\mu, F_\al^\otimes}.
\]

Recall the density $f_\alpha^{(0)}$ in the construction of $F_\alpha^\otimes$. We claim that
\begin{equation}\label{eq:continuity-H-theorem1}
{\bf H}\!\left(f_\alpha^{(0)}(\cdot-c),f_\alpha^{(0)}\right)
\longrightarrow 0,
\qquad\text{as }c\to0.
\end{equation}
Indeed, write
$f_\alpha^{(0)}(x) =
K_\alpha \exp\left(-a_\alpha |x|^\alpha\right)$
for some constants \(K_\alpha,a_\alpha>0\). For every fixed \(x\in\mathbb R\), continuity of \(f_\alpha^{(0)}\) gives
\[
\left[
\sqrt{f_\alpha^{(0)}(x-c)}
-
\sqrt{f_\alpha^{(0)}(x)}
\right]^2
\longrightarrow0,
\qquad\text{as }c\to0.
\]
Moreover, for \(|c|\le1\), $|x-c|\ge (|x|-1)_+$. Therefore, we have
\[
\begin{aligned}
\left[
\sqrt{f_\alpha^{(0)}(x-c)}
-
\sqrt{f_\alpha^{(0)}(x)}
\right]^2
&\le
f_\alpha^{(0)}(x-c)+f_\alpha^{(0)}(x)\\
&\le
K_\alpha
\left[
\exp\left(-a_\alpha (|x|-1)_{+}^\alpha\right)
+
\exp\left(-a_\alpha |x|^\alpha\right)
\right],
\end{aligned}
\]
where the last expression is integrable over \(\mathbb R\). Hence, by the dominated convergence theorem,
\[
{\bf H}^2\!\left(
f_\alpha^{(0)}(\cdot-c),
f_\alpha^{(0)}
\right)
=
\int_{\mathbb R}
\left[
\sqrt{f_\alpha^{(0)}(x-c)}
-
\sqrt{f_\alpha^{(0)}(x)}
\right]^2dx
\longrightarrow0,
\qquad\text{as }c\to0.
\]
This proves \Cref{eq:continuity-H-theorem1}.
Consequently, for any fixed \(\alpha_1\in(0,1)\), we may choose \(c>0\) sufficiently small, depending only on \(\alpha\) and \(\alpha_1\), such that
\[
\mathrm{TV}(P_1,P_2)
\le
{\bf H}(P_1,P_2)
=
{\bf H}\!\left(
f_\alpha^{(0)}(\cdot-c),
f_\alpha^{(0)}
\right)
\le \alpha_1.
\]

Moreover,
\[
\mu\!\left(L(\theta) \ge c |\eta_1|, \theta\in\Theta_s\right) = 1.
\]
Thus, applying \Cref{lem: general tool} with $t = \frac12 c|\eta_1|$ yields
\[
\inf_{\hat{T}}\sup_{\theta\in \Theta_s}
\mathbb{P}\!\left(|\hat{T}-L(\theta)| \ge \frac12 c|\eta_1|\right)
\;\ge\; \frac{1-\al_1}{2} > 0.
\]

\textit{Case 2: $s \le C_1$.}
We first show that it suffices to prove that
\begin{equation}\label{eq:small s}
\inf_{\hat{T}}\sup_{\theta\in\Theta_s}
\sup_{{\caP}_\xi \in {\cal G}_{\al,\tau}^\otimes}
\mathbb{E}_{\theta,{\caP}_\xi}\!\left(|\hat{T}-L(\theta)|^2\right)
\;\gtrsim\; \lambda_{o}^2,
\qquad \text{if } \lambda_{o} \ge C_3 |\eta_1|
\end{equation}
for some sufficiently large constant $C_3 > 0$.

We already know that
\[
\inf_{\hat{T}}\sup_{\theta\in\Theta_s}
\sup_{{\caP}_\xi \in {\cal G}_{\al,\tau}^\otimes}
\mathbb{E}_{\theta,{\caP}_\xi}\!\left(|\hat{T}-L(\theta)|^2\right)
\;\gtrsim\; |\eta_1|^2.
\]
Note that $\lambda_{o}^2 s^2 \asymp \lambda_{o}^2$ when $s \le C_1$.
Define
\[
\pi_j = \frac{2}{s} \cdot \frac{|\eta_j| e^{-\beta_+/|\eta_j|^\al}}{\sqrt{\sum_{j=1}^d \eta_j^2 e^{-\beta_+/|\eta_j|^\al}}}.
\]
From \eqref{eq:solution nonadaptive}, we have $\sum_{j=1}^d \pi_j \le 1$. Then
\[
\nu^2 = \sum_{j=1}^d \eta_j^2 e^{-\beta_+ / |\eta_j|^\al}
= \sum_{j=1}^d \frac{s\pi_j}{2} |\eta_j| \nu
\;\le\; \frac{s|\eta_1|}{2}\nu
\;\;\Rightarrow\;\; \nu \lesssim |\eta_1|.
\]
Hence, if $\lambda_{o} \lesssim |\eta_1|$, we have
$\lambda_{o}^2 s^2 + \nu^2 \lesssim |\eta_1|^2$.

Therefore, it remains to prove \eqref{eq:small s}.

To this end, we have $\lambda_{o} \ge C_3 |\eta_1|>0$, and thus $\beta_+=\beta>0$ and $\sum_{j=1}^d \pi_j=1$.

Let $J$ be a random index taking values in $[d]$ such that
$\mathbb{P}(J=j)=\pi_j$ for every $j\in[d]$, and define
\[
\theta_j
=
c\frac{\lambda_o}{\eta_j}\bbo\{J=j\},
\qquad j=1,\ldots,d,
\]

where $c > 0$ is a constant to be specified later. Let $\mu$ denote the distribution of $\theta$. Then
\[
\mu\!\left(L(\theta) \ge c \lambda_{o}, \theta\in\Theta_s\right) = 1.
\]

We now compute the $\chi^2$-divergence between
$P_1 = \mathbb{P}_{0, F_\al^\otimes}$
and
$P_2 = \mathbb{P}_{\mu, F_\al^\otimes}$.
Conditioned on $p=j$, we have for $x=(x_1,\ldots,x_d)\in\mathbb{R}^d$:
\[
\mathbb{P}_{\mu,F_\al^\otimes}(x\mid p=j)
= f_\al^{(0)}\!\left(x_j - c\frac{\lambda_{o}}{\eta_j}\right)
\prod_{k\neq j} f_\al^{(0)}(x_k).
\]
Hence,
\begin{equation}\label{eq:nonadaptive-lower-bounded-sparsity-case2-divergence}
    \begin{aligned}
1+\chi^2(P_2\|P_1)
&= \int_{\mathbb{R}^d} \Bigg[\sum_{j=1}^d \pi_j
   \frac{f_\al^{(0)}\!\left(x_j - c\frac{\lambda_{o}}{\eta_j}\right)}
        {f_\al^{(0)}(x_j)}\Bigg]^2 \,\diff P_1 \\
&= \sum_{i\neq j} \pi_i\pi_j
    \int_{\mathbb{R}^2} f_\al^{(0)}\!\left(x_i - c\frac{\lambda_{o}}{\eta_i}\right)
                        f_\al^{(0)}\!\left(x_j - c\frac{\lambda_{o}}{\eta_j}\right)\,\diff x_i\,\diff x_j \\
&\quad + \sum_{j=1}^d \pi_j^2
    \int_{\mathbb{R}} \frac{f_\al^{(0)}\!\left(x_j - c\frac{\lambda_{o}}{\eta_j}\right)^2}{f_\al^{(0)}(x_j)}\,\diff x_j \\
&\le 1 + \sum_{j=1}^d \pi_j^2
   \Bigg[C_{\al,1} \exp\!\left(\Big|\frac{c}{C_{\al,2}} \cdot \frac{\lambda_{o}}{\eta_j}\Big|^\al\right)-1\Bigg] \\
&\le 1 + \frac{(2/s)^2 C_{\al,1} \sum_{j=1}^d \eta_j^2 \exp\!\left(\big((c/C_{\al,2})^\al-2\big)\beta_+/|\eta_j|^\al\right)}
             {\sum_{j=1}^d \eta_j^2 e^{-\beta_+/|\eta_j|^\al}}.
\end{aligned}
\end{equation}

Choose $c>0$ sufficiently small that
$\delta:=1-\left(\frac{c}{C_{\al,2}}\right)^\al>0$.
Set $t_j=\beta_+/|\eta_j|^\al$. Since
$\lambda_o=\beta_+^{1/\al}$ and
$|\eta_j|\le |\eta_1|$, the condition
$\lambda_o\ge C_3|\eta_1|$ implies that
$$
t_j
=\frac{\beta_+}{|\eta_j|^\al}
\ge \left(\frac{\lambda_o}{|\eta_1|}\right)^\al
\ge C_3^\al,
\qquad j\in[d].
$$
Consequently, \eqref{eq:nonadaptive-lower-bounded-sparsity-case2-divergence} yields
\begin{align*}
\chi^2(P_2\|P_1)
&\le \frac{4C_{\al,1}}{s^2}
\frac{\sum_{j=1}^d
\eta_j^2e^{-t_j}e^{-\delta t_j}}
{\sum_{j=1}^d\eta_j^2e^{-t_j}}\\
&\le \frac{4C_{\al,1}}{s^2}
e^{-\delta C_3^\al}
\le 4C_{\al,1}e^{-\delta C_3^\al},
\end{align*}
where the second inequality follows because the ratio is a weighted
average of $e^{-\delta t_j}$ with nonnegative weights
$\eta_j^2e^{-t_j}$, and the last inequality uses $s\ge1$.
For any fixed $\al_1\in(0,1)$, we may therefore choose $C_3>0$
sufficiently large that
$$
4C_{\al,1}e^{-\delta C_3^\al}\le4\al_1^2.
$$
It follows that
$$
\mathrm{TV}(P_1,P_2)
\le\frac12\sqrt{\chi^2(P_2\|P_1)}
\le\al_1.
$$
Applying \Cref{lem: general tool} with
$t=\frac12c\lambda_o$ completes the proof.

\end{proof}

\subsection{Proof of the lower bound in \Cref{thm: testing}}\label{sec:test lower}
Based on the following lemma, the proof of the lower bound is essentially the same as that of \Cref{thm: nonadaptive linear}.
\begin{lemma}[{\cite[Lemma 3]{collier2017minimax}}]\label{lem: test}
If $\mu_1$ and $\mu_2$ are probability measures supported on
$\Theta_{s,0}$ and $\Theta_s(\rho)$, respectively, then
\[
\inf_{\Delta} \left\{ \sup_{\theta\in \Theta_{s,0}}\bbP_{\theta,{\cal N}^\otimes}(\Delta = 1) + \sup_{\theta \in \Theta_s(\rho)} \bbP_{\theta,{\caN}^\otimes}(\Delta = 0) \right\}
\;\ge\; 1 - {\rm TV}(\mathbb{P}_{\mu_1,{\cal N}^\otimes}, \mathbb{P}_{\mu_2,{\cal N}^\otimes}),
\]
where $\inf_\Delta$ denotes the infimum over all $\{0,1\}$-valued statistics.
\end{lemma}
Since $s \ge 2$, we can select some $\theta_0 \in \Theta_{s,0}$ such that
$\|\theta_0\|_0 \le \lfloor s/2 \rfloor$.
We then let $\mu_1$ be the point mass at $\theta_0$,
and construct $\mu_2$ as the distribution obtained by shifting $\theta_0$
with the least favorable alternatives used in the proof of
\Cref{thm: nonadaptive linear}.

\subsection{Proofs of lemmas in \Cref{subsec:proof exponential lower}}

\begin{proof}[Proof of \Cref{lem: tv bound for exponential lower sparse}]

Write
$P_1=\widetilde P_1*\mu$ and $P_2=\widetilde P_2*\mu$, where
$\widetilde P_1$ is the distribution of
$(\sigma_1\delta_1-\pi_1\gamma_1,\ldots,
\sigma_d\delta_d-\pi_d\gamma_d)$ and $\widetilde P_2$ is the distribution of
$(\delta_1,\ldots,\delta_d)$. The contraction of total variation under a
common convolution gives
\begin{equation}\label{eq:nonsym-tv-hellinger}
    {\rm TV}(P_1,P_2)
\le {\rm TV}(\widetilde P_1,\widetilde P_2)
\le {\bf H}(\widetilde P_1,\widetilde P_2).
\end{equation}

It therefore suffices to bound the Hellinger distance between
$\widetilde P_1$ and $\widetilde P_2$.

For each $j\in[d]$, define
$$
h_j^2
:=
\int_{\bbR}
\left(
\sqrt{f_0(x)}
-
\sqrt{\frac{f_0((x+\pi_j\gamma_j)/\sigma_j)}{\sigma_j}}
\right)^2\diff x.
$$
The product identity for the Hellinger distance yields
\begin{equation}\label{eq:nonsym-helliger-product-sum}
    {\bf H}^2(\widetilde P_1,\widetilde P_2)
=
2\left\{1-\prod_{j=1}^d\left(1-\frac{h_j^2}{2}\right)\right\}
\le \sum_{j=1}^d h_j^2.
\end{equation}

The two densities in $h_j$ differ in both scale and location, whereas
\Cref{lem: hellinger control} controls a change along a specified parametric
family. We therefore insert the density obtained by changing only the scale
of $f_0$. Specifically, let
$$
p_t(x)=\frac{1}{t}f_0\left(\frac{x}{t}\right),
\qquad t>0,
$$
be the scale family generated by $f_0$. It connects the baseline density
$p_1=f_0$ to the intermediate density $p_{\sigma_j}$. For the remaining
location change, keep the scale fixed at $\sigma_j$ and define
$$
q_{j,u}(x)=\frac{1}{\sigma_j}
f_0\left(\frac{x+u}{\sigma_j}\right),
\qquad u\in\bbR.
$$
This location family starts at $q_{j,0}=p_{\sigma_j}$ and ends at the target
density $q_{j,\pi_j\gamma_j}$ in the definition of $h_j$. Thus, the path from
$f_0$ to $q_{j,\pi_j\gamma_j}$ is decomposed into a scale-only change from
$p_1$ to $p_{\sigma_j}$, followed by a location-only change from $q_{j,0}$ to
$q_{j,\pi_j\gamma_j}$. The squared triangle inequality therefore bounds
$h_j^2/2$ by the sum of the squared Hellinger distances along these two
one-parameter families.

Write $g_0=\sqrt{f_0}$. Since $I_{f_0}<\infty$, $g_0$ has a
square-integrable weak derivative and
$\int_{\bbR}(g_0'(x))^2\,\diff x
= I_{f_0}/4<\infty$.
Because $g_0$ is supported on $[-3/2,3/2]$, both $g_0'$ and
$g_0/2+xg_0'$ are square integrable. The derivatives of the square-root
densities along the location and scale parameters are rescaled versions of
these two functions, respectively. Hence both families are regular.
The Fisher information of the scale family $\{p_t:t>0\}$ and that of the
location family $\{q_{j,u}:u\in\bbR\}$ are, respectively,
\begin{equation}\label{eq:nonsym-information-transform}
I_{\rm sc}(t)=\frac{J_{f_0}}{t^2},
\qquad
I_{\rm loc}(u;\sigma_j)=\frac{I_{f_0}}{\sigma_j^2},
\end{equation}
where
$$
J_{f_0}
:=
\int_{\bbR}
\left(1+x\frac{f_0'(x)}{f_0(x)}\right)^2f_0(x)\,\diff x.
$$
These can be obtained by
differentiating the two
log-densities:
\[
\partial_t\log p_t(x)
=
-\frac{1}{t}
\left\{
1+\frac{x}{t}\frac{f_0'(x/t)}{f_0(x/t)}
\right\},
\quad \text{ and } \quad
\partial_u\log q_{j,u}(x)
=
\frac{1}{\sigma_j}
\frac{f_0'((x+u)/\sigma_j)}{f_0((x+u)/\sigma_j)}.
\]
All ratios involving \(f_0'(x)/f_0(x)\) are interpreted as zero on
\(\{x:f_0(x)=0\}\). Since \(f_0\) is supported on \([-3/2,3/2]\),
\[
J_{f_0}
\le
2+2\int_{\bbR}x^2\frac{(f_0'(x))^2}{f_0(x)}\,\diff x
\le
2+\frac{9}{2}I_{f_0}<\infty.
\]Set \(K_{f_0}:=I_{f_0}\vee J_{f_0}\). By
\Cref{eq:lower-exp-gamma-bound}, \(\sigma_j^2\ge3/4\), and hence
\(\sigma_j\in[\sqrt{3}/2,1]\).
Define \(\mathcal I_j=[\min\{0,\pi_j\gamma_j\},\max\{0,\pi_j\gamma_j\}]\).
By \eqref{eq:nonsym-information-transform}, we have
\[
\sup_{t\in[\sigma_j,1]}I_{\rm sc}(t)
\le \frac{4}{3}K_{f_0},
\qquad
\sup_{u\in\mathcal I_j}I_{\rm loc}(u;\sigma_j)
\le \frac{4}{3}K_{f_0}.
\]
Applying \Cref{lem: hellinger control} to the scale path from \(1\) to
\(\sigma_j\) and to the location path over \(\mathcal I_j\), and then using
the squared triangle inequality, gives

\begin{equation}
    \begin{aligned}
\frac{1}{2}h_j^2
&\le
\int_{\bbR}
\left(
\sqrt{f_0(x)}-
\sqrt{\frac{f_0(x/\sigma_j)}{\sigma_j}}
\right)^2\diff x \\
&\quad+
\int_{\bbR}
\left(
\sqrt{\frac{f_0(x/\sigma_j)}{\sigma_j}}-
\sqrt{\frac{f_0((x+\pi_j\gamma_j)/\sigma_j)}{\sigma_j}}
\right)^2\diff x \\
&\le
\frac{(1-\sigma_j)^2}{4}
\sup_{t\in[\sigma_j,1]}I_{\rm sc}(t)
+
\frac{\pi_j^2\gamma_j^2}{4}
\sup_{u\in\mathcal I_j}I_{\rm loc}(u;\sigma_j) \\
&\le
\frac{K_{f_0}}{3}
\left\{(1-\sigma_j)^2+\pi_j^2\gamma_j^2\right\} \\
&\le
\frac{K_{f_0}}{3}
\left\{(1-\sigma_j^2)^2+\pi_j^2\gamma_j^2\right\} \\
&\le
\frac{K_{f_0}}{3}
(\gamma_j^2+\gamma_j^4)\pi_j^2.
\end{aligned}\label{eq:nonsym-hj-bound}
\end{equation}

The last inequality uses
$1-\sigma_j^2=\gamma_j^2\pi_j(1-\pi_j)$.

To finish the proof, we bound the sum of $h_j^2$.
Since $\pi_j=0$ for $j<s^2$, only the indices $j\ge s^2$ contribute. For
these indices, set $a_j:=|\lambda_{\cal H}/\eta_j|$. The reduction at the
beginning of \Cref{subsec:proof exponential lower} gives $a_j\ge\sqrt{C}\ge1$, and
\eqref{eq:exponential lower condition} becomes
$$
\sum_{j\ge s^2}a_j^4\exp(-2a_j^\al)\le\bar C.
$$
Using $\gamma_j^2=c_2^2a_j^2$, $\pi_j^2=c^2\exp(-2a_j^\al)$, and $c_2\le1$,
we obtain
\begin{align*}
\sum_{j=1}^d(\gamma_j^2+\gamma_j^4)\pi_j^2
&=
c^2c_2^2\sum_{j\ge s^2}
\left(a_j^2+c_2^2a_j^4\right)
\exp(-2a_j^\al) \\
&\le
2c^2c_2^2
\sum_{j\ge s^2}a_j^4\exp(-2a_j^\al) \\
&\le 2\bar Cc^2c_2^2,
\end{align*}
where $(a_j^2+c_2^2a_j^4)\le2a_j^4$ follows from $a_j\ge1$ and
$c_2\le1$.
By \eqref{eq:nonsym-tv-hellinger}, \eqref{eq:nonsym-helliger-product-sum}, and \eqref{eq:nonsym-hj-bound}, we have
$$
{\rm TV}^2(P_1,P_2) \le
{\bf H}^2(\widetilde P_1,\widetilde P_2)
\le
\frac{2K_{f_0}}{3}
\sum_{j=1}^d(\gamma_j^2+\gamma_j^4)\pi_j^2
\le
\frac{4K_{f_0}\bar C}{3}c^2c_2^2.
$$
Thus the lemma holds with $\bar C'=2\sqrt{K_{f_0}\bar C/3}$.

\end{proof}

\section{Proof of Examples in \Cref{sec:example}}\label{sec:pf-example}
\subsection{Homogeneous loading vector}
In this case, we have $\eta_j = 1$ for all $j \in [d]$.
Furthermore, \eqref{eq:solution linear adaptive} gives
\[
\beta_*(s)
=
2\log\left(\frac{2\sqrt{d\log(es)}}{s}\right).
\]
Hence $\lambda_*(s)>0$ if and only if
$s^2<4d\log(es)$, and therefore
$s_0\asymp\sqrt{d\log(ed)}$.

\begin{itemize}
  \item \textit{$\Phi_{\rm o}(s;\eta)$:} From \Cref{eq:solution nonadaptive}, we obtain
  \[
  \frac{s}{2} \;=\; \frac{d \exp(-\beta)}{\sqrt{d \exp(-\beta)}}
  \;=\; \sqrt{d}\,\exp(-\beta/2),
  \]
  which implies
  \[
  \beta = 2 \log\!\left(\frac{2\sqrt{d}}{s}\right).
  \]
  \begin{itemize}
    \item If $s \ge 2\sqrt{d}$, then $\lambda_{o} = 0$ and
    \[
    \Phi_{\rm o}(s;\eta) \;\asymp\; d \;\asymp\; s^2 \log^{2/\al}\!\left(1+\frac{d^{\al/2}}{s^\al}\right).
    \]
    \item If $s < 2\sqrt{d}$, then
    \[
    \Phi_{\rm o}(s;\eta) \;\asymp\; d \exp(-\beta) + s^2 \lambda_{o}^2
    \;\asymp\; s^2 \log^{2/\al}\!\left(\frac{2\sqrt{d}}{s}\right)
    \;\asymp\; s^2 \log^{2/\al}\!\left(1+\frac{d^{\al/2}}{s^\al}\right).
    \]
  \end{itemize}
  Therefore,
  \[
  \Phi_{\rm o}(s;\eta) \;\asymp\; s^2 \log^{2/\al}\!\left(1+\frac{d^{\al/2}}{s^\al}\right).
  \]

\item \textit{$\Phi_{\rm adp}(s;\eta)$:}
We first consider $1\le s\le s_0$ and write
$\ell_s=\log(es)$. From \eqref{eq:phi star equation},
\[
\Phi_*(s;\eta)
\asymp
s^2\log^{2/\al}\!\left(
1+\frac{d^{\al/2}\ell_s^{\al/2}}{s^\al}
\right).
\]
Let $x=s/\ell_s$ and $A=(d/\ell_s)^{\al/2}$. Since
$\ell_s\le s$, we have $x\ge1$.
Let $h(u):=u\log(1+A/u)$ for  $u\in(0,\infty)$. It is elementary to verify that $h$ is non-decreasing.
Consequently,
\begin{align*}
\Phi_*(s;\eta)
\asymp
s^2\log^{2/\al}\!\left(
1+\frac{d^{\al/2}\ell_s^{\al/2}}{s^\al}
\right)
&=
\ell_s^2
\left[
x^\al\log\left(1+\frac{A}{x^\al}\right)
\right]^{2/\al}\\
&=
\ell_s^2
\left[h(x^\al)
\right]^{2/\al}\\
&\ge \ell_s^2
\left[h(1)
\right]^{2/\al} \\
& = \ell_s^2\log^{2/\al}(1+A)\\
&\asymp
\Phi_*(1;\eta)\ell_s^2,
\end{align*}
where the last comparison
follows from
$\ell_s\le\log(ed)$ and
$\Phi_*(1;\eta)=\Phi_{\rm o}(1;\eta)\asymp\log^{2/\al}(1+d^{\al/2})$.
When $0<\al<2$, the preceding bound shows that
\[
\Phi_*(s;\eta)
\gtrsim
\Phi_*(1;\eta)\log^2(es),
\]
so the second term in \eqref{eq:threshold adaptive} is
dominated by $\Phi_*(s;\eta)$. When $\al\ge2$, we have
$\Phi_{\rm adp}(s;\eta)=\Phi_*(s;\eta)$ by definition.
This shows that for $s\le s_0$,
\[
\Phi_{\rm adp}(s;\eta)\asymp\Phi_*(s;\eta).
\]

For $s>s_0$, the definition gives
$\Phi_{\rm adp}(s;\eta)=\Phi_{\rm adp}(s_0;\eta)$.
Recall that $s_0\asymp\sqrt{d\log(ed)}$ in the homogeneous case.
Hence, for all $s_0<s\le d$, we have
\[
\Phi_{\rm adp}(s;\eta)
=\Phi_{\rm adp}(s_0;\eta)
\asymp d\log(ed).
\]
Furthermore, for these $s$, let $u_s=
\frac{d^{\al/2}\log^{\al/2}(es)}{s^\al}$.
Since $s_0\asymp\sqrt{d\log(ed)}$, we have $u_s\lesssim1$, $\log(1+u_s) \asymp u_s$, and
$\log(es)\asymp\log(ed)$. Hence
\[
s^2\log^{2/\al}(1+u_s)
\asymp s^2u_s^{2/\al}
=d\log(es)
\asymp d\log(ed).
\]

Combining both cases, we conclude that for every $\al>0$,
\[
\Phi_{\rm adp}(s;\eta)\asymp
s^2\log^{2/\al}\!\left(
1+\frac{d^{\al/2}\log^{\al/2}(es)}{s^\al}
\right).
\]

  \item \textit{$\Phi_{\rm ns}(s;\eta)$:} From \Cref{thm: exponential upper}, the upper bound is
  \[
  \sum_{j\le j_3(s)} \eta_j^2 \;=\;  j_3(s) = \big\lceil s^2 \log^{2/\al}(ed/s) \big\rceil \wedge d.
  \]
  \begin{itemize}
    \item If $s \lesssim {\sqrt{d}}/{\log^{2/\al}(ed)}$, then by \eqref{eq:solution exponential 0},
    \[
    \lambda_{\cal H} \;=\; \log^{1/\al}\!\left(\frac{d-s^2+1}{s}\right)
    \;\asymp\; \log^{1/\al} d.
    \]
    In this regime, condition \eqref{eq:exponential lower condition} holds, and thus the asymmetric lower bound is
    \(
    s^2 \log^{2/\al}(ed),
    \)
    as given in \eqref{eq:exponential lower asymmetric}.
    \item If $s \gtrsim \sqrt{d}$, then the lower bound reduces to
    \(
    d,
    \)
    from \eqref{eq:exponential lower symmetric}.
  \end{itemize}
\end{itemize}
\subsection{Two-phase loading vector}
In this case, we have $\eta_j = d^{\gamma_\lambda}$ for $1 \le j \le d^{\gamma_d}$ and $\eta_j = 1$ for $d^{\gamma_d} < j \le d$.
Let $\eta^{(1)},\eta^{(2)}$ denote the two subvectors of $\eta$, i.e.,
\[
\eta^{(1)} \in \mathbb{R}^{d^{\gamma_d}}, \quad \eta^{(1)}_j = d^{\gamma_\lambda}, \ \forall j \le d^{\gamma_d},
\]
and
\[
\eta^{(2)} \in \mathbb{R}^{d-d^{\gamma_d}}, \quad \eta^{(2)}_j = 1, \ \forall j \le d-d^{\gamma_d}.
\]

Following the analysis in the homogeneous case, we obtain
\[
\Phi_{\rm o}(s;\eta^{(1)})\;\asymp\;
d^{2\gamma_\lambda}s^2\log^{2/\al}\!\left(1+\frac{d^{\gamma_d\al/2}}{s^\al}\right)
\]
and
\[
\Phi_{\rm o}(s;\eta^{(2)})\;\asymp\;
s^2 \log^{2/\al}\!\left(1+\frac{(d-d^{\gamma_d})^{\al/2}}{s^\al}\right).
\]
Write $\Theta^{(1)}=\{\theta\in \Theta_{s}: \theta_j=0, \forall j>d^{\gamma_d} \}$ and $\Theta^{(2)}=\{\theta\in \Theta_{s}: \theta_j=0, \forall j\le d^{\gamma_d} \}$.
We then have $\Theta^{(1)} \cup \Theta^{(2)}\subset \Theta_{s}$.
It is straightforward to see from the definition of minimax rate that $\Phi_{\rm o}(s;\eta)\;\gtrsim\; \Phi_{\rm o}(s;\eta^{(i)})$ for $i=1,2$.

We also have $\Theta_{s}\subset \Theta^{(1)} +  \Theta^{(2)}$:
For any $\theta$, we can write $\theta=\theta^{(1)}+\theta^{(2)}$ so that $\theta^{(i)}\in \Theta^{(i)}$. It follows that  $L(\theta)=\langle\eta, \theta\rangle=\langle \eta^{(1)}, \theta^{(1)}\rangle + \langle \eta^{(2)}, \theta^{(2)}\rangle$.
Therefore, for any minimax estimators for $\langle \eta^{(1)}, \theta^{(1)}\rangle$ and $\langle \eta^{(2)}, \theta^{(2)}\rangle$, the sum of their worst-case risks (multiplied by 2) is an upper bound on the minimax risk for $L(\theta)$.
This implies that $\Phi_{\rm o}(s;\eta)\lesssim \Phi_{\rm o}(s;\eta^{(1)})+\Phi_{\rm o}(s;\eta^{(2)})$.

To sum up, the minimax rate for estimating $\langle\eta, \theta\rangle$ satisfies
\[
\Phi_{\rm o}(s;\eta)\;\asymp\; \Phi_{\rm o}(s;\eta^{(1)}) + \Phi_{\rm o}(s;\eta^{(2)}).
\]
Therefore,
\begin{equation}\label{eq:example two phase overall minimax}
\Phi_{\rm o}(s;\eta)\;\asymp\; d^{2\gamma_\lambda}s^2\log^{2/\al}\!\left(1+\frac{d^{\gamma_d\al/2}}{s^\al}\right)
+ s^2 \log^{2/\al}\!\left(1+\frac{d^{\al/2}}{s^\al}\right), \quad s\in [1,d].
\end{equation}

We further simplify \Cref{eq:example two phase overall minimax} as follows.
\begin{enumerate}
    \item When $s\lesssim d^{\gamma_d/2}$, both the logarithmic factors are equivalent to $\log (d)$ so the first term dominates.
    \item When $ d^{\gamma_d/2}\ll s \ll  d^{\gamma_\lambda+\gamma_d/2}$ (which is $\ll \sqrt{d}$ by assumption), the second logarithmic factor remains the same as $\log(d)$ but the first factor becomes $\asymp \left(\frac{d^{\gamma_d/2}}{s}\right)^{\al}$.
    In this case, the two terms are $d^{2\gamma_\lambda+\gamma_d} $ and $s^2 \log^{2/\al}(d)$. This suggests that the first term dominates if and only if $s^2\lesssim \frac{d^{2\gamma_\lambda+\gamma_d} }{\log^{2/\al}(d)}$.
\end{enumerate}

\medskip
We next derive the adaptive rate.

We first determine the order of $s_0$.
At $\beta=0$, the
left-hand side of \eqref{eq:solution linear adaptive} is
\[
\frac{
d^{\gamma_\lambda+\gamma_d}
+d-d^{\gamma_d}
}{
\sqrt{
d^{2\gamma_\lambda+\gamma_d}
+d-d^{\gamma_d}
}
}
\asymp \sqrt{d},
\]
where we have used $2\gamma_\lambda+\gamma_d<1$.
As in the homogeneous case, the cutoff $\lambda_*(s)>0$
is therefore characterized by
\[
\frac{s}{\sqrt{\log(es)}}\lesssim\sqrt{d}.
\]
Since $s\mapsto s/\sqrt{\log(es)}$ is increasing on $[1,\infty)$,
the boundary value satisfies
$s_0^2\asymp d\log(es_0)$.

Moreover, $\log(es_0)\asymp\log(ed)$, and hence
\[
s_0\asymp\sqrt{d\log(ed)}.
\]

In the following, write $\ell_s=\log(es)$ and we consider $s\le s_0$ and $s>s_0$ separately.

\textbf{Part 1.}
We first consider
$1\le s\le s_0$.  By
\eqref{eq:phi star equation} and
\eqref{eq:example two phase overall minimax},
\begin{align*}
\Phi_*(s;\eta)
\asymp
\ell_s\,
\Phi_{\rm o}\!\left(
\frac{s}{\sqrt{\ell_s}};\eta
\right)
\asymp
d^{2\gamma_\lambda}s^2
\log^{2/\al}\!\left(
1+
\frac{
d^{\gamma_d\al/2}\ell_s^{\al/2}
}{
s^\al
}
\right)
+
s^2
\log^{2/\al}\!\left(
1+
\frac{
d^{\al/2}\ell_s^{\al/2}
}{
s^\al
}
\right).
\end{align*}

When $\al\ge2$, we have
$\Phi_{\rm adp}(s;\eta)=\Phi_*(s;\eta)$ by definition.

Next, suppose that $0<\al<2$.
Applying the monotonicity
argument used for the homogeneous loading vector separately
to $\eta^{(1)}$ and $\eta^{(2)}$ yields
\begin{align*}
&d^{2\gamma_\lambda}s^2
\log^{2/\al}\!\left(
1+
\frac{
d^{\gamma_d\al/2}\ell_s^{\al/2}
}{
s^\al
}
\right)
\gtrsim
\Phi_{\rm o}(1;\eta^{(1)})\ell_s^2,\\
&s^2
\log^{2/\al}\!\left(
1+
\frac{
d^{\al/2}\ell_s^{\al/2}
}{
s^\al
}
\right)
\gtrsim
\Phi_{\rm o}(1;\eta^{(2)})\ell_s^2.
\end{align*}
Here we have used $d-d^{\gamma_d}\asymp d$ and
\[
\log\left(1+
\left(\frac{d^{\gamma_d}}{\ell_s}\right)^{\al/2}
\right)
\asymp
\log\left(1+d^{\gamma_d\al/2}\right),
\qquad
\log\left(1+
\left(\frac{d}{\ell_s}\right)^{\al/2}
\right)
\asymp
\log\left(1+d^{\al/2}\right),
\]
which follow from $\ell_s\le\log(ed)$ and
$\gamma_d>0$.

Hence, we have
\begin{align*}
\Phi_*(s;\eta)
&\gtrsim
\left\{
\Phi_{\rm o}(1;\eta^{(1)})
+
\Phi_{\rm o}(1;\eta^{(2)})
\right\}\ell_s^2\\
&\asymp
\Phi_{\rm o}(1;\eta)\ell_s^2
=
\Phi_*(1;\eta)\log^2(es),
\end{align*}
where the second comparison follows from
\eqref{eq:example two phase overall minimax} with $s=1$.
Thus, when $0<\al<2$, the second term in
\eqref{eq:threshold adaptive} is dominated by
$\Phi_*(s;\eta)$, and thus $\Phi_{\rm adp}(s;\eta)\asymp \Phi_*(s;\eta)$.

Consequently, for $1\le s\le s_0$ and any $\al>0$, we have
\begin{align*}
\Phi_{\rm adp}(s;\eta)
&\asymp
d^{2\gamma_\lambda}s^2
\log^{2/\al}\!\left(
1+
\frac{
d^{\gamma_d\al/2}\log^{\al/2}(es)
}{
s^\al
}
\right)\\
&\quad+
s^2
\log^{2/\al}\!\left(
1+
\frac{
d^{\al/2}\log^{\al/2}(es)
}{
s^\al
}
\right).
\end{align*}

\textbf{Part 2.}
We now consider $s>s_0$.
The definition gives
$\Phi_{\rm adp}(s;\eta)=\Phi_{\rm adp}(s_0;\eta)$.
Evaluating the expression obtained in Part~1 at $s=s_0$
and using $s_0^2\asymp d\log(ed)$, we obtain
\begin{align*}
\Phi_{\rm adp}(s_0;\eta)
&\asymp
d^{2\gamma_\lambda+\gamma_d}\log(ed)
+d\log(ed)\\
&\asymp d\log(ed),
\end{align*}
where the last comparison follows from
$2\gamma_\lambda+\gamma_d<1$. Therefore,
\[
\Phi_{\rm adp}(s;\eta)\asymp d\log(ed),
\qquad s_0<s\le d.
\]

For these $s$, let
$u_{1,s}=d^{\gamma_d\al/2}\log^{\al/2}(es)/s^\al$
and $u_{2,s}=d^{\al/2}\log^{\al/2}(es)/s^\al$.
Since $s_0\asymp\sqrt{d\log(ed)}$, we have
$u_{1,s}\le u_{2,s}\lesssim1$,
$\log(1+u_{i,s})\asymp u_{i,s}$ for $i=1,2$, and
$\log(es)\asymp\log(ed)$. Applying the same calculation
as in the homogeneous case, we can obtain

\begin{align*}
&d^{2\gamma_\lambda}s^2
\log^{2/\al}\!\left(
1+
\frac{
d^{\gamma_d\al/2}\log^{\al/2}(es)
}{
s^\al
}
\right)
\asymp
d^{2\gamma_\lambda+\gamma_d}\log(ed),\\
&s^2
\log^{2/\al}\!\left(
1+
\frac{
d^{\al/2}\log^{\al/2}(es)
}{
s^\al
}
\right)
\asymp
d\log(ed).
\end{align*}
The second term dominates
because $2\gamma_\lambda+\gamma_d<1$.

\textbf{Synthesis.}
Combining the cases
$s\le s_0$ and $s>s_0$, we obtain
\begin{align*}
\Phi_{\rm adp}(s;\eta)
&\asymp
d^{2\gamma_\lambda}s^2
\log^{2/\al}\!\left(
1+
\frac{
d^{\gamma_d\al/2}\log^{\al/2}(es)
}{
s^\al
}
\right)\\
&\quad+
s^2
\log^{2/\al}\!\left(
1+
\frac{
d^{\al/2}\log^{\al/2}(es)
}{
s^\al
}
\right).
\end{align*}

\subsection{Exponentially decaying loading vector}
Some of the calculations for this example are borrowed from Example E.4 in \cite{chhor2024sparse}.

Recall that $j_0 = \min\{\, j \mid \eta_j < 1/2 \,\}$.
We first show that for any $\gamma\ge 1$, the following relation holds for $\beta\ge 0$:
\begin{equation}\label{eq:example exp decay relation}
    \sum_{j \ge j_0} \eta_j^\gamma \exp \left(-\beta / \eta_j^\al\right) \lesssim \sum_{j < j_0} \eta_j^\gamma \exp \left(-\beta / \eta_j^\al\right).
\end{equation}
We write $j_*=j_0-1$.
If $j_*=d$, the relation is obvious.

If $j_*<d$, then $\varphi(j_*)>\log(2)$. Furthermore, we will make use of the convexity of $\varphi$ to derive
\begin{equation}\label{eq:example-exp-convex-j_*}
    \varphi(j_*+\ell)\ge\varphi(j_*)+\ell g_*,
\qquad \ell\ge0,
\end{equation}
where we pick $g_*\in\partial\varphi(j_*)$.

For any $\beta\ge 0$, we have
$$
\begin{aligned}
\sum_{j \ge j_0} \eta_j^\gamma \exp \left(-\beta / \eta_j^\al\right) = & \sum_{j \ge j_0}  \exp \left(-\beta / \eta_j^\al - \gamma \varphi(j-1)\right) \\
(\because \varphi \text{ is non-decreasing and } \beta\ge 0) \qquad \le  &  \sum_{j \ge j_*}  \exp \left(-\beta / \eta_{j_0}^\al - \gamma \varphi(j) \right)\\
(\because  \text{inequality in \eqref{eq:example-exp-convex-j_*}})\qquad
\le  &  \sum_{\ell \ge 0}  \exp \left(-\beta / \eta_{j_0}^\al - \gamma \varphi(j_*) - \gamma \ell g_*\right) \\
& \le \frac{\exp \left(-\beta / \eta_{j_0}^\al - \gamma \varphi(j_*) \right)}{1 - \exp \left(-\gamma  g_* ) \right)}\\
(\because \text{definition of } j_*) \qquad
&\le \frac{\exp \left(-\beta / \eta_{j_0}^\al  \right)}{1 - \exp \left(-\gamma  g_* ) \right)} \times 2^{-\gamma}.
\end{aligned}
$$

On the other hand, since $\varphi$ is non-decreasing and convex with
$\varphi(0)=0$, for every $0\le j\le j_*$ we have
\[
0\le\varphi(j)
\le
\frac{j}{j_*}\varphi(j_*)
\le
\varphi(j_*).
\]

Therefore,
$$
\begin{aligned}
\sum_{1\le j <  j_0} \eta_j^\gamma \exp \left(-\beta / \eta_j^\al\right)= & \sum_{0\le j \le j_*-1}  \exp \left(-\beta / \eta_{j+1}^\al - \gamma \varphi(j)\right) \\
(\because  \beta\ge 0) \qquad  \ge  &  \sum_{j \le j_*-1}  \exp \left(-\beta / \eta_{j_0}^\al - \gamma j \varphi(j_*)/j_*\right) \\
= & \exp \left(-\beta / \eta_{j_0}^\al\right) \frac{1-\exp\left(-\gamma  j_*\varphi(j_*)/j_*\right)}{1-\exp\left(-\gamma  \varphi(j_*)/j_*\right)}  \\
(\because \text{definition of } j_*) \qquad
\ge & \exp \left(-\beta / \eta_{j_0}^\al\right) \frac{1-2^{-\gamma}}{1-\exp\left(-\gamma  \varphi(j_*)/j_*\right)} \\
(\because \varphi(j_*)/j_*\le g_*) \qquad
\ge & \exp \left(-\beta / \eta_{j_0}^\al\right) \frac{1-2^{-\gamma}}{1-\exp\left(-\gamma  g_*\right)}.
\end{aligned}
$$

Therefore, for any $\beta\ge 0$, we have
\begin{equation}\label{eq:example exp decay equiv ratio}
    \frac{\sum_{j \ge j_0} \eta_j^\gamma \exp \left(-\beta / \eta_j^\al\right)}{\sum_{1\le j <  j_0} \eta_j^\gamma \exp \left(-\beta / \eta_j^\al\right)}\le 2^\gamma-1.
\end{equation}

We take $\gamma =1$ and $\gamma=2$ to simplify the numerator and denominator in \Cref{eq:solution nonadaptive} into $$
\sum_{j<j_0} |\eta_j| \exp(-\beta / |\eta_j|^\al), \text{ and } \sqrt{\sum_{j<j_0} \eta_j^2 \exp(-\beta / |\eta_j|^\al)}
$$
up to a constant. Furthermore, for any $j<j_0$, we have $1/2<\eta_j\le 1$. Therefore, the solution $\beta=\beta(s;\eta)$ to \Cref{eq:solution nonadaptive} satisfies that
\begin{equation}\label{eq:two phase}
    \frac{s}{2} \asymp \frac{\sum_{j<j_0} |\eta_j| \exp(-\beta / |\eta_j|^\al)}{\sqrt{\sum_{j<j_0} \eta_j^2 \exp(-\beta / |\eta_j|^\al)}}  \asymp \frac{\sum_{j<j_0}  \exp(-\beta / |\eta_j|^\al)}{\sqrt{\sum_{j<j_0}  \exp(-\beta / |\eta_j|^\al)}},
\end{equation}
which is lower bounded by $\sqrt{j_0 \exp(-2^\al \beta) }$ and upper bounded $\sqrt{j_0 \exp(-\beta) }$ up to absolute constants.

We discuss as follows:
\begin{enumerate}
    \item Suppose $\beta> 1$, by \eqref{eq:two phase}, we have $$\sum_{j<j_0}\exp(-\beta/|\eta_j|^\al)\asymp s^2.$$
    Since $\eta_j\in (1/2,1)$ for any $j<j_0$, we have $\lambda_o=\beta^{1/\al}\asymp \log^{1/\al}(1+j_0^{\al/2}/s^\al).$
    By \Cref{eq:example exp decay equiv ratio} with $\gamma=2$, we have
    $$\nu^2(s)\asymp \sum_{j< j_0}\eta_j^2 \exp(-\beta/|\eta_j|^\al)\asymp \sum_{j<j_0}\exp(-\beta/|\eta_j|^\al)\asymp s^2.$$
    Therefore, we have
    $$\Phi_{\rm o}(s;\eta)\asymp s^2\lambda_o^2+\nu^2\asymp s^2\lambda_o^2\asymp s^2\log^{2/\al}(1+j_0^{\al/2}/s^\al).$$

\item Suppose $0\le\beta\le1$.
The bounds in
\eqref{eq:two phase} imply that $s^2\asymp j_0$.
By \Cref{eq:example exp decay equiv ratio} with $\gamma=2$, we have $$\nu^2 \asymp \sum_{j < j_0} \eta_j^2\exp(-\beta/|\eta_j|^\al) \asymp \sum_{j<j_0}\exp(-\beta/|\eta_j|^\al)\asymp j_0.$$
Moreover, $\lambda_o^2=\beta^{2/\al}\lesssim1$.
Therefore,
\[
\Phi_{\rm o}(s;\eta)
\asymp s^2\lambda_o^2+\nu^2(s)
\asymp j_0
\asymp
s^2\log^{2/\al}\!\left(
1+\frac{j_0^{\al/2}}{s^\al}
\right).
\]
\item Suppose $\beta<0$. Then $\beta_+=0$ and hence
$\lambda_o=0$.
As in the case with $0\le \beta \le 1$, we have
\[
\nu^2(s)
=
\sum_{j=1}^d\eta_j^2
\exp\left\{
-\left(\frac{\lambda_o}{\eta_j}\right)^\al
\right\}
=
\sum_{j=1}^d\eta_j^2 \asymp j_0.
\]
Since $\lambda_o=0$, we
conclude that
\[
\Phi_{\rm o}(s;\eta)
=
\left(\lambda_os+\nu(s)\right)^2
=
\nu^2(s)
\asymp j_0.
\]

\end{enumerate}

Combining the three cases above, we obtain, for every real
$u\in[1,d]$,
\begin{equation}\label{eq:exp-decay-effective-dimension}
\Phi_{\rm o}(u;\eta)
\asymp
\Phi_{\rm o}(u\wedge j_0;\mathbf{1}_{j_0}),
\qquad 1\le u\le d.
\end{equation}

The same three cases show that the transition between
$\beta(u)>0$ and $\beta(u)\le0$ occurs at
$u\asymp\sqrt{j_0}$.

We next derive the adaptive rate.

Let $s_{0,\eta}$ and $s_{0,h}$ denote the adaptive cutoff values
associated with $\eta$ and $\mathbf{1}_{j_0}$, respectively.
Since \eqref{eq:solution linear adaptive} is obtained from
\eqref{eq:solution nonadaptive} by replacing $u$ with
$s/\sqrt{\log(es)}$, we have
\[
s_{0,\eta}
\asymp
s_{0,h}
\asymp
\sqrt{j_0\log(ej_0)}.
\]

\textbf{Regime 1.}
If $s\le s_{0,\eta}\wedge s_{0,h}$, then
\eqref{eq:phi star equation} gives
\begin{align*}
\Phi_*(s;\eta)
&\asymp
\log(es)\,
\Phi_{\rm o}\!\left(
\frac{s}{\sqrt{\log(es)}};\eta
\right)\\
&\asymp
\log(es)\,
\Phi_{\rm o}\!\left(
\frac{s}{\sqrt{\log(es)}};\mathbf{1}_{j_0}
\right)
\asymp
\Phi_*(s;\mathbf{1}_{j_0}).
\end{align*}

\textbf{Regime 2.}
If $s\ge s_{0,\eta}\vee s_{0,h}$, both quantities have reached
their constant values. Since the case $\beta<0$ gives
$\sum_{j=1}^d\eta_j^2\asymp j_0$, we have
\begin{equation}\label{eq:exp-decay-phi-star-saturation}
\Phi_*(s;\eta)
\asymp
\Phi_*(s \wedge j_0 ;\mathbf{1}_{j_0})
\asymp
j_0\log(ej_0),
\qquad
s\ge s_{0,\eta}\vee s_{0,h}.
\end{equation}

\textbf{Regime 3.}
If $s$ lies between the two cutoff values, then
\begin{equation}\label{eq:exp-decay-intermediate-scales}
s\asymp\sqrt{j_0\log(ej_0)},
\qquad
\log(es)\asymp\log(ej_0),
\qquad
\frac{s}{\sqrt{\log(es)}}\asymp\sqrt{j_0}.
\end{equation}
Applying \eqref{eq:exp-decay-effective-dimension} with
$u=s/\sqrt{\log(es)}$ and then using the homogeneous rate formula
in \eqref{eq:homogeneous case}, we obtain
\begin{align*}
\Phi_{\rm o}\!\left(
\frac{s}{\sqrt{\log(es)}};\eta
\right)
&\asymp
\Phi_{\rm o}\!\left(
\frac{s}{\sqrt{\log(es)}};\mathbf{1}_{j_0}
\right)\\
&\asymp j_0,
\end{align*}
where the last comparison uses
$s/\sqrt{\log(es)}\asymp\sqrt{j_0}$.

If $s_{0,\eta}\le s<s_{0,h}$, then
$\Phi_*(s;\eta)$ has reached its constant value, whereas
$\Phi_*(s;\mathbf{1}_{j_0})$ is still given by
\eqref{eq:phi star equation}. Therefore,
\begin{align*}
\Phi_*(s;\eta)
&=
\Phi_*(s_{0,\eta};\eta)
\asymp j_0\log(ej_0),\\
\Phi_*(s;\mathbf{1}_{j_0})
&\asymp
\log(es)\,
\Phi_{\rm o}(\frac{s}{\sqrt{\log(es)}};\mathbf{1}_{j_0})
\asymp j_0\log(ej_0),
\end{align*}
where the first line uses \eqref{eq:exp-decay-phi-star-saturation} and the second line
 uses \eqref{eq:exp-decay-intermediate-scales}.

If $s_{0,h}\le s<s_{0,\eta}$, then
$\Phi_*(s;\mathbf{1}_{j_0})$ has reached its constant value,
whereas $\Phi_*(s;\eta)$ is still given by
\eqref{eq:phi star equation}. Using again \eqref{eq:exp-decay-phi-star-saturation} and \eqref{eq:exp-decay-intermediate-scales}, we obtain
\begin{align*}
\Phi_*(s;\mathbf{1}_{j_0})
&=
\Phi_*(s_{0,h};\mathbf{1}_{j_0})
\asymp j_0\log(ej_0),\\
\Phi_*(s;\eta)
&\asymp
\log(es)\,
\Phi_{\rm o}(\frac{s}{\sqrt{\log(es)}};\eta)
\asymp j_0\log(ej_0).
\end{align*}

Thus, in either ordering of $s_{0,\eta}$ and $s_{0,h}$,
\[
\Phi_*(s;\eta)
\asymp
\Phi_*(s;\mathbf{1}_{j_0})
\asymp
j_0\log(ej_0).
\]

\textbf{Synthesis.}
Combining the three ranges of $s$, we obtain
\[
\Phi_*(s;\eta)
\asymp
\Phi_*(s\wedge j_0;\mathbf{1}_{j_0}),
\qquad 1\le s\le d.
\]

When $\al\ge2$, the desired adaptive comparison follows directly
from the definition of $\Phi_{\rm adp}$. When $0<\al<2$, we also
have
\[
\Phi_*(1;\eta)
\asymp
\Phi_*(1;\mathbf{1}_{j_0})
\]
and
\[
\log\!\left(e(s\wedge s_{0,\eta})\right)
\asymp
\log\!\left(e(s\wedge s_{0,h})\right).
\]
Hence the two terms in the maximum in
\eqref{eq:threshold adaptive} are equivalent, and therefore
\[
\Phi_{\rm adp}(s;\eta)
\asymp
\Phi_{\rm adp}(s\wedge j_0;\mathbf{1}_{j_0}),
\qquad 1\le s\le d.
\]

\subsection{Constant-sparsity regime}\label{app:constant-sparsity}
\begin{corollary}[Constant-sparsity regime]
\label{cor:constant-sparsity}
Suppose that the sparsity level satisfies
\[
s \le \overline{s}
\]
for some absolute constant \(\overline{s}\ge 1\) independent of \(d\). Then the minimax risk satisfies
\[
\Phi_{\rm o}(s;\eta)\asymp \max_{j\in[d]} \eta_j^2 \log^{2/\al}(1+j).
\]
\end{corollary}

\begin{proof}
Since \(s\) is bounded by an absolute constant, it is enough to establish matching upper and lower bounds up to constants depending only on \(\overline{s}\) and \(\al\).

\textbf{Upper bound.}\\
Consider the thresholding estimator
\[
\hat L=\sum_{j=1}^d \eta_j y_j \bbo\{|y_j|>t_j\},
\qquad
t_j=A_1 \log^{1/\al}(1+j),
\]
where \(A_1>0\) is a sufficiently large constant.

Following the decomposition used in the proof of the upper bound above, for any \(\theta\in \Theta_s\), we have
\begin{align*}
    \bbE_\theta (\hat{L}-L)^2
    &=\bbE_\theta\left[\sum_{j\in S}\eta_j\xi_j-\sum_{j\in S}\eta_jy_j\bbo\{|y_j|\le t_j\}+\sum_{j\in S^c}\eta_j \xi_j \bbo\{|y_j|>t_j\}\right]^2\\
    &\le 3\left[\sum_{j\in S}\bbE_\theta\eta_j^2\xi_j^2+s^2\max_{j\in S}\eta_j^2t_j^2+\sum_{j\in S^c}\bbE_\theta\eta_j^2 \xi_j^2 \bbo\{|\xi_j|>t_j\}\right].
\end{align*}
Using the same argument as before, together with the Cauchy--Schwarz inequality and the tail bounds of the noise, we obtain
\begin{align*}
    \bbE_\theta (\hat{L}-L)^2
    \lesssim \sum_{j\in S}\eta_j^2+s^2\max_{j\in S}\eta_j^2t_j^2+\sum_{j\in S^c}\eta_j^2 \exp\left(-\frac{A_1^\al}{\tau^\al}\log(1+j)\right).
\end{align*}
When \(A_1\) is chosen sufficiently large, the last term satisfies
\[
\sum_{j\in S^c}\eta_j^2 \exp\left(-\frac{A_1^\al}{\tau^\al}\log(1+j)\right)
\lesssim \eta_1^2\sum_{j=1}^d \frac{1}{(1+j)^2}
\lesssim \eta_1^2.
\]
Since \(s\le \overline{s}\), it follows that
\[
\sum_{j\in S}\eta_j^2+s^2\max_{j\in S}\eta_j^2t_j^2
\lesssim \max_{j\in S}\eta_j^2 \log^{2/\al}(1+j),
\]
and therefore
\[
\bbE_\theta (\hat{L}-L)^2
\lesssim
\max_{j\in[d]}\eta_j^2 \log^{2/\al}(1+j).
\]
This proves the desired upper bound.

\textbf{Lower bound.}\\
Let
\[
j_* \in \arg\max_{j\in[d]} \eta_j^2 \log^{2/\al}(1+j).
\]
We distinguish two cases.

\medskip
\noindent
\emph{Case 1: \(j_*\) is sufficiently large.}

Let $P_0$ denote the population under $\theta= \mathbf{0}_d$ and let $p_0$ denote the density of $P_0$.
Let $\mu_0$ be the point mass at $\mathbf{0}_d$. Under $\mu_0$,
\[
\mu_0\!\left(
L(\theta)\le0,\ \theta\in\Theta_s
\right)=1.
\]

Consider the collection of one-sparse vectors
\[
\theta^{(k)}
=
a_0\,\mathrm{sgn}(\eta_k)
\log^{1/\al}(1+j_*)e_k,
\qquad 1\le k\le j_*,
\]

where \(a_0>0\) is a sufficiently small constant to be determined.
We have
\begin{equation}\label{eq:constant-sparsity-L-separation}
    L(\theta^{(k)})-L(\mathbf{0}_d)
=
a_0|\eta_k|\log^{1/\al}(1+j_*)
\ge
a_0|\eta_{j_*}|\log^{1/\al}(1+j_*),
\end{equation}
where the inequality follows from
$|\eta_k|\ge|\eta_{j_*}|$ for $k\le j_*$.
Therefore, all alternatives are separated from the null in the same direction.

For $1\le k\le j_*$, let $P_k$ denote the distribution of
the observations under $\theta^{(k)}$, and let $p_k$
denote the densities of $P_k$. By the
definition of $\theta^{(k)}$,
\[
p_k(x)
=
f_\al^{(0)}\!\left(
x_k-a_0\mathrm{sgn}(\eta_k)
\log^{1/\al}(1+j_*)
\right)
\prod_{\ell\neq k}f_\al^{(0)}(x_\ell),
\qquad x\in\mathbb R^d.
\]

Let \(\mu\) be the uniform prior on
\(\{\theta^{(1)},\dots,\theta^{(j_*)}\}\), and let \(P_\mu\) denote the corresponding uniform mixture distribution with density $p_\mu(x)$.
Then, the mixture density
satisfies
\[
\frac{p_\mu(x)}{p_0(x)}
=
\frac{1}{j_*}\sum_{k=1}^{j_*}r_k(x_k),
\]
where
\[
r_k(t)
=
\frac{
f_\al^{(0)}\!\left(
t-a_0\mathrm{sgn}(\eta_k)
\log^{1/\al}(1+j_*)
\right)
}{
f_\al^{(0)}(t)
}.
\]
Under $P_0$, the random variables $r_k(X_k)$ are independent
and satisfy $\mathbb E_0r_k(X_k)=1$. Therefore,
\begin{align*}
1+\chi^2(P_\mu\|P_0)
&=
\mathbb E_0\left[
\frac{1}{j_*}\sum_{k=1}^{j_*}r_k(X_k)
\right]^2\\
&=
1+\frac{1}{j_*^2}
\sum_{k=1}^{j_*}
\chi^2(f_{\al,k}^*\|f_\al^{(0)}),
\end{align*}
where
\[
f_{\al,k}^*(t)
=
f_\al^{(0)}\!\left(
t-a_0\mathrm{sgn}(\eta_k)
\log^{1/\al}(1+j_*)
\right).
\]

Let
\[
f_\al^*(t)
=
f_\al^{(0)}\!\left(
t-a_0\log^{1/\al}(1+j_*)
\right).
\]
Since $f_\al^{(0)}$ is symmetric, positive and negative location
shifts of the same magnitude have the same chi-square divergence.
Thus,
\[
\chi^2(f_{\al,k}^*\|f_\al^{(0)})
=
\chi^2(f_\al^*\|f_\al^{(0)}),
\qquad 1\le k\le j_*.
\]
Consequently, using $1+x\le e^x$ and
\Cref{lem: chi square},
\begin{align*}
1+\chi^2(P_\mu\|P_0)
&=
1+\frac{1}{j_*}
\chi^2(f_\al^*\|f_\al^{(0)})\\
&\le
\exp\left\{
\frac{1}{j_*}
\chi^2(f_\al^*\|f_\al^{(0)})
\right\}\\
&\le
\exp\left[
\frac{1}{j_*}
\left\{
C_{\al,1}
\exp\left(
a_0^\al
\frac{\log(1+j_*)}{C_{\al,2}^\al}
\right)-1
\right\}
\right].
\end{align*}

Hence, once \(j_*\) is sufficiently large, choosing \(a_0\) small enough ensures that
\[
\chi^2(P_\mu,P_0)\le 1/2.
\]

By \eqref{eq:constant-sparsity-L-separation}, we have
\[
\mu\!\left(
L(\theta)\ge
a_0|\eta_{j_*}|\log^{1/\al}(1+j_*),
\ \theta\in\Theta_s
\right)=1.
\]

Moreover, $\mathrm{TV}(P_\mu,P_0)
\le
\frac12\sqrt{\chi^2(P_\mu\|P_0)}
\le
\frac{1}{2\sqrt{2}}$.
Applying \Cref{lem: general tool} with
\[
P_1=P_0,\qquad
P_2=P_\mu,\qquad
c=0,\qquad
t=\frac{a_0}{2}|\eta_{j_*}|
\log^{1/\al}(1+j_*),
\]
we conclude that
\[
\inf_{\hat T}
\sup_{\theta\in\Theta_s}
\sup_{{\cal P}_\xi\in{\cal G}_{\al,\tau}^{\otimes}}
\mathbb E_{\theta,{\cal P}_\xi}
\left(\hat T-L(\theta)\right)^2
\gtrsim
\eta_{j_*}^2\log^{2/\al}(1+j_*).
\]

\medskip
\noindent
\emph{Case 2: \(j_*\) is bounded.}
If \(j_*\le J\) for some absolute constant \(J\), then
\[
\max_{j\in[d]} \eta_j^2 \log^{2/\al}(1+j)
\asymp \eta_1^2.
\]
In this regime, the desired lower bound follows directly from the argument in \Cref{sec: add proof lower}.

Combining the two cases with the upper bound completes the proof.
\end{proof}

\section{Comparison with Chhor, Mukherjee, and Sen (2024)}
\label{app:comparison-chhor}

This appendix clarifies the relation between our lower-bound argument and that of \cite{chhor2024sparse}, while emphasizing the new ingredients required in our setting. In the Gaussian case, our defining equation \eqref{eq:solution nonadaptive} formally resembles Equation~(5) of \cite{chhor2024sparse}. However, the underlying problems are fundamentally different. \cite{chhor2024sparse} study testing of a sign-invariant alternative, whereas our goal is to estimate the signed heterogeneous linear functional
\[
L(\theta)=\sum_{j=1}^d \eta_j\theta_j.
\]
This difference changes the calibration target, and therefore leads to a different prior construction and a different \(\chi^2\)-divergence calculation.

Moreover, our analysis extends beyond the Gaussian case through the new \(\chi^2\)-control for generalized Gaussian noise established in Lemma~\ref{lem: chi square}. These are the main points at which our argument departs from and extends \cite{chhor2024sparse}.

\subsection{The lower-bound construction in \cite{chhor2024sparse}}

\cite{chhor2024sparse} consider sparse signal detection in a heteroscedastic Gaussian sequence model for $d$ observations:
\begin{equation}\label{eq:chhor model}
Z_j\sim N(\mu_j,\sigma_j^2).
\end{equation}

Let $\mu=(\mu_1,\ldots, \mu_d)^\top$.
They consider testing
$$
H_0: \mu=0 \quad\text{ against } \quad H_1: \|\mu\|_t \ge \epsilon, \|\mu\|_0 \le s,
$$
where $\|\cdot\|_t$ is the $L^t$ norm and $s$ is the sparsity level.
For bounded $s$ (or so-called the ``case of constant sparsity''), the analysis is based on the construction of two 1-sparse priors; this part is  straightforward and completed in their Appendix A.1.

For $s$ unbounded, their lower bound is more interesting and is based on the random-sparsity prior constructed as follows:
\begin{equation}\label{eq: prior in  chhor}
    \mu_j=b_j\omega_j\gamma_j,\qquad j=1,\dots,d,
\end{equation}
where the mutually independent random variables satisfy
\[
b_j\sim \mathrm{Ber}(\pi_j),\qquad \omega_j\sim \mathrm{Rad}(1/2),
\]
while $\pi_j$ and $\gamma_j$ are nonnegative constants to be determined.
Thus, conditional on being active (i.e., $b_j=1$), the \(j\)-th coordinate takes the values \(+\gamma_j\) and \(-\gamma_j\) with equal probability.

The parameters \(\{\pi_j,\gamma_j: j=1,\dots,d \}\) are calibrated through the optimization problem
\begin{equation}\label{eq:chhor-optimization-app}
\begin{aligned}
\max_{\gamma_{1:d},\pi_{1:d}}\quad & \sum_{j=1}^d \gamma_j^t\pi_j\\
\text{s.t.}\quad
& \sum_{j=1}^d \pi_j = s/2,\\
& \pi_j\in[0,1],\quad \gamma_j\ge 0,\qquad j\in[d],\\
& \sum_{j=1}^d \pi_j^2 \sinh^2\!\left(\frac{\gamma_j^2}{2\sigma_j^2}\right)\le c,
\end{aligned}
\end{equation}
where \(c>0\) is a sufficiently small constant.
In view of \eqref{eq: prior in  chhor}, the objective of this optimization problem satisfies that
\[
\sum_{j=1}^d \gamma_j^t\pi_j =\mathbb{E}\|\mu\|_t^t,
\]
the first constraint controls the expected sparsity of the random $\mu$, and the last constraint controls the \(\chi^2\)-distance between the induced mixture and the null distributions of data.

Although \eqref{eq:chhor-optimization-app} is formulated for all $t \in[1, \infty]$, the lower-bound calibration in \cite{chhor2024sparse} takes different forms in the regimes $t \in[1,2]$ and $t \ge 2$. Since our defining equation \eqref{eq:solution nonadaptive} , in the Gaussian case formally resembles Equation~(5) of \cite{chhor2024sparse}, the most relevant direct comparison is with their $t \ge 2$ random-sparsity construction. Their treatment of $t \in[1,2]$, developed in their Section 3 and based on the more involved calibration around their equations (23)-(28), is structurally different and is not the part directly mirrored by our equation.

\subsection{Our lower bound separates a different functional}

We consider estimation of the signed linear functional
$
L(\theta)=\sum_{j=1}^d \eta_j\theta_j
$
under the homoscedastic Gaussian sequence model
\[
Y_j\sim N(\theta_j,\sigma^2),\qquad j=1,\dots,d.
\]
Since \(\eta_j\neq 0\) for all \(j\), the transformed observations
$
Z_j=\eta_jY_j
$
satisfy
\[
Z_j\sim N(\eta_j\theta_j,\sigma^2\eta_j^2),
\]
so the transformed model has the same heteroscedastic Gaussian form in \eqref{eq:chhor model} if we equate $\mu_j=\eta_j\theta_j$ for all $j=1,\ldots, d$.

At this level, the structural similarity with \cite{chhor2024sparse} becomes apparent. However, the relevant functional is different. In \cite{chhor2024sparse}, the lower bound is calibrated to separate the sign-invariant functional \(\|\mu\|_t\). In our setting, the lower bound must separate the mean shift of the signed functional \(L(\theta)\). Consequently, if we use a sparse prior of the form
\[
\theta_j=\mathrm{sign}(\eta_j)b_j\gamma_j,\qquad b_j\sim \mathrm{Ber}(\pi_j),
\]
then
\[
\mathbb{E}_{\mu}[L(\theta)]
=
\sum_{j=1}^d |\eta_j|\gamma_j\pi_j.
\]

To contrast with \eqref{eq:chhor-optimization-app}, it is natural to consider the following one-sided calibration heuristic in our setting:
\begin{equation}\label{eq:our-optimization-app}
\begin{aligned}
\max_{\gamma_{1:d},\pi_{1:d}}\quad & \sum_{j=1}^d |\eta_j|\gamma_j\,\pi_j\\
\text{s.t.}\quad
& \sum_{j=1}^d \pi_j \le s/2,\\
& \pi_j\in[0,1],\quad \gamma_j\ge 0,\qquad j\in[d],\\
& \text{the induced mixture remains statistically close to the null.}
\end{aligned}
\end{equation}
The constraint on the statistical closeness is descriptive at this point and will be made accurate later.

Comparing \eqref{eq:our-optimization-app} with \eqref{eq:chhor-optimization-app} from \cite{chhor2024sparse}, the first essential difference lies in the objective: in \eqref{eq:chhor-optimization-app}, the objective is sign-invariant, whereas the objective in \eqref{eq:our-optimization-app} is sign-sensitive and depends explicitly on the loading vector \(\eta\).
As a result, the lower-bound prior in our problem must be aligned with the signs of \(\eta\), and therefore cannot be sign-symmetric.

Indeed, if one were to use the sign-symmetric prior
\[
\theta_j=b_j\omega_j\gamma_j,\qquad \omega_j\sim \mathrm{Rad}(1/2),
\]
then
\[
\mathbb{E}_{\mu}[L(\theta)]
=
\sum_{j=1}^d \eta_j\,\mathbb{E}[b_j]\mathbb{E}[\omega_j]\gamma_j
=0.
\]
Hence such a prior would generate no separation (in the prior mean) for the signed linear functional, and therefore cannot be used in our one-sided mean-separation construction.
This is already a substantive difference from the prior used in \cite{chhor2024sparse}.

\subsection{Even in the Gaussian case, the divergence constraint is different}

The difference is not limited to the target functional. Even when \(\al=2\), our \(\chi^2\)-constraint associated with our one-sided prior differs from that in \cite{chhor2024sparse}.

For one coordinate, let
\[
P_{0,j}=N(0,\sigma^2),
\qquad
Q_j=(1-\pi_j)N(0,\sigma^2)+\pi_jN(\gamma_j \sgn(\eta_j),\sigma^2).
\]
Then
\[
\frac{dQ_j}{dP_{0,j}}(x)
=
1-\pi_j+\pi_j\exp\!\left(\frac{\gamma_j \sgn(\eta_j) x}{\sigma^2}-\frac{\gamma_j^2}{2\sigma^2}\right),
\]
and a direct calculation yields
\[
1+\chi^2(Q_j\|P_{0,j})
=
1+\pi_j^2\left(\exp\!\left(\frac{\gamma_j^2}{\sigma^2}\right)-1\right).
\]
By independence across coordinates, a sufficient condition for bounded \(\chi^2\)-divergence is therefore
\[
\sum_{j=1}^d \pi_j^2\left(\exp\!\left(\frac{\gamma_j^2}{\sigma^2}\right)-1\right)\lesssim 1.
\]
Accordingly, in the Gaussian case, our one-sided calibration heuristic becomes
\begin{equation}\label{eq:our-gaussian-optimization-app}
\begin{aligned}
\max_{\gamma,\pi}\quad & \sum_{j=1}^d |\eta_j|\gamma_j\,\pi_j\\
\text{s.t.}\quad
& \sum_{j=1}^d \pi_j \le s/2,\\
& \pi_j\in[0,1],\quad \gamma_j\ge 0,\qquad j\in[d],\\
& \sum_{j=1}^d \pi_j^2\left(\exp\!\left(\frac{\gamma_j^2}{\sigma^2}\right)-1\right)\le c.
\end{aligned}
\end{equation}

This should be compared with \eqref{eq:chhor-optimization-app}, whose divergence term involves
\[
\sinh^2\!\left(\frac{\gamma_j^2}{2\sigma_j^2}\right).
\]
The two constraints are not equivalent. In particular, as \(x\to 0\),
\[
\sinh^2(x)\asymp x^2,
\qquad\text{whereas}\qquad
e^x-1\asymp x.
\]
Hence, even in the Gaussian case, the local behavior of the two divergence terms is different, which leads to different calibrations of \((\pi_j,\gamma_j)\). This is the second essential difference from \cite{chhor2024sparse}. The coincidence of the defining equation does not mean that the underlying optimization problem is the same.

\subsection{Our novelty beyond the Gaussian case}

The preceding comparison concerns the Gaussian case. Our result goes further and treats symmetric generalized Gaussian noise.
The key difference lies in the bounded $\chi^2$-divergence constraint in \eqref{eq:our-gaussian-optimization-app},
i.e.,
$$\sum_{j=1}^d \pi_j^2\left(\exp\!\left(\frac{\gamma_j^2}{\sigma^2}\right)-1\right)\lesssim 1. $$
This divergence admits a closed-form expression in the Gaussian case, but it is no longer available for symmetric generalized Gaussian noise.
Therefore, new techniques are needed.

By a direct calculation and independence across coordinates,
it is possible to see that
a sufficient condition for bounded \(\chi^2\)-divergence in our product-prior construction is
\[
\sum_{j=1}^d \pi_j^2\chi^2(P_{1,j}\|P_{0,j}) \lesssim 1.
\]
The key new ingredient is Lemma~\ref{lem: chi square}, which shows that
\[
1+\chi^2(f_\al^{(1)}\|f_\al^{(0)})
\le
C_{\al,1}\exp\!\left(
\left|\frac{\gamma}{C_{\al,2}}\right|^\al
\right),
\]
where \(C_{\al,2}\) is the constant defined in Lemma~\ref{lem: chi square}.
Applying Lemma~\ref{lem: chi square} coordinatewise yields the concrete sufficient condition for bounded \(\chi^2\)-divergence
\[
\sum_{j=1}^d
\pi_j^2
\left[
C_{\al,1}\exp\!\left(
\left|\frac{\gamma_j}{\sigma C_{\al,2}}\right|^\al
\right)-1
\right]
\lesssim 1.
\]
This is exactly where the generalized-Gaussian analysis enters the proof.

Hence, the generalized-Gaussian extension is not obtained by a straightforward modification of the Gaussian proof.
Rather, it is the new \(\chi^2\)-control in Lemma~\ref{lem: chi square} that allows us to carry out the same one-sided sparse-prior strategy sharply beyond the Gaussian case.

\subsection{Summary}

The relation with \cite{chhor2024sparse} can be summarized as follows. In the Gaussian case, our defining equation \eqref{eq:solution nonadaptive} coincides formally with Equation~(5) of \cite{chhor2024sparse}, because both arise from a similar sparse-prior lower-bound strategy.
However, our problem requires a different lower-bound construction even when \(\al=2\).
The target is a signed heterogeneous linear functional rather than a sign-invariant testing alternative; this leads to a different prior, a different divergence constraint, and a different calibration problem.
Beyond that, our extension to generalized Gaussian noise relies on the new \(\chi^2\)-control in Lemma~\ref{lem: chi square}, which constitutes another technical novelty relative to \cite{chhor2024sparse}.

\end{document}